%% file: GePUP-ES.tex
\documentclass[envcountsect,smallextended]{svjour3}       
\smartqed  
\usepackage{graphicx}
\usepackage{algorithmic}
\usepackage{amsfonts}
\usepackage{amsmath}
\usepackage{amssymb}
\usepackage{bm}
\usepackage[sans]{dsfont}
\usepackage{enumitem}
\usepackage{epstopdf}
\usepackage{graphicx}
\usepackage{hyperref}
\usepackage{mathrsfs}
\usepackage{multirow}
\usepackage{subfigure}
\usepackage{tabularx}
\usepackage{tikz}
\usepackage{url}
\usepackage{verbatim}

\journalname{}
\newcommand{\dif}{\mathrm{d}}
\newcommand{\Dim}{{\scriptsize \textup{D}}}

\newcommand{\Div}{{\mathbf{D}}}
\newcommand{\Grad}{{\mathbf{G}}}
\newcommand{\Iden}{{\mathbf{I}}}
\newcommand{\Lapl}{{\mathbf{L}}}
\newcommand{\Proj}{{\mathbf{P}}}
\newcommand{\cProj}{{\mathcal{P}}}
\newcommand{\cProjLH}{{\mathscr{P}}}

\newcommand{\innerProd}[1]{\left\langle #1 \right\rangle}

\spnewtheorem{thm}{Theorem}[section]{\bf}{\it}
\spnewtheorem{lem}[thm]{Lemma}{\bf}{\it}
\spnewtheorem{coro}[thm]{Corollary}{\bf}{\it}
\spnewtheorem{defn}[thm]{Definition}{\bf}{\rm}
\spnewtheorem{exm}[thm]{Example}{\bf}{\rm}
\spnewtheorem{rem}[thm]{Remark}{\bf}{\rm}
\spnewtheorem{ntn}[thm]{Notation}{\bf}{\it}

\begin{document}

\title{GePUP-ES:
	High-order Energy-stable Projection Methods
	for Incompressible Navier-Stokes Equations
	with No-slip Conditions
}

\titlerunning{GePUP-ES: High-order Energy-stable Projection Methods}        

\author{Yang Li
	\and
	Xu Wu
	\and
	Jiatu Yan
	\and
	Jiang Yang
	\and
	Qinghai Zhang
	\and
	Shubo Zhao
	\thanks{
		Jiatu Yan and Shubo Zhao are co-first authors with equal contributions.}
}

\authorrunning{
	Y. Li, X. Wu, J. Yan, J. Yang, Q. Zhang, and S. Zhao
} 

\institute{
	L. Yang, J. Yan \at
	School of Mathematical Sciences,
	Zhejiang University,
	Hangzhou, Zhejiang, 310058,
	China
	\and
	X. Wu and J. Yang \at
	SUSTech International Center for Mathematics,
	Southern University of Science and Technology,
	Shenzhen, Guangdong, 518055,
	China
	\and
	X. Wu \at
	Department of Applied Mathematics,
	The Hong Kong Polytechnic University,
	Kowloon, Hong Kong,
	China
	\and
	J. Yang \at
	Department of Mathematics,
	Southern University of Science and Technology,
	Shenzhen, Guangdong, 518055,
	China
	\and
	Q. Zhang and S. Zhao \at
    School of Mathematical Sciences,
	Zhejiang University,
	Hangzhou, Zhejiang, 310058,
	China;
	College of Mathematics and System Sciences,
	Xinjiang University,
	Urumqi, Xinjiang, 830046,
	China
	\and
	Q. Zhang (Corresponding author) \email{qinghai@zju.edu.cn} \at
	Institute of Fundamental and Transdisciplinary Research,
	Zhejiang University,
	Hangzhou,
	Zhejiang,
	310058, China
}

\date{}

\maketitle

\begin{abstract}
	\input{sec/abstract.tex}
\end{abstract}

\input{sec/introduction.tex}


\input{sec/preliminaries.tex}

\input{sec/analysis.tex}

\input{sec/SAV.tex}

\input{sec/algorithms.tex}

\input{sec/tests.tex}

\input{sec/conclusions.tex}

\begin{acknowledgements}
  This work was supported by the grant 12272346
  from the National Natural Science Foundation of China.
  We acknowledge helpful discussions with Prof. Buyang Li
  at the Hong Kong Polytechic University. 
  We are also grateful to one anonymous referee, 
  whose insightful comments and suggestions
  lead to an improvement of this paper.
\end{acknowledgements}

\section*{Declarations}

\subsection*{Data Availability}

The datasets generated during and/or analysed during the current study
are available from the corresponding author on reasonable request.

\subsection*{Conﬂict of interest}

The authors declare that they have no conﬂict of interest.

%
%

\bibliographystyle{spmpsci}      
\bibliography{bib/GePUP-ES.bib}

\end{document}

%% file: sec/abstract.tex
Inspired by the unconstrained PPE (UPPE) formulation
 [Liu, Liu, \& Pego 2007 Comm. Pure Appl. Math., 60 pp. 1443],
 we previously proposed the GePUP formulation
 [Zhang 2016 J. Sci. Comput., 67 pp. 1134]
 for numerically solving incompressible Navier-Stokes equations (INSE)
 on no-slip domains.
In this paper, we propose GePUP-E and GePUP-ES,
 variants of GePUP
 that feature (a) electric boundary conditions 
 with no explicit enforcement of the no-penetration condition,
 (b) equivalence to the no-slip INSE,
 (c) exponential decay of
 the divergence of an initially non-solenoidal velocity,
 and (d) monotonic decrease of the kinetic energy.
Different from UPPE,
 the GePUP-E and GePUP-ES formulations
 are of strong forms and are designed
 for finite volume/difference methods
 under the framework of method of lines.
Furthermore, 
 we develop 
 semi-discrete algorithms that preserve (c) and (d)
 and fully discrete algorithms
 that are fourth-order accurate for velocity both in time and in space.
These algorithms employ algebraically stable time integrators
 in a black-box manner
 and only consist of solving a sequence of linear equations
 in each time step.
Results of numerical tests confirm our analysis. 
 

%% file: sec/introduction.tex
\section{Introduction}
\label{sec:introduction}

The incompressible Navier-Stokes equations (INSE)
 with no-slip conditions
 govern an enormous range of real-world phenomena
 such as blood flow, turbulence,
 atmosphere and ocean currents; 
 they read 
 \begin{subequations}
   \label{eq:INS}
   \begin{align}
     \frac{\partial \mathbf{u}}{\partial t}
     + \mathbf{u}\cdot\nabla\mathbf{u}
     &= \mathbf{g} -\nabla p +  \nu\Delta \mathbf{u}
       \quad \text{ in } \Omega,
     \\
     \nabla \cdot \mathbf{u} &= 0  \qquad\qquad\qquad\quad
                               \text{ in } \Omega,
     \\
     \label{eq:no-slip}
     \mathbf{u} &= \mathbf{0} \qquad\qquad\qquad\quad
                  \text{ on } \partial\Omega,
   \end{align}
 \end{subequations}
 where $t\in [0, +\infty)$ is time,
 $\Omega$ a \emph{domain}, i.e.,
 a bounded connected open subset
 of $\mathbb{R}^{\Dim}$,
 $\partial\Omega$ the domain boundary, 
 $\mathbf{g}$ the external force,
 $p$ the pressure,
 $\mathbf{u}$  the velocity,
 and $\nu$ the kinematic viscosity.
In addition to efficiently, accurately, and faithfully
 reproducing the physical processes modeled by the equations, 
 we confront four major challenges (FMC)
 of numerically solving the INSE, 
 {
   \setlist[enumerate]{label=(FMC-\arabic{enumi})}
   \begin{enumerate}
   \item How to fulfill the solenoidal condition (\ref{eq:INS}b)
     and other physical constraints such as
     the monotonic decrease of the kinetic energy?
   \item How to ensure various types of numerical stability? 
   \item How to obtain high-order convergence both in time and in space?
   \item How to decouple time integration from spatial discretization 
     so that (i) the entire solver is constituted
     by orthogonal modules for these aspects,
     and (ii) solution methods for each aspect
     can be employed in a black-box manner and thus easily changed 
     to make the entire INSE solver versatile?
   \end{enumerate}
 }

(FMC-1) concerns
 the prominent feature of mass conservation: 
 neither source nor sink exists anywhere inside the domain.
It is well known that a violation of (\ref{eq:INS}b),
 even with small errors, 
 might lead to qualitatively different flow patterns,
 especially for large Reynolds numbers.
Another important physical constraint
 to be fulfilled in this work
 is the monotonic decrease of the total kinetic energy
 as defined in (\ref{eq:kineticEnergy}).
 
In (FMC-2), 
 a crucial and indispensable type of numerical stability
 is the eigenvalue stability for the main evolutionary variable,
 which is typically the velocity.
In addition,
 preserving the monotonicity of kinetic energy
 is equivalent to ensuring numerical stability
 on the velocity with respect to the 2-norm. 
 

Challenge (FMC-3) concerns accuracy and efficiency. 
Near no-slip boundaries,
 flows at high Reynolds numbers  
 tend to develop structures of multiple length scales and time scales.
A numerical method should resolve
 all scales that are relevant to the important physics.
Compared with fourth- and higher-order methods, 
 first- and second-order methods have 
 simpler algorithms and cheaper computations, 
 but towards a given accuracy
 the computational resources may be rapidly exhausted.
It is shown both theoretically and numerically
 in \cite[Sec. 7]{zhang16:_GePUP} that
 fourth-order methods may have a large
 efficiency advantage over second-order methods.

Challenge (FMC-3) also concerns faithfully simulating
 flows where velocity derivatives such as vorticity
 crucially affect the physics.
For first-order finite volume/difference methods, 
 the computed velocity converges,
 but the vorticity does not, nor does the velocity gradient tensor.
Consequently,
 the $O(1)$ error in $\nabla \mathbf{u}$
 may lead to structures different from that
 of the original flow.
In other words,
 it is not clear whether or not solutions of a first-order method
 have converged to the \emph{right} physics.
Similar suspicions apply to second-order methods
 for flows where second derivatives of the velocity
 are important.

Challenge (FMC-4) concerns versatility and user-friendliness
 of the numerical solver.
To cater for the problem at hand, 
 it is often desirable to change the time integrator
 from one to another. 
For example, flows with high viscosity are usually stiff
 while those with small viscosity are not; 
 accordingly,
 an implicit time integrator should be used in the former case 
 while an explicit one is usually suitable for the latter.
If the internal details of a time integrator are coupled 
 into the INSE solver in a boilerplate manner,
 it would be difficult and very inconvenient
 to change the time integrator; see also the discussion
 in the paragraph under (\ref{eq:saddleStructureSystem}). 
Hence 
 a time integrator should be treated as a black box: 
 for the ordinary differential equation (ODE)
 $\frac{\dif U}{\dif t}=f(U,t)$,
 we should only need to feed into the time integrator
 the initial condition $U^n$ and samples of $f$
 at a number of time instances 
 to get the solution $U^{n+1}$ 
 from the black box.
 
This versatility further leads to user-friendliness.
Analogous to orthogonal bases of a vector space,
 the mutually independent policies 
 span a space of solvers, 
 where each solver can be conveniently assembled
 by selecting a module for each constituting policy.
For example,
 a specific INSE solver is formed by choosing
 semi-implicit Runge-Kutta (RK) for time integration,
 finite volume for spatial discretization,
 fourth-order for accuracy and so on; 
 see Table \ref{tab:GePUP-SAV-SDIRK-policies}. 

\subsection{Previous methods related to this work}

In the original projection method
 independently proposed by
 Chorin \cite{Chorin_1968_numerical}
 and Temam \cite{temam69:_sur_navier_stokes},
 the initial condition $\mathbf{u}^n\approx \mathbf{u}(t^n)$
 is first advanced
 to an auxiliary velocity $\mathbf{u}^*$
 without worrying about the pressure gradient term
 and then $\mathbf{u}^*$ is projected
 to the divergence-free space to obtain $\mathbf{u}^{n+1}$, 
 \begin{subequations}
   \label{eq:1stOrderProjection}
   \begin{align}
     \frac{\mathbf{u}^*-\mathbf{u}^n}{k} &= -\mathbf{C}(\mathbf{u}^*, \mathbf{u}^n)+\mathbf{g}^n+\nu\Lapl\mathbf{u}^*, \\
     \mathbf{u}^{n+1} &= \Proj\mathbf{u}^*,
   \end{align}
 \end{subequations}
 where $k$ is the time step size,
 $t^n$ the starting time of the $n$th step, 
 \mbox{$\mathbf{g}^n\approx \mathbf{g}(t^n)$}, 
 $\mathbf{C}(\mathbf{u}^*, \mathbf{u}^n)\approx
 [(\mathbf{u}\cdot\nabla)\mathbf{u}](t^n)$, 
 and 
 $\Lapl$ and $\Proj$ 
 discrete approximations of the Laplacian $\Delta$
 and the Leray-Helmholtz projection $\cProjLH$, respectively; 
 see Section \ref{sec:leray-helmh-proj}.

\subsubsection{Second-order methods with fractional time stepping}
\label{sec:projectionMethods2}

The original projection method is first-order accurate
 and its improvement to the second order
 has been the aim of many subsequent works;
 see, e.g., 
  \cite{Kim.Moin_1985_application,bell89:_secon_order_projec_method_incom,orszag86:_bound_condit_incom_flows,e.03:_gauge,guermond06:_overview_of_projection_methods_for_incompressible_flows,brown01:_accurate_projection_methods_for_ins}
  and references therein.
A common basis of many second-order methods
 is the temporal discretization of (\ref{eq:INS}) with
 the trapezoidal rule,
 \begin{subequations}
   \label{eq:2ndOrderProjectionIdea}
   \begin{align}
     \frac{\mathbf{u}^{n+1}-\mathbf{u}^n}{k}+\nabla p^{n+\frac{1}{2}} &= -[(\mathbf{u}\cdot\nabla)\mathbf{u}]^{n+\frac{1}{2}}+\mathbf{g}^{n+\frac{1}{2}}+\frac{\nu}{2}\Delta(\mathbf{u}^{n+1}+\mathbf{u}^n), \\
     \nabla\cdot \mathbf{u}^{n+1} &= 0,
   \end{align}
 \end{subequations}
 where  $p^{n+\frac{1}{2}}\approx p(t^{n+\frac{1}{2}})$,
 $\mathbf{g}^{n+\frac{1}{2}} \approx \mathbf{g}(t^{n+\frac{1}{2}})$,
 and $[(\mathbf{u}\cdot\nabla)\mathbf{u}]^{n+\frac{1}{2}}\approx [(\mathbf{u}\cdot\nabla)\mathbf{u}](t^{n+\frac{1}{2}})$
 are numerical approximations
 at $t^{n+\frac{1}{2}}:=\frac{1}{2}(t^n+t^{n+1})$. 

Replacing 
 the gradient $\nabla$, the divergence $\nabla\cdot$, 
 and the Laplacian $\Delta$ in (\ref{eq:2ndOrderProjectionIdea}) 
 respectively with their second-order discrete counterparts
 $\Grad$, $\Div$, and $\Lapl$ yields
 \begin{equation}
   \label{eq:saddleStructureSystem}
   A
   \begin{bmatrix}
     \mathbf{u}^{n+1} \\ p^{n+\frac{1}{2}}
   \end{bmatrix}
   :=
   \begin{bmatrix}
     \frac{1}{k}\Iden - \frac{\nu}{2} \Lapl & \Grad
     \\
     -\Div & \mathbf{0}
   \end{bmatrix}
   \begin{bmatrix}
     \mathbf{u}^{n+1} \\ p^{n+\frac{1}{2}}
   \end{bmatrix}
   = \mathbf{F}.
 \end{equation}
Since $\Grad^T=-\Div$, the matrix $A$
 has a saddle point structure
 and the above method is often called the saddle point approach.
Despite its simplicity,
 this approach has two main disadvantages.
First, the spatial discretization and time integration
 are coupled in a boilerplate manner
 and thus a change of either part would 
 demand a complete rederivation of the matrix $A$.
For fourth- or higher-order accuracy in time integration,
 it is often too complicated to have an explicit expression
 of the matrix $A$, 
 as $A$ contains all \emph{internal} details of the time integrator. 
Consequently,
 it is highly difficult for this approach
 to address challenges (FMC-3,4).
Second, 
 it is challenging \cite{benzi05:_numer}
 to efficiently solve the linear system
 (\ref{eq:saddleStructureSystem})
 since all the velocity components
 and the pressure are coupled into a big unknown vector; 
 in contrast, Chorin's projection method
 only requires the solutions of linear systems
 with the unknowns as either the pressure
 or a velocity component.

In the fractional-stepping projection methods
\cite{Kim.Moin_1985_application,bell89:_secon_order_projec_method_incom,brown01:_accurate_projection_methods_for_ins}, 
 one replaces $\mathbf{u}^{n+1}$ and $p^{n+\frac{1}{2}}$
 in (\ref{eq:2ndOrderProjectionIdea}a)
 respectively with $\mathbf{u}^*$ and $q$, 
 solves for the auxiliary velocity $\mathbf{u}^*$
 with some boundary condition
 $\mathbf{B}(\mathbf{u}^*) = \mathbf{0}$, 
 obtains $\mathbf{u}^{n+1}$ by
 the projection $\mathbf{u}^{n+1}=\mathbf{u}^*-k\nabla \phi^{n+1}$,
 and updates the pressure
 with $p^{n+\frac{1}{2}} = q + \mathbf{U}(\phi^{n+1})$.

 Fractional-stepping projection methods
  have been very successful. 
 However, the choices of $q$, $\mathbf{B}(\mathbf{u}^*)$,
 and $\mathbf{U}(\phi^{n+1})$ 
 are coupled according to internal details of the time integrator
 \cite{brown01:_accurate_projection_methods_for_ins}.
Consequently,
 switching from one time integrator to another
 calls for a new derivation.
Furthermore,
 although appearing divorced,
 the velocity and the pressure are still implicitly coupled
 by the boundary condition of $\mathbf{u}^*$, 
 with the coupling determined not by physics
 but still by internal details of the time integrator.
Hence these methods are not suitable
 for tackling the challenges (FMC-3,4) either.

\subsubsection{The formulation via the pressure Poisson equation (PPE) }
 \label{sec:PPEmethods}
 
As a specialization of Newton's second law, 
 the momentum equation (\ref{eq:INS}a) can be rewritten as
 \begin{equation}
    \label{eq:INSEa-EulerianAccel}
    \mathbf{a}^* = \mathbf{a} + \nabla p, 
  \end{equation}
where the \emph{Eulerian acceleration}s are vector functions 
  \begin{equation}
    \label{eq:EulerianAccelerations}
    \mathbf{a} := \frac{\partial \mathbf{u}}{\partial t},\qquad
    \mathbf{a}^* := -\mathbf{u}\cdot\nabla\mathbf{u}
    +\mathbf{g} +  \nu\Delta \mathbf{u}.
  \end{equation}

The PPE describes an \emph{instantaneous} relation
 between the pressure and the velocity in the INSE
 and, on no-slip domains, has the form
 \begin{subequations}
   \label{eq:PPE}
   \begin{align}
     \Delta p &= \nabla \cdot
                (\mathbf{g}-\mathbf{u}\cdot\nabla\mathbf{u}
                + \nu \Delta \mathbf{u})
                \qquad \text{ in } \Omega,
     \\
     \mathbf{n}\cdot \nabla p
              &= \mathbf{n}\cdot (\mathbf{g}+\nu\Delta\mathbf{u})
                \quad\qquad\qquad \quad\text{ on } \partial \Omega, 
   \end{align}
 \end{subequations}
 where (\ref{eq:PPE}b) follows from
 the normal component of (\ref{eq:INS}a) and the no-slip conditions
 (\ref{eq:no-slip})
 while (\ref{eq:PPE}a) from the divergence of (\ref{eq:INS}a) 
 and the divergence-free condition (\ref{eq:INS}b).
For the PPE with other boundary conditions,
 (\ref{eq:PPE}b) should be replaced with
 the normal component of (\ref{eq:INSEa-EulerianAccel}). 
As explained in Section \ref{sec:leray-helmh-proj},
 the pressure gradient is uniquely
 determined from $\mathbf{a}^*$ by
 $\nabla p = ({\mathcal I}-\cProjLH) \mathbf{a}^*$
 where ${\mathcal I}$ and $\cProjLH$
 are the identity operator
 and the Leray-Helmholtz projection, respectively.
Thus neither the initial condition nor the boundary condition
 of the pressure $p$ is needed in the INSE.

(\ref{eq:INS}a), (\ref{eq:no-slip}),
 (\ref{eq:PPE}), and the additional boundary condition
 $\nabla\cdot\mathbf{u}=0$ on $\partial \Omega$
 are collectively called
 \emph{the PPE formulation of the INSE on no-slip domains}
 \cite{gresho87:_navier_stokes}.
In terms of computation, however, 
 the PPE formulation 
 has a decisive advantage over the original INSE.
If (\ref{eq:INS}a) is discretized in time
 with (\ref{eq:INS}b) as a constraint,
 the resulting index-2 differential algebraic system
 may suffer from large order reductions
 \cite{sanderse12:_accur_runge_kutta_navier_stokes}.
In contrast,
 replacing the divergence-free constraint
 with the PPE avoids this difficulty.
As such,
 the PPE formulation allows
 the time integrator to be treated as a black box
 and thus to be easily changed; 
 indeed, the pressure is an implicit function of $\mathbf{u}$
 and its interaction with $\mathbf{u}$
 is completely decoupled from
 internal details of the time integrator.
Also,
 there is no need to introduce
 nonphysical auxiliary variables. 
These advantages of the PPE formulation lead to
 quite a number of successful numerical methods \cite{kleiser80:_treat_incom_bound_condit_d,gresho87:_navier_stokes,henshaw94:_fourt_accur_method_incom_navier,johnston04:_accur_navier_stokes,liu07:_stabil_conver_effic_navier_stokes,liu10_stable_accurate_pressure_unsteady_incompressible_viscous_flow,shirokoff11:_navier_stokes}.

Unfortunately,
 as observed by Liu, Liu \& Pego
 \cite{liu07:_stabil_conver_effic_navier_stokes}, 
 (\ref{eq:INS}a) and (\ref{eq:PPE}a)
 yield
 \begin{equation}
   \label{eq:degenerateVelDivPPE}
   \frac{\partial \nabla\cdot\mathbf{u}}{\partial t} = 0;
 \end{equation}
 this degenerate equation implies that
 in the PPE formulation 
 we have no control over $\nabla\cdot\mathbf{u}$
 and its evolution is up to the particularities
 of the numerical schemes.
Our tests show that a fourth-order finite-volume
 method-of-lines (MOL) discretization of the PPE formulation
 is unstable, 
 with the computed velocity divergence growing indefinitely
 near the domain boundary. 

\subsubsection{The formulation of unconstrained PPE (UPPE) }
\label{sec:UPPE}

The application of the Leray-Helmholtz projection $\cProjLH$
 to (\ref{eq:INSEa-EulerianAccel}) yields
 \begin{equation}
   \label{eq:UPPE1}
   \frac{\partial \mathbf{u}}{\partial t}
   - \cProjLH \mathbf{a}^*
   = \nu\nabla(\nabla\cdot \mathbf{u}),
 \end{equation}
 where the zero right-hand side (RHS) is 
 added for stability reasons
 \cite{liu07:_stabil_conver_effic_navier_stokes}.
The divergence of (\ref{eq:UPPE1})
 and the second identity in (\ref{eq:ProjectionProperties})
 give
 \begin{equation}
   \label{eq:heatEqInUPPE}
   \frac{\partial(\nabla\cdot\mathbf{u})}{\partial t}
   = \nu\Delta(\nabla\cdot\mathbf{u}), 
 \end{equation}
 which, by the maximum principle of the heat equation, 
 dictates an exponential decay of 
 $\nabla\cdot \mathbf{u}$ in $\Omega$.
A juxtaposition of (\ref{eq:heatEqInUPPE})
 and (\ref{eq:degenerateVelDivPPE}) 
 exposes a prominent advantage of 
 (\ref{eq:UPPE1}) 
 that any divergence residue is now well under control.

Via the identity $\nabla(\nabla\cdot \mathbf{u})
 =\Delta({\mathcal I}-\cProjLH)\mathbf{u}$
 and the Laplace-Leray commutator
 $[\Delta, \cProjLH] := \Delta \cProjLH - \cProjLH \Delta$, 
 cf. Section \ref{sec:lapl-leray-comm}, 
 Liu, Liu \& Pego \cite{liu07:_stabil_conver_effic_navier_stokes}
 rewrote (\ref{eq:UPPE1}) as 
\begin{equation}
  \label{eq:UPPE2}
  \frac{\partial \mathbf{u}}{\partial t}
  + \cProjLH (\mathbf{u}\cdot\nabla \mathbf{u} -\mathbf{g})
  + \nu [\Delta, \cProjLH]\mathbf{u}
  = \nu\Delta \mathbf{u}, 
\end{equation}
 which provides a fresh viewpoint
 of the INSE
 as a controlled perturbation of the vector diffusion equation
 $\frac{\partial \mathbf{u}}{\partial t}
  = \nu\Delta \mathbf{u}$.
For $\mathbf{u}\in H^2 \cap H_0^1(\Omega, \mathbb{R}^{\Dim})$
 with ${\mathcal C}^3$ boundary $\partial \Omega$,
 they gave a sharp bound on $\|[\Delta, \cProjLH]\mathbf{u}\|$
 in terms of $\|\Delta \mathbf{u}\|$
 and 
 proved the unconditional stability and convergence
 of a first-order scheme,
 \begin{subequations}
   \label{eq:UPPE-1stOrder}
   \begin{align}
     \innerProd{\nabla p^n, \nabla \phi}
     &= \innerProd{\mathbf{g}^{n}- \mathbf{u}^n\cdot\nabla \mathbf{u}^{n}
       +\nu\Delta \mathbf{u}^{n} - \nu\nabla\nabla\cdot\mathbf{u}^{n},
       \nabla \phi}, 
     \\
     \frac{\mathbf{u}^{n+1}-\mathbf{u}^n}{k}+\nabla p^{n}
     &= \mathbf{g}^{n}- \mathbf{u}^n\cdot\nabla \mathbf{u}^{n}
       +\nu\Delta\mathbf{u}^{n+1} \qquad \text{ in } \Omega,\\
     \mathbf{u}^{n+1} &= \mathbf{0} \qquad\qquad\qquad\qquad\qquad\qquad
                        \text{ on } \partial\Omega, 
   \end{align}
 \end{subequations}
 where (\ref{eq:UPPE-1stOrder}a) is the PPE in weak form with
 $\phi \in H^1(\Omega)$
 and  $\left\langle\mathbf{u}, \mathbf{v}\right\rangle  :=
 \int_{\Omega} \mathbf{u}\cdot \mathbf{v}\,\dif V$. 

Based on the UPPE formulation,
 a slip-corrected projection method
 \cite{liu10_stable_accurate_pressure_unsteady_incompressible_viscous_flow}
 is developed with third-order accuracy
 both in time and in space.
 
%
From (\ref{eq:UPPE2}) and the contents
 in Sections \ref{sec:lapl-leray-comm} and \ref{sec:stokes-pressure}, 
 a strong form of UPPE can be deduced
\cite{liu10_stable_accurate_pressure_unsteady_incompressible_viscous_flow} as
  \begin{subequations}
    \label{eq:UPPEstrong}
    \begin{align}
      \frac{\partial \mathbf{u}}{\partial t}
      + \mathbf{u}\cdot\nabla\mathbf{u}
      &= \mathbf{g} -\nabla p +  \nu\Delta \mathbf{u}
        \qquad\qquad\qquad \text{ in } \Omega,
      \\
      \mathbf{u} &= \mathbf{0} \qquad\qquad\qquad\qquad\qquad\qquad
                   \text{ on } \partial\Omega,
      \\
      \Delta p
      &= \nabla\cdot(\mathbf{g}-\mathbf{u}\cdot\nabla\mathbf{u})
        \qquad\qquad\quad\ \  \text{ in } \Omega, \\
     \mathbf{n}\cdot\nabla p &=
                               \mathbf{n}\cdot(\mathbf{g}+\nu\Delta\mathbf{u}-\nu\nabla\nabla\cdot\mathbf{u})\qquad
                               \text{ on } \partial\Omega. 
    \end{align}
  \end{subequations}
  
The PPE (\ref{eq:PPE})
 and the UPPE (\ref{eq:UPPEstrong}c,d)
 have slightly different forms
 and nonetheless a crucial distinction:
 it follows from the divergence of (\ref{eq:UPPEstrong}a)
 that (\ref{eq:UPPEstrong}c) leads to (\ref{eq:heatEqInUPPE})
 whereas 
 (\ref{eq:PPE}a) leads to (\ref{eq:degenerateVelDivPPE}). 
  
Unfortunately,
 (\ref{eq:UPPEstrong}) is not yet suitable
 for the design of MOL-type finite volume
 and finite difference methods, 
 due to two main reasons (TMR).
 {\setlist[enumerate]{label=(TMR-\arabic{enumi})}
 \begin{enumerate}
 \item The Leray-Helmholtz projection is absent in (\ref{eq:UPPEstrong})
   and thus any projection on the velocity in an MOL algorithm
   would be a mismatch of the numerical algorithm
   to the governing equations.
   Of course one can replace the velocity $\mathbf{u}$
   with $\cProjLH \mathbf{u}$ in (\ref{eq:UPPEstrong}),
   but which should be replaced?
   In other words, which $\mathbf{u}$'s in (\ref{eq:UPPEstrong})
   should be projected in MOL?
 \item It is difficult for a discrete projection $\mathbf{P}$
   with fourth- and higher-order accuracy
   to satisfy all properties of the Leray-Helmholtz projection $\cProjLH$
   in (\ref{eq:ProjectionProperties}).
   In particular, the discretely projected velocity
   may not be divergence-free.
   Then how does the approximation error of
   $\mathbf{P}$ to $\cProjLH$
   affect the stability of the ODE system
   under the MOL framework?
   It is neither clear nor trivial
   how to answer this question with (\ref{eq:UPPEstrong}). 
 \end{enumerate}
}

\subsubsection{Generic projection
  and unconstrained PPE (GePUP)}
\label{sec:GePUP-intro}

A \emph{generic projection}
 is a linear operator $\cProj$ on a vector space
 satisfying 
 \begin{equation}
   \label{eq:cProj}
   \cProj \mathbf{u} = \mathbf{w} := \mathbf{u} -\nabla\phi,
 \end{equation}
 where $\phi$ is a scalar function
 and $\nabla\cdot\mathbf{w}=0$ may or may not hold.
Since $\phi$ is not specified in terms of $\mathbf{w}$,
 (\ref{eq:cProj}) is not a precise definition of $\cProj$, 
 but rather a characterization of a family of operators, 
 which, in particular, includes
 the Leray-Helmholtz projection $\cProjLH$. 
 $\cProj$ can be used to perturb $\mathbf{u}$
 to some non-solenoidal velocity $\mathbf{w}$
 and is thus more flexible than $\cProjLH$
 in characterizing discrete projections
 that fail to fulfill the divergence-free constraint exactly.

To accommodate the fact
 that the discrete velocity might not be divergence-free,
 we switch the evolutionary variable
 to a non-solenoidal velocity $\mathbf{w}=\cProj \mathbf{u}$
 instead of the divergence-free velocity $\mathbf{u}$
 in the UPPE formulation (\ref{eq:UPPEstrong}).
More precisely,
 the evolutionary variable $\mathbf{u}$
 in the time-derivative term $\frac{\partial \mathbf{u}}{\partial t}$
 is perturbed to a non-solenoidal velocity $\mathbf{w}:=\mathbf{u}-\nabla \phi$
 where $\phi$ is some scalar function; 
 meanwhile in the diffusion term we change $\mathbf{u}$ to $\mathbf{w}$ 
 to set up a mechanism that drives the divergence towards zero.
Then, there is no need to worry about
 the influence of $\nabla\cdot\mathbf{w}\ne 0$ on numerical stability
 because the evolution of $\mathbf{w}$
 is not subject to the divergence-free constraint.
These ideas lead to the \emph{GePUP formulation} \cite{zhang16:_GePUP}: 
\begin{subequations}
  \label{eq:GePUP}
  \begin{alignat}{2}
    \frac{\partial \mathbf{w}}{\partial t} &= \mathbf{g}-\mathbf{u}\cdot \nabla \mathbf{u}-\nabla q+\nu\Delta \mathbf{w} &\quad& \text{in } \Omega, \\
    \mathbf{w} &= \mathbf{0},\ \
    \mathbf{u}\cdot \boldsymbol{\tau} = 0, \ \ 
    && \text{on } \partial \Omega, \\
    \mathbf{u} &= \cProjLH\mathbf{w} && \text{in } \Omega, \\
    \mathbf{u}\cdot \mathbf{n} &= 0 && \text{on } \partial\Omega, \\
    \Delta q &= \nabla\cdot(\mathbf{g}-\mathbf{u}\cdot\nabla \mathbf{u}) && \text{in } \Omega, \\
    \mathbf{n}\cdot \nabla q &= \mathbf{n}\cdot(\mathbf{g}+\nu\Delta
    \mathbf{u}-\nu\nabla\nabla\cdot \mathbf{u}) && \text{on } \partial
    \Omega, 
  \end{alignat}
\end{subequations}
where (\ref{eq:GePUP}e) and the divergence of (\ref{eq:GePUP}) yield
 $\frac{\partial\left(\nabla\cdot \mathbf{w}\right)}{\partial t}
 = \nu\Delta\left(\nabla\cdot \mathbf{w}\right)$.
Then either $\nabla\cdot\mathbf{w}=0$
 or $\mathbf{n}\cdot \nabla\nabla\cdot \mathbf{w}=0$
 on $\partial \Omega$ 
 drives $\nabla \cdot \mathbf{w}$ towards zero. 
 
 
 
\subsection{The contribution of this work}

We couple GePUP with electric boundary conditions
 \cite{shirokoff11:_navier_stokes,rosales21:_high_poiss_navier_stokes} 
 and a scalar auxiliary variable (SAV)
 \cite{shen18:SAV,shen19:SAV_Review}
 to propose GePUP-E and GePUP-ES,
 variants of GePUP that enforce the solenoidal conditions,
 preserve energy stability,
 decouple time integration from spatial discretization,
 and lead to versatile algorithms
 that are fourth-order accurate both in time and in space.
The letter `E' in the acronyms stands for the electric boundary
 conditions
 while the letter `S' for the SAV approach; 
 altogether `ES' also stands for energy stability.

GePUP-E and GePUP-ES answer all the challenges in (FMC-1,2,3,4).
 {\setlist[enumerate]{label=(\Alph{enumi})}
\begin{enumerate}
\item We reformulate the INSE of two variables
   into the GePUP-E formulation (\ref{eq:GePUPe})
   of three variables,
   embedding the solution manifold $\mathcal{M}(\mathbf{u},p)$ of the INSE
   in the higher-dimensional solution manifold
   $\mathcal{N}(\mathbf{w},\mathbf{u},q)$ of GePUP-E.
  More importantly, we equip this embedding
   with a divergence-decaying mechanism
   that drives any deviating solution in
   $\mathcal{N}(\mathbf{w},\mathbf{u},q)$ 
   back to $\mathcal{M}(\mathbf{u},p)$.
  We prove the equivalence of INSE and GePUP-E,
   the convergence of the non-solenoidal velocity $\mathbf{w}$ to
   the divergence-free velocity $\mathbf{u}$,
   the exponential decay of the divergence $\nabla\cdot\mathbf{w}$,
   and the monotonic decrease of the kinetic energy.
  As such, GePUP-E resolves the difficulties in (TMR-1,2). 
\item By coupling GePUP-E to SAV
  \cite{shen18:SAV,shen19:SAV_Review}, 
  we propose the GePUP-ES formulation in (\ref{eq:GePUPSAV}), 
  design a family of semi-discrete GePUP-ES algorithms, 
  and prove their energy stability in Theorem \ref{thm:GePUPSAVRKEnergyDecay}.
\item Based on (B), 
  we further propose a family of fully discrete INSE solvers, 
  named GePUP-ES-SDIRK, 
  to answer all challenges in (FMC-1,2,3,4).
\end{enumerate}
}
(A), (B), and (C) are elaborated in Sections
 \ref{sec:gepupe-sav}, \ref{sec:GePUP-ES},
 and \ref{sec:algorithms}, respectively.
In Section \ref{sec:preliminaries}, 
 we introduce notation
 to make this paper somewhat self-contained.
We test GePUP-ES-SDIRK in Section \ref{sec:tests}
 and draw conclusions in Section \ref{sec:conclusions}. 



%% file: sec/preliminaries.tex
\section{Preliminaries}
\label{sec:preliminaries}

Throughout this paper, 
we denote by $\left\langle\cdot, \cdot\right\rangle$
 the $L^2$ inner product
 of vector- (or scalar-) valued functions $\mathbf{u}$ and $\mathbf{v}$
 over $\Omega$, 
 $\left\langle\mathbf{u}, \mathbf{v}\right\rangle  :=
 \int_{\Omega} \mathbf{u}\cdot \mathbf{v}\,\dif V$,
 and by $\|\cdot\|$ the induced $L^2$ norm
 $\|\mathbf{u}\| :=
 \sqrt{\left\langle \mathbf{u},\mathbf{u}\right\rangle}$.
  
We start with a well-known result on boundary value problems (BVPs).

\begin{thm}[Solvability of BVPs with pure Neumann conditions]
  \label{thm:ExistenceAndUniquenessOfPoissonNeumannBVP}
  Suppose $f$ and $g$ are two sufficiently smooth functions.
  Then there exists a unique solution
  (up to an additive constant) for the Neumann BVP
  \begin{subequations}
    \label{eq:PoissonNeumannBVP}
    \begin{align}
      \Delta \phi &= f\qquad \text{ in } \Omega; \\
      \mathbf{n}\cdot \nabla \phi &= g\qquad \text{ on } \partial \Omega
    \end{align}
  \end{subequations}
  if and only if
  $\int_{\Omega}f\,\dif V = \int_{\partial \Omega}g\,\dif A$.
\end{thm}
\begin{proof}
  See \cite[page 409]{taylor11:_partial_differ_equat_i}.
\qed 
\end{proof}

\subsection{The Leray-Helmholtz projection $\cProjLH$}
\label{sec:leray-helmh-proj}

\begin{thm}[Helmholtz decomposition]
  \label{thm:HelmholtzHodgeDecomposition}
  A continuously differentiable 
  vector field $\mathbf{v}^*$ in a domain $\Omega$ can be
  uniquely decomposed into a divergence-free part $\mathbf{v}$
  and a curl-free part $\nabla \phi$, 
  \begin{equation}
    \label{eq:HelmholtzHodgeDecomposition}
    \mathbf{v}^* = \mathbf{v}+\nabla \phi,
  \end{equation}
  where $\mathbf{v}\cdot \mathbf{n}$
  is given \emph{a priori} on $\partial \Omega$
  and satisfies $\oint_{\partial\Omega}\mathbf{v}\cdot\mathbf{n}=0$.
\end{thm}
\begin{proof}
  The decomposition can be realized
  by solving 
  \begin{equation}
    \label{eq:NeumannBVPforHelmDecomp}
    \left\{
    \begin{array}{rll}
      \Delta \phi &= \nabla \cdot \mathbf{v}^* & \text{ in } \Omega, 
      \\
      \mathbf{n}\cdot \nabla \phi &= \mathbf{n}\cdot
                                    (\mathbf{v}^*-\mathbf{v})& 
                                    \text{ on } \partial \Omega, 
    \end{array}
    \right.
  \end{equation}
  since Theorem \ref{thm:ExistenceAndUniquenessOfPoissonNeumannBVP}
  uniquely determines $\nabla \phi$
  with $\oint_{\partial\Omega}\mathbf{v}\cdot\mathbf{n}=0$. 
\qed 
\end{proof}


The \emph{Leray-Helmholtz projection} $\cProjLH$ 
 is an idempotent operator
that maps a vector field $\mathbf{v}^*$
to its divergence-free part $\mathbf{v}$, 
c.f. the decomposition (\ref{eq:HelmholtzHodgeDecomposition}),
i.e., 
\begin{equation}
  \label{eq:cProjLH}
  \cProjLH \mathbf{v}^* := \mathbf{v} = \mathbf{v}^* -\nabla\phi.
\end{equation}

The proof of Theorem \ref{thm:HelmholtzHodgeDecomposition} implies
 the constructive form
 $\cProjLH = {\mathcal I} - \nabla (\Delta_n)^{-1} \nabla\cdot$, 
 where $(\Delta_n)^{-1}$ denotes solving (\ref{eq:NeumannBVPforHelmDecomp}).
 For a ${\mathcal C}^1$ vector field $\mathbf{v}^*$
 and a ${\mathcal C}^1$ scalar field $\phi$, we have
\begin{equation}
  \label{eq:ProjectionProperties}
  \cProjLH^2=\cProjLH,\quad
  \nabla\cdot\cProjLH\mathbf{v}^* =0,\quad
  \cProjLH\nabla\phi = \mathbf{0}.
\end{equation}

\subsection{The Laplace-Leray commutator
  $\Delta \cProjLH - \cProjLH \Delta$}
\label{sec:lapl-leray-comm}

On periodic domains, $\Delta$ and $\cProjLH$ commute.
However, one main difficulty for no-slip domains is the fact that 
$\Delta \cProjLH - \cProjLH \Delta \ne \mathbf{0}$. 
In this subsection we rephrase several results
 in \cite{liu07:_stabil_conver_effic_navier_stokes}.

\begin{lem}
  \label{lem:propertiesOfB}
  The \emph{divergence-gradient commutator} defined as
  \begin{equation}
    \label{eq:opB}
    {\mathcal B} = [\nabla\cdot, \nabla] :=\Delta - \nabla \nabla\cdot 
  \end{equation}  
  satisfies $\nabla\cdot{\mathcal B} = 0$, 
  $\Delta \cProjLH = {\mathcal B}$, and
  in three dimensions ${\mathcal B} = -\nabla \times \nabla \times$. 
\end{lem}
\begin{proof}
  $\nabla\cdot{\mathcal B} = 0$ follows from (\ref{eq:opB})
  and $\Delta \nabla\cdot = \nabla\cdot \Delta$
  while ${\mathcal B}\mathbf{v}^*
  = -\nabla \times \nabla \times \mathbf{v}^*$
  from the tensor notation and the epsilon-delta relation. 
  $\Delta \cProjLH = {\mathcal B}$ holds because
  \begin{equation}
    \label{eq:deriveCommutatorB}
    \Delta \cProjLH \mathbf{v}^*
    = \Delta(\mathbf{v}^* - \nabla\phi)
    = \Delta \mathbf{v}^* - \nabla\Delta\phi
    = \Delta \mathbf{v}^* - \nabla \nabla\cdot \mathbf{v}^*,
  \end{equation}
  where we have applied (\ref{eq:cProjLH}),
  the commutativity of $\Delta$ and $\nabla$,
  and (\ref{eq:opB}).
\qed 
\end{proof}

With $\cProjLH\nabla\phi = \mathbf{0}$ in
(\ref{eq:ProjectionProperties}), 
the first and third terms
in (\ref{eq:deriveCommutatorB})
lead to 
\begin{equation}
  \label{eq:LerayLaplacian}
  \mathscr{P}\Delta\mathscr{P} = \mathscr{P}\Delta
\end{equation}
because $\mathscr{P}\Delta\mathscr{P}\mathbf{v}^* = \mathscr{P}\Delta
\mathbf{v}^*$ holds for any sufficiently smooth vector field $\mathbf{v}^*$.

Then Lemma \ref{lem:propertiesOfB} and (\ref{eq:LerayLaplacian}) give

\begin{coro}
  \label{coro:commutatorLL}
  The \emph{Laplace-Leray commutator} is 
  \begin{equation}
    \label{eq:commutatorLL}
    [\Delta, \cProjLH] := \Delta \cProjLH - \cProjLH \Delta
    = ({\mathcal I}-\cProjLH)\Delta \cProjLH
    = ({\mathcal I}-\cProjLH) {\mathcal B}
    = -({\mathcal I}-\cProjLH)(\nabla \times \nabla \times), 
  \end{equation}
  where ${\mathcal I}$ is the identity operator
  and the last equality holds only in three dimensions.
\end{coro}

\subsection{The Stokes pressure}
\label{sec:stokes-pressure}

By (\ref{eq:commutatorLL}),
the action of the Laplace-Leray commutator on
any vector field $\mathbf{v}^*$
results in the gradient of some scalar field.
In the case of $\mathbf{v}^*$
being the velocity $\mathbf{u}$ in the INSE, 
the scalar is known as the \emph{Stokes pressure} \cite{liu07:_stabil_conver_effic_navier_stokes}:
\begin{equation}
  \label{eq:stokesPressure}
  \nabla p_s := (\Delta \cProjLH - \cProjLH \Delta)
  \mathbf{u}.
\end{equation}
It follows from (\ref{eq:commutatorLL}) and (\ref{eq:stokesPressure})
that
\begin{displaymath}
  {\mathcal B}\mathbf{u}
  = \cProjLH {\mathcal B}\mathbf{u} + \nabla p_s.
\end{displaymath}
Then (\ref{eq:ProjectionProperties}) and $\nabla\cdot{\mathcal B} = 0$
in Lemma \ref{lem:propertiesOfB}
yield $\Delta p_s=0$,
i.e., the Stokes pressure is harmonic.
Interestingly, the vector field $\nabla p_s$
is both divergence-free and curl-free.

Define another scalar $p_c$ as 
\begin{equation}
  \label{eq:pc}
  \nabla p_c := 
  ({\mathcal I} - \cProjLH)\left(\mathbf{g}
    -\mathbf{u}\cdot\nabla\mathbf{u}\right).
\end{equation}
Apply the Leray-Helmholtz projection to (\ref{eq:INS}a),
use the commutator (\ref{eq:commutatorLL}),
invoke the definitions
(\ref{eq:stokesPressure}) and (\ref{eq:pc}),
and we have
\begin{displaymath}
  \nabla p = \nabla p_c + \nu \nabla p_s.
\end{displaymath}

The pressure gradient in the INSE consists of two parts:
$\nabla p_c$ balances the divergence of the forcing term
and the nonlinear convection term
while $\nabla p_s$ accounts for the Laplace-Leray commutator.
In the two limiting cases of $\nu\rightarrow 0$
and $\nu\rightarrow +\infty$,
the pressure gradient is 
dominated by $\nabla p_c$ and $\nu \nabla p_s$, respectively. 

\subsection{Vector identities}
\label{sec:some-vect-ident}

See a standard text on differentiable manifolds
 such as \cite{munkres97:_analy_manif} for a proof of

\begin{thm}[Gauss-Green]
  \label{thm:GaussGreen}
  A scalar or vector function
  $u\in \mathcal{C}^1(\overline{\Omega})$ satisfies
  \begin{displaymath}
    \int_{\Omega} \frac{\partial u}{\partial x_i}\dif V = \int_{\partial \Omega}un_i\,\dif A,
  \end{displaymath}
  where $n_i$ is the $i$th component of the unit outward normal $\mathbf{n}$ of $\partial \Omega$.
\end{thm}

Apply Theorem \ref{thm:GaussGreen} to $uv$ and we have

\begin{lem}[Integration-by-parts]
  \label{lem:integrationByParts}
  For $u,v\in \mathcal{C}^1(\overline{\Omega})$, we have
  \begin{equation}
    \label{eq:integrationByParts}
    \int_{\Omega} \frac{\partial u}{\partial x_i}v\,\dif V = -\int_{\Omega}u \frac{\partial v}{\partial x_i}\dif V + \int_{\partial \Omega}uvn_i\,\dif A,
  \end{equation}
  where $n_i$ is the $i$th component of the unit outward normal $\mathbf{n}$ of $\partial \Omega$.
\end{lem}

Replace $v$ in (\ref{eq:integrationByParts})
 with $\frac{\partial v}{\partial x_i}$, 
sum over $i$,
and we have

\begin{lem}[Green's formula]
  For $u,v\in \mathcal{C}^2(\overline{\Omega})$, 
  we have
  \begin{equation}
    \label{eq:GreenFormula}
    \int_{\Omega}\nabla u\cdot \nabla v\,\dif V = -\int_{\Omega}u\Delta v\,\dif V + \int_{\partial \Omega}u \frac{\partial v}{\partial \mathbf{n}}\dif A.
  \end{equation}
\end{lem}

\subsection{B-stable and algebraically stable RK methods}
\label{sec:algebraicStabilityRK}

To solve an ODE system $\mathbf{u}' = \mathbf{f}(\mathbf{u},t)$, 
an \emph{$s$-stage RK method} is a 
one-step method of the form
\begin{equation}
  \label{eq:RungeKutta}
  \renewcommand{\arraystretch}{1.2}
  \left\{
    \begin{array}{l}
      \mathbf{y}_i = \mathbf{f}(\mathbf{U}^n+ k\sum_{j=1}^s a_{i,j}\mathbf{y}_j, t_n+ c_ik),
      \\
      \mathbf{U}^{n+1} = \mathbf{U}^n + k \sum_{j=1}^sb_j\mathbf{y}_j, 
    \end{array}
  \right.
\end{equation}
where $i=1, 2, \ldots, s$
and the coefficients $a_{i,j}$, $b_j$, $c_i$ are real.

A function $\mathbf{f}: \mathbb{R}^{\Dim}\times[0, +\infty)\to \mathbb{R}^{\Dim}$
is \emph{one-sided Lipschitz continuous} if
\begin{equation}
  \label{eq:oneSidedLipschitz}
  \forall t\ge 0, \forall \mathbf{u},\mathbf{v}\in \mathbb{R}^{\Dim},\ \ 
  \left\langle \mathbf{u}-\mathbf{v},\ 
    \mathbf{f}(\mathbf{u},t)-\mathbf{f}(\mathbf{v},t)\right\rangle
  \le \mu\|\mathbf{u}-\mathbf{v}\|^2,
\end{equation}
where $\mu$ is the \emph{one-sided Lipschitz constant} of $\mathbf{f}$.
The ODE system 
is \emph{contractive} or \emph{monotone}
if $\mathbf{f}$ satisfies (\ref{eq:oneSidedLipschitz})
with $\mu=0$.

A contractive ODE system is \emph{dissipative}:
for any solutions $\mathbf{u}(t)$ and $\mathbf{v}(t)$,
the norm $\|\mathbf{u}(t)-\mathbf{v}(t)\|$
decreases monotonically as $t$ increases.
In other words,
different solution trajectories of a contractive ODE system
never depart from each other, 
and hence small perturbations remain small.
This leads to

\begin{defn}[B-stability \cite{butcher75:_runge_kutta}]
  \label{def:B-Stability}
  A one-step method 
  is \emph{B-stable} 
  if, for any contractive ODE system, 
  each pair of numerical solutions
  $\mathbf{U}^n$ and $\mathbf{V}^n$ satisfy
  \begin{displaymath}
    \forall n=0,1,\ldots,\quad
    \|\mathbf{U}^{n+1}-\mathbf{V}^{n+1}\| \le \|\mathbf{U}^n-\mathbf{V}^n\|.
  \end{displaymath}
\end{defn}

It can be shown that B-stable methods are A-stable.

\begin{defn}
  \label{def:RKAlgebraicStability}
  An RK method is \emph{algebraically stable} if
  \begin{itemize}
  \item the RK weights $b_1, b_2, \ldots, b_s$ are nonnegative, 
  \item the following symmetric matrix $M\in \mathbb{R}^{s\times s}$ 
    is positive semidefinite:
    \begin{equation}
      \label{eq:algebraicStabilityMatrix}
      m_{i,j} = b_ia_{i,j}+b_ja_{j,i}-b_ib_j. 
    \end{equation}
  \end{itemize}
\end{defn}

An algebraically stable RK method is B-stable
  and thus A-stable \cite{hairer96:_solvin_ordin_differ_equat_ii}.



%% file: sec/analysis.tex
\section{The GePUP-E formulation}
\label{sec:gepupe-sav}

Boundary conditions of the evolutionary variable $\mathbf{w}$
 play a crucial role in
 establishing a stable numerical scheme.
As mentioned in Subsection \ref{sec:GePUP-intro},
 the GePUP formulation (\ref{eq:GePUP}) leads to the heat equation
 $\frac{\partial\left(\nabla\cdot \mathbf{w}\right)}{\partial t}
 = \nu\Delta\left(\nabla\cdot \mathbf{w}\right)$.
However, 
 neither the homogeneous Dirichlet
 nor the homogeneous Neumann condition
 is explicitly imposed on $\nabla\cdot\mathbf{w}$ in (\ref{eq:GePUP}),
 thus it is difficult to prove
 the dacay of $\nabla\cdot\mathbf{w}$.
To fix this glitch,
 we draw inspiration from the excellent work
 of Rosales, Shirokoff, and their colleagues
 \cite{shirokoff11:_navier_stokes,rosales21:_high_poiss_navier_stokes}
 to adapt ``electric'' boundary conditions into GePUP, 
 proposing 

\begin{defn}
  \label{def:GePUPe-equations}
  The \emph{GePUP-E formulation of INSE on no-slip domains} is
  \begin{subequations}
    \label{eq:GePUPe}
    \begin{alignat}{2}
      \frac{\partial \mathbf{w}}{\partial t} &= \mathbf{g}-\mathbf{u}\cdot \nabla \mathbf{u}-\nabla q+\nu\Delta \mathbf{w} &\quad& \text{in } \Omega, \\
      \mathbf{w}\cdot \boldsymbol{\tau} &= 0, \ \
      \nabla \cdot \mathbf{w} = 0 \ \
      && \text{on } \partial \Omega, \\
      \mathbf{u} &= \cProjLH\mathbf{w} && \text{in } \Omega, \\
      \mathbf{u}\cdot \mathbf{n} &= 0 && \text{on } \partial\Omega, \\
      \Delta q &= \nabla\cdot(\mathbf{g}-\mathbf{u}\cdot\nabla \mathbf{u}) && \text{in } \Omega, \\
      \mathbf{n}\cdot \nabla q &= \mathbf{n}\cdot(\mathbf{g}
       - \mathbf{u}\cdot\nabla\mathbf{u}
      +\nu\Delta\mathbf{w})
      + \lambda\mathbf{n}\cdot\mathbf{w} && \text{on } \partial \Omega,
    \end{alignat}
  \end{subequations}
  where 
  $\mathbf{u}$ is the divergence-free velocity in (\ref{eq:INS}),
  $\mathbf{w}=\mathbf{u}-\nabla \phi$ a non-solenoidal velocity
  for some scalar function $\phi$, 
  $\mathbf{n}$ and $\boldsymbol{\tau}$
  the unit normal and unit tangent vector
  of $\partial \Omega$, respectively,
  and $\lambda$ a nonnegative penalty parameter.
  The two velocities $\mathbf{w}$ and $\mathbf{u}$ have the same initial condition
  in $\overline{\Omega}$, the closure of $\Omega$, i.e.,
  \begin{equation}
    \label{eq:GePUP-initialW}
    \forall \mathbf{x}\in \overline{\Omega},\quad
    \mathbf{w}(\mathbf{x}, t_0) = \mathbf{u}(\mathbf{x}, t_0).
  \end{equation}
\end{defn}

(\ref{eq:GePUPe}b) is different from (\ref{eq:GePUP}b). 
First,
 the boundary condition $\mathbf{n}\cdot\mathbf{w}=0$ 
 in (\ref{eq:GePUP}b) is removed 
 and the term $\lambda \mathbf{n}\cdot\mathbf{w}$ is added 
 in (\ref{eq:GePUPe}f)
 so that any nonzero $\mathbf{n}\cdot\mathbf{w}$
 decays exponentially towards zero; 
 see Lemma \ref{lem:decayOfNormalVelGePUPe}. 
Second, 
the boundary condition $\nabla \cdot \mathbf{w}= 0$
 is added in (\ref{eq:GePUPe}b) to set up
 an exponential decay of $\nabla \cdot \mathbf{w}$;
 see Theorem \ref{thm:decayOfVelDivGePUP}. 
Third, the boundary condition
 $\mathbf{w}\cdot\boldsymbol{\tau}=0$
 is added in (\ref{eq:GePUPe}b) to
 close the vector diffusion equation governing the evolution of $\mathbf{w}$.
Lastly, the boundary condition
$\mathbf{u}\cdot\boldsymbol{\tau}=0$ in (\ref{eq:GePUP}b)
is removed because, 
as a perturbed version of $\mathbf{u}$,
the non-solenoidal velocity $\mathbf{w}$
\emph{is} $\mathbf{u}$
when the initial condition (\ref{eq:GePUP-initialW}) is imposed;
see Lemma \ref{lem:wISu}. 
Another initial condition weaker than (\ref{eq:GePUP-initialW}) is 
 \begin{equation}
   \label{eq:GePUP-initialW-bdry}
   \forall \mathbf{x}\in \partial{\Omega},\quad
   \mathbf{w}(\mathbf{x}, t_0) = \mathbf{u}(\mathbf{x}, t_0).
 \end{equation}

To connect $\mathbf{w}$ to $\mathbf{u}$, 
(\ref{eq:GePUPe}c) appears to be the most natural choice.

Compared to the formulation in
 \cite{shirokoff11:_navier_stokes,rosales21:_high_poiss_navier_stokes},
 the GePUP-E formulation (\ref{eq:GePUPe}) facilitates
 the design and analysis of numerical schemes
 that treat the nonlinear convection term and the pressure gradient term explicitly;
 see Definition \ref{def:GePUP-ES-RK}.

\begin{lem}
  The Neumann BVP (\ref{eq:GePUPe}c,d) admits a unique solution
   of $\mathbf{u}$.
\end{lem}
\begin{proof}
  (\ref{eq:GePUPe}c,d) and
   the definition $\mathbf{w} = \mathbf{u}-\nabla\phi$ yield 
   \begin{equation}
     \label{eq:BVPprojectW}
     \begin{array}{rl}
       \Delta\phi &= -\nabla\cdot \mathbf{w}
                    \qquad \text{ in } \Omega,
       \\
       \mathbf{n}\cdot \nabla \phi &= -\mathbf{n}\cdot\mathbf{w}
                                     \qquad    \text{ on } \partial \Omega.
     \end{array}
   \end{equation}
  Then the proof is completed by
   Theorem \ref{thm:ExistenceAndUniquenessOfPoissonNeumannBVP}
   and the divergence theorem.
\qed 
\end{proof}



\subsection{The exponential decay
  of $\mathbf{n}\cdot\mathbf{w}$
  and $\nabla\cdot\mathbf{w}$
  in (\ref{eq:GePUPe})}

Although the no-penetration condition
 $\mathbf{n}\cdot\mathbf{w}=0$ is not explicitly stated
 in (\ref{eq:GePUPe}),
 the exponential decay of
 $\mathbf{n}\cdot\mathbf{w}$ on $\partial\Omega$
 is guaranteed by

\begin{lem}
  \label{lem:decayOfNormalVelGePUPe}
  The GePUP-E formulation (\ref{eq:GePUPe}) satisfies
   \begin{equation}
     \label{eq:normalVel}
     \frac{\partial (\mathbf{n}\cdot\mathbf{w})}{\partial t} = -\lambda\mathbf{n}\cdot\mathbf{w} \quad \mathnormal{on~}\partial\Omega
   \end{equation}
   and thus 
   $\mathbf{n}\cdot\mathbf{w}(t)
   = e^{-\lambda (t-t_0)}\mathbf{n}\cdot\mathbf{w}(t_0)$ holds on $\partial\Omega$.
   In particular, we have
   \begin{displaymath}
     \label{eq:normalVelZero}
     \mathbf{n}\cdot\mathbf{w}(t_0) = 0\ \Rightarrow\ 
     \forall t> t_0,\ \mathbf{n} \cdot \mathbf{w}(t)=0.
   \end{displaymath}
\end{lem}
\begin{proof}
  (\ref{eq:normalVel}) follows immediately from
  (\ref{eq:GePUPe}f) and 
  the normal component of (\ref{eq:GePUPe}a).
\qed 
\end{proof}


The exponential decay of $\nabla \cdot \mathbf{w}$
 is guaranteed by
 
\begin{thm}
  \label{thm:decayOfVelDivGePUP}
  The GePUP-E formulation (\ref{eq:GePUPe}) satisfies
   \begin{equation}
     \label{eq:heatEqOfDivW}
     \frac{\partial\left(\nabla\cdot \mathbf{w}\right)}{\partial t}
     = \nu\Delta\left(\nabla\cdot \mathbf{w}\right)
     \quad \text{ in } \Omega, 
   \end{equation}
   which implies 
   \begin{equation}
     \label{eq:DecayOfDivW}
     \frac{\dif}{\dif t}\|\nabla\cdot \mathbf{w}\|^2
     = -2\nu\left\|\nabla\left(\nabla\cdot \mathbf{w}\right)\right\|^2, 
   \end{equation}
   which further implies
   \begin{equation}
     \label{eq:GePUP-DecayOfDivW}
     \|\nabla\cdot\mathbf{w}(t)\| \le
     e^{-\nu C(t-t_0)} \|\nabla\cdot\mathbf{w}(t_0)\|
   \end{equation}
   where $t_0$ is the initial time
   and $C$ a positive constant independent of $\mathbf{w}$.
\end{thm}
\begin{proof}
  (\ref{eq:heatEqOfDivW}) follows from
   the divergence of (\ref{eq:GePUPe}a), 
   (\ref{eq:GePUPe}e), and the commutativity of 
   $\nabla\cdot$ and $\Delta$. 
  Multiply (\ref{eq:heatEqOfDivW}) with
   $\nabla\cdot \mathbf{w}$, integrate over $\Omega$,
   and we have
   \begin{displaymath}
     \begin{array}{rl}
     & \frac{1}{2}\frac{\dif}{\dif t}\|\nabla\cdot \mathbf{w}\|^2
       = \nu\left\langle\Delta\left(\nabla\cdot \mathbf{w}\right),
       \nabla\cdot \mathbf{w} \right\rangle
     \\
     =& -\nu\left\|\nabla\left(\nabla\cdot \mathbf{w}\right)\right\|^2
        + \nu\oint_{\partial\Omega}\left(\nabla\cdot \mathbf{w}\right)
        \frac{\partial\left(\nabla\cdot \mathbf{w}\right)}{\partial \mathbf{n}}\dif A 
        = -\nu\left\|\nabla\left(\nabla\cdot \mathbf{w}\right)\right\|^2,
     \end{array}
   \end{displaymath}
   where the second step follows from Green's formula
   (\ref{eq:GreenFormula})
   with $u=v=\nabla\cdot\mathbf{w}$ and
   the last from (\ref{eq:GePUPe}b).
 Then the Poincar\'e-Friedrichs inequality
  and (\ref{eq:DecayOfDivW}) imply
  $\frac{\dif}{\dif t}\|\nabla\cdot \mathbf{w}\|^2 \le
  -2\nu C\|\nabla\cdot \mathbf{w}\|^2$
  and then (\ref{eq:GePUP-DecayOfDivW}) follows from
  \begin{displaymath}
    \begin{array}{l}
    \frac{\dif}{\dif t}\left(e^{2\nu Ct}\|\nabla\cdot\mathbf{w}\|^2
    \right)
    = e^{2\nu Ct}\left( 2\nu C \|\nabla\cdot\mathbf{w}\|^2
      + \frac{\dif}{\dif t} \|\nabla\cdot\mathbf{w}\|^2 \right)
    \le 0. 
    \end{array}
\quad\qed
  \end{displaymath}
\end{proof}

Due to the exponential decay in (\ref{eq:GePUP-DecayOfDivW}),
 the convergence of $\mathbf{w}$ to $\mathbf{u}$
 should be sufficiently fast for most practical applications.
In numerical simulations via GePUP-E,
 it is suggested to 
 set the initial condition of $\mathbf{w}$
 to that of $\mathbf{u}$. 
Then the mechanism of exponential decay
 in Theorem \ref{thm:decayOfVelDivGePUP} 
 will suppress divergence residue
 caused by truncation errors of spatial operators;
 see Theorem \ref{thm:wConverge2u}.

\subsection{Recovering INSE from GePUP-E}
\label{sec:equiv-GePUP-INSE}

\begin{lem}
  \label{lem:wISu}
  With the initial condition (\ref{eq:GePUP-initialW}),
   the GePUP-E formulation (\ref{eq:GePUPe})
   recovers the INSE in (\ref{eq:INS}) for all $t\ge t_0$.
\end{lem}
\begin{proof}
  (\ref{eq:GePUP-initialW}) and (\ref{eq:GePUP-DecayOfDivW})
   yield $\nabla\cdot\mathbf{w}(t)=0$.
  With $\mathbf{w}\cdot\mathbf{n} = 0$ in
   Lemma \ref{lem:decayOfNormalVelGePUPe},
   the BVP (\ref{eq:BVPprojectW})
   reduces to a Laplace equation with
   homogeneous Neumann conditions,
   for which 
   $\phi$ being a constant in $\Omega$
   is a particular solution.
  Thus $\nabla\phi = \mathbf{0}$
   and 
   $\mathbf{w}= \mathbf{u}$.
\qed 
\end{proof}

Lemma \ref{lem:wISu} implies
 the well-posedness of (\ref{eq:GePUPe}e,f).
In comparison to (\ref{eq:GePUP}f),
 (\ref{eq:GePUPe}f) contains the convection term
 so that the well-posedness of the Neumann BVP (\ref{eq:GePUPe}e,f)
 is independent of the boundary condition of $\mathbf{u}$. 
 
\begin{lem}
  \label{lem:GePUP-PPE-Wellposed}
  With the initial condition (\ref{eq:GePUP-initialW}),
  the Neumann BVP (\ref{eq:GePUPe}e,f) admits a unique solution
  of $\nabla q$.
\end{lem}
\begin{proof}
  The divergence theorem
  and the commutativity of $\nabla\cdot$ and $\Delta$ imply
   \begin{displaymath}
     \oint_{\partial\Omega}
     \left(
       \mathbf{n}\cdot\nu\Delta\mathbf{w}
       + \lambda \mathbf{n}\cdot\mathbf{w}
     \right) \dif A
     =
     \int_{\Omega} 
     \left(
       \nu\Delta\nabla\cdot\mathbf{w} + \lambda \nabla\cdot\mathbf{w}
     \right) \dif V =0,
   \end{displaymath}
   where the second step follows from Lemma \ref{lem:wISu}. 
   The divergence theorem gives
   \begin{displaymath}
     \int_{\Omega}\nabla\cdot\left(
       \mathbf{g}-\mathbf{u}\cdot\nabla \mathbf{u}\right)\dif V
     = \oint_{\partial\Omega} \mathbf{n}\cdot
     (\mathbf{g}-\mathbf{u}\cdot\nabla \mathbf{u}) \,\dif A.
   \end{displaymath}
   Then the rest of the proof follows from
   Theorem \ref{thm:ExistenceAndUniquenessOfPoissonNeumannBVP}.
\qed 
\end{proof}

\begin{thm}
  \label{thm:wConverge2u}
  GePUP-E (\ref{eq:GePUPe})
  with the initial condition (\ref{eq:GePUP-initialW-bdry}) 
  satisfies
  \begin{displaymath}
    \forall \epsilon  >  0,\ \exists t^* >  t_0,\  
    \text{ s.t. } \forall t > t^*,\quad 
    \left\{
      \begin{array}{l}
        \sup_{\mathbf{x}\in\Omega}\left|\nabla \cdot \mathbf{w}(\mathbf{x},t)\right| <  \epsilon
        \\
        \left\|\mathbf{w}(\mathbf{x},t)-\mathbf{u}(\mathbf{x},t)\right\|  <  \epsilon,
        \\
        \left\|q(\mathbf{x}, t) - p(\mathbf{x}, t) \right\| < \epsilon.
      \end{array}
      \right.
  \end{displaymath}
\end{thm}
\begin{proof}
  The first limit
  follows from (\ref{eq:GePUP-DecayOfDivW}).
  The second limit then follows
  from arguments similar to those in the proof of Lemma \ref{lem:wISu}.
  The third limit is a consequence of the second limit.
\qed 
\end{proof}

\subsection{The equivalence of GePUP-E with INSE}
\label{sec:derive-GePUP-INSE}

To derive GePUP-E from INSE,
 we split $\mathbf{u}$ into $\nabla\phi$ and $\mathbf{w}$
 and suppress $\nabla\cdot\mathbf{w}$
 via a heat equation.

\begin{lem}
  \label{lem:GePUPfromINSE}
  The GePUP-E formulation (\ref{eq:GePUPe}) is derived
  from the INSE in (\ref{eq:INS})
   by splitting the velocity $\mathbf{u}$
   as $\mathbf{u}=\mathbf{w}+\nabla\phi$
   and demanding that $\mathbf{w}$ satisfy (\ref{eq:GePUPe}b) and (\ref{eq:heatEqOfDivW}).
\end{lem}
\begin{proof}
  For an arbitrary vector field $\mathbf{u}$,
   (\ref{eq:cProj}) and (\ref{eq:opB}) yield
   \begin{displaymath}
     \Delta \cProj \mathbf{u}
     = \Delta(\mathbf{u} - \nabla\phi)
     = \Delta \mathbf{u} - \nabla\nabla\cdot\nabla\phi
     = {\mathcal B} \mathbf{u} +\nabla
     \nabla\cdot\cProj\mathbf{u}, 
   \end{displaymath}
   and thus we have
   \begin{equation}
     \label{eq:LaplProj}
     \Delta \cProj 
     = \Delta - \nabla \nabla\cdot + \nabla \nabla\cdot\cProj.
   \end{equation}
  
We deduce (\ref{eq:GePUPe}a) from 
\begin{equation}
  \label{eq:deriveGePUP1}
  \begin{array}{rl}
    \frac{\partial \mathbf{w}}{\partial t}
    =& \frac{\partial \cProj\mathbf{u}}{\partial t}
       = \cProj\frac{\partial \mathbf{u}}{\partial t}
       = \cProj\mathbf{a} - \nu(\Delta \mathbf{u} -
       \nabla\nabla\cdot \mathbf{u} 
       +\nabla\nabla\cdot \cProj\mathbf{u}
       - \Delta \cProj\mathbf{u})
    \\
    =& \cProj\mathbf{a} - \mathbf{a}^{*}
       + \mathbf{g}-\mathbf{u}\cdot\nabla \mathbf{u}
       +\nu\nabla\nabla\cdot \mathbf{u}
       -\nu\nabla\nabla\cdot \mathbf{w}
       +\nu\Delta \mathbf{w}
    \\
    =& \mathbf{g}-\mathbf{u}\cdot \nabla \mathbf{u}-\nabla q
       +\nu\Delta \mathbf{w}, 
  \end{array}
\end{equation}
 where the third step follows from (\ref{eq:LaplProj}),
 the fourth from (\ref{eq:EulerianAccelerations}), 
 and the last from 
 \begin{equation}
   \label{eq:deriveGePUPe_gradP}
   \nabla q := \mathbf{a}^{*}-\cProj\mathbf{a}-\nu\nabla\nabla\cdot
   \mathbf{u}+\nu\nabla\nabla\cdot \mathbf{w}; 
 \end{equation}
 the above RHS is indeed a gradient
 because of (\ref{eq:cProj}) and (\ref{eq:INSEa-EulerianAccel}).

(\ref{eq:GePUPe}e) follows from (\ref{eq:heatEqOfDivW}) 
 and the divergence of \eqref{eq:deriveGePUPe_gradP}, i.e., 
 \begin{displaymath}
   \Delta q
   = \nabla\cdot\left(\mathbf{g}-\mathbf{u}\cdot\nabla
     \mathbf{u}\right) -
   \left[\frac{\partial \nabla\cdot\mathbf{w}}{\partial t}
     - \nu (\Delta \nabla\cdot\mathbf{w})
   \right]
   \qquad \text{ in } \Omega.
 \end{displaymath}

(\ref{eq:GePUPe}f) follows from
 Lemma \ref{lem:decayOfNormalVelGePUPe} and 
 the normal components of (\ref{eq:GePUPe}a) on $\partial\Omega$, 
 i.e., $\mathbf{n}\cdot \nabla q = \mathbf{n}\cdot(\mathbf{g}
   - \mathbf{u}\cdot\nabla\mathbf{u}
   +\nu\Delta\mathbf{w})
   + \frac{\partial(\mathbf{n}\cdot\mathbf{w})}{\partial t}$.
   \qed
%
\end{proof}

An alternative interpretation of (\ref{eq:deriveGePUP1})
 might be illuminating.
A given scalar field $\phi$ furnishes
 a specific generic projection
 $\cProj_{\phi}\mathbf{u}=\mathbf{u}-\nabla\phi$
 that perturbs $\mathbf{u}$ to be non-solenoidal. 
For any $\phi$, the composite operator
 $\cProjLH\circ\cProj_{\phi}$ is the identity
 on the divergence-free vector space, 
 i.e., $\nabla\cdot\mathbf{u}=0$ implies
 $\cProjLH\circ\cProj_{\phi} \mathbf{u}=\cProjLH\mathbf{w} = \mathbf{u}$
 for any ${\mathcal C}^1$ scalar field $\phi$. 
If the evolution of $\mathbf{w}=\cProj_{\phi}\mathbf{u}$
 in (\ref{eq:deriveGePUP1}) did not
 have the exponential decay of $\nabla\cdot\mathbf{w}$,
 we would have to apply $\cProjLH$
 to $\frac{\partial \mathbf{w}}{\partial t}$
 to recover the INSE.
However, 
 thanks to the exponential decay of $\nabla\cdot\mathbf{w}$
 in Theorem \ref{thm:decayOfVelDivGePUP}, 
 $\mathbf{w}$ converges to $\mathbf{u}$
 in Lemma \ref{thm:wConverge2u} 
 and thus there is no need 
 to recover $\mathbf{u}$ from $\mathbf{w}$. 

\begin{thm}
  \label{thm:GePUPeEquivINSE}
  GePUP-E in Definition \ref{def:GePUPe-equations}
  is equivalent to the INSE in (\ref{eq:INS}).
\end{thm}
\begin{proof}
  This follows directly from Lemma \ref{lem:wISu} and Lemma \ref{lem:GePUPfromINSE}.
\qed 
\end{proof}

\subsection{The monotonic decrease of the total kinetic energy}
\label{sec:monot-decr-kinet}


\begin{lem}
  \label{lem:orthogonalityRelation}
  A vector field $\mathbf{v}\in\mathcal{C}^1(\Omega)$
   with $\nabla\cdot\mathbf{v}=0$ in $\Omega$
   and 
   $\mathbf{v}\cdot\mathbf{n}=0$ on $\partial\Omega$ 
   is orthogonal to the gradient field of any scalar function
   $\phi\in{\mathcal C}^1(\Omega)$,
   i.e., 
   \begin{equation}
     \label{eq:orthogonalityRelation}
     \left\langle \mathbf{v},\nabla\phi\right\rangle = \int_{\Omega}\mathbf{v}\cdot \nabla \phi\,\dif V = 0.
   \end{equation}
\end{lem}
\begin{proof}
  By the chain rule, 
   we have
   \begin{displaymath}
     \int_{\Omega}\mathbf{v}\cdot \nabla \phi\,\dif V
     = \int_{\Omega}\nabla\cdot(\phi\mathbf{v})\dif V - \int_{\Omega}\phi\nabla\cdot \mathbf{v}\,\dif V 
     = \int_{\partial \Omega}\phi\mathbf{v}\cdot \mathbf{n}\,\dif A = 0,
   \end{displaymath}
   where the second step follows from the divergence theorem and
   $\nabla\cdot\mathbf{v}=0$ in $\Omega$ 
   and the last from $\mathbf{v}\cdot\mathbf{n}=0$ on
   $\partial\Omega$. 
\qed 
\end{proof}

As a benefit of Lemma \ref{lem:orthogonalityRelation},
 the accuracy of the computed velocity 
 is largely decoupled from that of the pressure gradient; 
this orthogonality condition can be enforced to machine precision
 in the FV formulation. 

\begin{defn}
  The \emph{kinetic energy} of a fluid
  with velocity $\mathbf{u}$ 
  is 
\begin{equation}
  \label{eq:kineticEnergy}
  E_{\text{kinetic}} := \frac{1}{2}\|\mathbf{u}\|^2 = \frac{1}{2}\int_{\Omega}\mathbf{u}\cdot \mathbf{u}\,\dif V.
\end{equation}
\end{defn}

\begin{thm}
  \label{thm:energyDissipationOfGePUP}
  Suppose the body force $\mathbf{g}$ in (\ref{eq:GePUPe}a)
  is conservative, i.e., $\mathbf{g} = -\nabla\varphi$
  for some scalar field $\varphi$ in $\Omega$.
  Then the evolution of the kinetic energy
  in the GePUP-E formulation
  in Definition \ref{def:GePUPe-equations} is governed by
  \begin{displaymath}
    \frac{\dif}{\dif t}E_{\text{kinetic}} = -\nu\|\nabla
    \mathbf{u}\|^2
    := -\nu \sum_{i=1}^{\Dim}\sum_{j=1}^{\Dim} \int_{\Omega}\left|\frac{\partial u_i}{\partial x_j}\right|^2\dif V.
  \end{displaymath}
\end{thm}
\begin{proof}
  Since $\varphi$ can be absorbed into $q$, 
  the assumption $\mathbf{g}=\mathbf{0}$
  has no loss of generality. 
  The inner product of $\mathbf{u}$
   and the momentum equation (\ref{eq:GePUPe}a) give
   \begin{equation*}
     (*):\quad
     \left\langle \frac{\partial \mathbf{w}}{\partial t}, \mathbf{u}\right\rangle
     = - \left\langle \mathbf{u}\cdot\nabla \mathbf{u},
       \mathbf{u}\right\rangle
     - \left\langle \nabla q, \mathbf{u}\right\rangle
     + \left\langle\nu\Delta \mathbf{w}, \mathbf{u}\right\rangle.
   \end{equation*}

  The left-hand side (LHS) is computed as 
   \begin{equation}
     \label{eq:w-InnerProd-u}
     \renewcommand{\arraystretch}{1.2}
     \begin{array}{rl}
       \left\langle \frac{\partial \mathbf{w}}{\partial t},
       \mathbf{u}\right\rangle
       &= \left\langle \frac{\partial (\mathbf{u}-\nabla\phi)}{\partial t},
         \mathbf{u} \right\rangle
         = \left\langle \frac{\partial \mathbf{u}}{\partial t},
         \mathbf{u} \right\rangle
         - \left\langle \nabla \frac{\partial\phi}{\partial t},
         \mathbf{u} \right\rangle
         = \left\langle \frac{\partial \mathbf{u}}{\partial t},
         \mathbf{u}\right\rangle
       \\
       & = \frac{\dif }{\dif
         t}\left(\frac{1}{2}\int_{\Omega}\mathbf{u}\cdot \mathbf{u}\,\dif
         V\right)
         = \frac{\dif }{\dif t}E_{\text{kinetic}}, 
     \end{array}
   \end{equation}
   where 
   the second step follows from
   the commutativity of $\partial_t$ and $\nabla$,
   the third from (\ref{eq:GePUPe}c,d) and
   Lemma \ref{lem:orthogonalityRelation},
   and the last from (\ref{eq:kineticEnergy}). 

  The first RHS term in $(*)$ vanishes because
   \begin{displaymath}
     \renewcommand{\arraystretch}{1.2}
     \begin{array}{rl}
       & \left\langle \mathbf{u}\cdot \nabla \mathbf{u},
         \mathbf{u}\right\rangle
         = \int_{\Omega}\mathbf{u}\cdot\left(\mathbf{u}\cdot \nabla \mathbf{u}\right)\dif V = \int_{\Omega}u_i\left(u_j \frac{\partial u_i}{\partial x_j}\right)\dif V 
         = \frac{1}{2}\int_{\Omega}u_j \frac{\partial (u_iu_i)}{\partial
         x_j}\dif V \\
       =& -\frac{1}{2}\int_{\Omega}u_iu_i \frac{\partial u_j}{\partial x_j}\dif V + \frac{1}{2}\int_{\partial \Omega}u_iu_iu_jn_j\,\dif A \\
       =& -\frac{1}{2}\int_{\Omega}|\mathbf{u}|^2\left(\nabla\cdot \mathbf{u}\right)\dif V +
          \frac{1}{2}\int_{\partial\Omega}|\mathbf{u}|^2\left(\mathbf{u}\cdot
          \mathbf{n}\right)\dif A = 0,
     \end{array}
   \end{displaymath}
   where the fourth step, 
   in Einstein summation convention,
   follows from Lemma \ref{lem:integrationByParts}.
  The second term in $(*)$ is also zero
   due to Lemma \ref{lem:orthogonalityRelation}.
%
  The third term is 
  \begin{equation}
    \label{eq:innerProdLapWandU}
    \begin{array}{rl}
      & \left\langle\Delta \mathbf{w}, \mathbf{u}\right\rangle
        = \left\langle \Delta (\mathbf{u}-\nabla\phi), \mathbf{u} \right\rangle
        = \left\langle \Delta \mathbf{u}, \mathbf{u} \right\rangle
        - \left\langle \nabla\Delta\phi, \mathbf{u} \right\rangle
        = \left\langle \Delta \mathbf{u}, \mathbf{u} \right\rangle
      \\
      =& \int_{\Omega}u_i\Delta u_i\,\dif V
         = -\int_{\Omega}\nabla u_i\cdot \nabla u_i\,\dif V
         +\int_{\partial\Omega}u_i
         \frac{\partial u_i}{\partial \mathbf{n}}\dif A
       = -\left\|\nabla \mathbf{u}\right\|^2,
    \end{array}
  \end{equation}
  where
  the first step follows from the definition of $\mathbf{w}$,
  the second from the commutativity of $\Delta$ and $\nabla$,
  the third from (\ref{eq:GePUPe}c,d) and
  Lemma \ref{lem:orthogonalityRelation},
  the penultimate from Green's formula (\ref{eq:GreenFormula}), 
  and the last from the no-slip condition of $\mathbf{u}$,
  which holds from Lemma \ref{lem:wISu} and (\ref{eq:GePUPe}b,d).
\qed 
\end{proof}

\subsection{Prominent features of the GePUP-E formulation}

{
\setlist[enumerate]{label=(\alph{enumi})}
\begin{enumerate}
 \item The sole evolutionary variable 
   is the non-solenoidal velocity $\mathbf{w}$,
   with $\mathbf{u}$ determined from $\mathbf{w}$
   via (\ref{eq:GePUPe}c,d)
   and $q$ from $\mathbf{u}$ and $\mathbf{w}$ 
   via (\ref{eq:GePUPe}e,f).
   This chain of determination
   $\mathbf{w}\rightarrow \mathbf{u} \rightarrow q$
   from Neumann BVPs
   is \emph{instantaneous}
   and has nothing to do with time integration. 
   Therefore, 
   a time integrator in MOL can be employed in a black-box manner. 
 \item There is no ambiguity on which velocities 
   should be projected and which should not in MOL;
   this resolves (TMR-1) discussed in Section \ref{sec:UPPE}. 
 \item Now that the main evolutionary variable $\mathbf{w}$
   in (\ref{eq:GePUPe}) 
   is formally non-solenoidal,
   the Leray-Helmholtz projection $\cProjLH$
   only comes into the system (\ref{eq:GePUPe}) on the RHS.
   Although still contributing to the local truncation error, 
   the approximation error of a discrete projection
   to $\cProjLH$ does not affect numerical stability of MOL; 
   this resolves (TMR-2) in Section \ref{sec:UPPE}. 
 \item (\ref{eq:GePUPe}) comes with the built-in mechanisms
   of exponential decays of velocity divergence and total kinetic energy, 
   which are conducive to the design of semi-discrete algorithms that
   ensure numerical stability
   and preserve physical structures of incompressible flows,
   c.f. Theorems
   \ref{thm:decayOfVelDiv-GePUP-SAV-RK} and
   \ref{thm:GePUPSAVRKEnergyDecay}.
 \end{enumerate}
}



%% file: sec/SAV.tex
\section{The GePUP-ES formulation}
\label{sec:GePUP-ES}

The SAV approach,
 as introduced in \cite{shen18:SAV,shen19:SAV_Review},
 has been proposed to
 develop time discretization schemes
 that are both efficient and stable for gradient flows.
This approach was originally designed to
 create schemes that are linear, decoupled, unconditionally stable,
 and up to second-order accurate.
It is also successfully extended to
address the Navier-Stokes equations in
\cite{li22:_new_sav_pressure_correction_inse,lin19:_numer_appx_inses_auxiliary_energy_variable,li20:_error_analysis_sav_mac}.

More recently,
based on the generalized SAV approach \cite{huang22:_new_class_imex_bdfk_sav}, 
Huang et al. have devised high-order consistent splitting schemes
for the Navier-Stokes equations,
with periodic boundary conditions in
\cite{huang21:_stability_error_analysis_higher_order_imex_ns_periodic_bc} and
no-slip boundary conditions in
\cite{wu22:_new_class_higher_order_decoupled_inses_rotating}.

In this section,
 we couple the GePUP-E formulation (\ref{eq:GePUPe})
  to the SAV approach
  introduced in \cite{shen18:SAV,shen19:SAV_Review} to
  deal with the nonlinear convection term
  $\mathbf{u}\cdot\nabla \mathbf{u}$
  so that unconditionally energy-stable numerical schemes
  can be constructed.
\begin{defn}
  \label{def:GePUPSAV}
  The \emph{GePUP-ES formulation of INSE on no-slip domains}
is
  \begin{subequations}
    \label{eq:GePUPSAV}
    \begin{alignat}{2}
      \frac{\partial \mathbf{w}}{\partial t}
      &= \mathbf{g}-r(t)\mathbf{u}\cdot \nabla \mathbf{u}
      -\nabla q+\nu\Delta \mathbf{w} &\quad& \text{in } \Omega,
      \\
      \mathbf{w}\cdot \boldsymbol{\tau} &= 0, \ \
      \nabla \cdot \mathbf{w} = 0 \ \
      && \text{on } \partial \Omega, \\
      \frac{\dif r}{\dif t} &= I_{cv}(\mathbf{u},\mathbf{u}), && \\
      \mathbf{u} &= \mathscr{P}\mathbf{w} && \text{in } \Omega, \\
      \mathbf{u}\cdot \mathbf{n} &= 0 && \text{on } \partial\Omega,
      \\
      \Delta q &= \nabla\cdot\left(\mathbf{g}
        - r(t)\mathbf{u}\cdot \nabla \mathbf{u}\right)
      && \text{in } \Omega, \\
      \mathbf{n}\cdot \nabla q &= \mathbf{n}\cdot(\mathbf{g}
      - r(t)\mathbf{u}\cdot\nabla\mathbf{u}
      +\nu\Delta\mathbf{w})
      + \lambda\mathbf{n}\cdot\mathbf{w} && \text{on } \partial \Omega,
    \end{alignat}
  \end{subequations}
  where
  $I_{cv}(\mathbf{u},\mathbf{v}):=\int_{\Omega}\left(\mathbf{u}\cdot
    \nabla \mathbf{u}\right)\cdot \mathbf{v}\,\dif V$,
  the SAV $r(t) \equiv 1$, 
  and $\lambda$ is a nonnegative penalty parameter. 
\end{defn}

The introduction of the SAV $r(t)\equiv 1$
 immediately implies $\frac{\dif r}{\dif t}=0$.
We \emph{can} define the evolution of $r(t)$
 as the ODE in (\ref{eq:GePUPSAV}c)
 because, for no-slip conditions, 
 we always have 
 $\left\langle \mathbf{u}\cdot\nabla\mathbf{u},
   \mathbf{u}\right\rangle = 0$. 
As a newly added evolutionary variable,
 the SAV $r(t)$ is a double-edged sword.
On the one hand,
 it leads to a tighter coupling
 between $\mathbf{u}$ and $\mathbf{w}$, 
 which makes it difficult
 to orchestrate an implicit or semi-implicit RK method
 as solving a sequence of \emph{linear} systems;
 see the discussions in Section \ref{sec:GePUP-SAV-SDIRK}.
On the other hand,
 it preserves the monotonic decrease of the kinetic energy; 
 see Theorems
 \ref{thm:modifiedEnergyDissipationGePUPSAV}
 and \ref{thm:GePUPSAVRKEnergyDecay}.
 
\begin{thm}
  \label{thm:modifiedEnergyDissipationGePUPSAV}
  Suppose 
  that the body force $\mathbf{g}$ is conservative.
  Then the energy dissipation of
  the GePUP-ES formulation (\ref{eq:GePUPSAV}) is governed by
  \begin{displaymath}
    \frac{\dif}{\dif t}\left(E_{\mathrm{kinetic}} + \frac{r^2}{2}\right)
    = -\nu\|\nabla \mathbf{u}\|^2.
  \end{displaymath}
\end{thm}
\begin{proof}
  Take inner product with $\mathbf{u}$ in (\ref{eq:GePUPSAV}a),
  multiply (\ref{eq:GePUPSAV}c) by $r(t)$,
  add up the resulting two equations,
  and we cancel the integral of the convection term
  to obtain
$\left\langle \frac{\partial \mathbf{w}}{\partial t},
  \mathbf{u}\right\rangle + r(t)r'(t)
= \left\langle \mathbf{g}-\nabla q+\nu\Delta \mathbf{w}, \mathbf{u}\right\rangle$.
  The rest of the proof follows from (\ref{eq:w-InnerProd-u}),
  (\ref{eq:innerProdLapWandU}),
  and Lemma \ref{lem:orthogonalityRelation}. 
\qed 
\end{proof}

Theorem \ref{thm:decayOfVelDivGePUP} also holds for GePUP-ES 
 with exactly the same proof.
Similarly, the GePUP-ES formulation \eqref{eq:GePUPSAV} retains 
 the advantage of GePUP-E
 that the temporal integration and spatial discretization
 are completely decoupled.
Hence
 the fourth-order finite-volume discrete operators
 in \cite{zhang16:_GePUP} can be reused 
 and we will focus on temporal integration hereafter.



%% file: sec/algorithms.tex
\section{Algorithms}
\label{sec:algorithms}

Based on the GePUP-ES formulation, 
 we construct numerical algorithms 
 to preserve the monotonic decrease
 of the modified kinetic energy, the divergence residue,
 and the magnitude of the normal velocity on the domain boundary.

\subsection{Semi-discrete GePUP-ES-RK schemes}
\label{sec:GePUP-ES-RK}

These schemes follow directly
 from discretizing
 the GePUP-ES formulation (\ref{eq:GePUPSAV}) in time 
 by the RK method (\ref{eq:RungeKutta}).

\begin{defn}
   \label{def:GePUP-ES-RK}
   A \emph{GePUP-ES-RK scheme
     for solving INSE with no-slip conditions} is
   a semidiscrete algorithm of the form 
  \begin{equation}
    \label{eq:GePUPSAVRKFinal}
    \renewcommand{\arraystretch}{1.2}
    \left\{
      \begin{array}{rlcl}
        \mathbf{w}^{n+1}
        &= \mathbf{u}^n
          +k\sum_{i=1}^{s}b_i\boldsymbol{\rho}^{(i)}
      &\quad& \text{in } \Omega,
      \\
      r^{n+1}
      &= r^n + k\sum_{i=1}^{s}b_i
        I_{cv}(\widetilde{\mathbf{u}}^{(i)}, \mathbf{u}^{(i)}), &&
      \\
      \mathbf{u}^{n+1}
      &= \mathscr{P}\mathbf{w}^{n+1} &\quad& \text{in } \Omega,
      \\
      \mathbf{n}\cdot \mathbf{u}^{n+1} &= 0
      &\quad& \text{on } \partial\Omega,
    \end{array}
    \right.
  \end{equation}
  where $s$ is the number of stages of the employed RK method $(A,
  \mathbf{b}, \mathbf{c})$,
  the integral $I_{cv}(\widetilde{\mathbf{u}}^{(i)}, \mathbf{u}^{(i)})$
  is the same as that in Definition \ref{def:GePUPSAV}, 
  the auxiliary velocity $\widetilde{\mathbf{u}}^{(i)}$ is a suitable
  explicit approximation to $\mathbf{u}(t^n+c_ik)$
  and
  \begin{subequations}
    \label{eq:GePUPeSAVRKIntermediate}
    \begin{alignat}{2}
      \Delta q^{(i)}
      &= \nabla\cdot\left(\mathbf{g}^{(i)}
        -r^{(i)}\widetilde{\mathbf{u}}^{(i)}
        \cdot\nabla\widetilde{\mathbf{u}}^{(i)}\right)
      &\quad &\text{in } \Omega,
      \\
      \mathbf{n}\cdot \nabla q^{(i)} &=
      \mathbf{n}\cdot\left(\mathbf{g}^{(i)}
        -r^{(i)}\widetilde{\mathbf{u}}^{(i)}
        \cdot\nabla\widetilde{\mathbf{u}}^{(i)}
        +\nu\Delta \mathbf{w}^{(i)}
      \right)
      + \lambda \mathbf{n}\cdot\mathbf{w}^{(i)}
      &\quad&\text{on } \partial\Omega,
      \\
      \boldsymbol{\rho}^{(i)} &:=\mathbf{g}^{(i)}-r^{(i)}\widetilde{\mathbf{u}}^{(i)}\cdot\nabla
      \widetilde{\mathbf{u}}^{(i)}-\nabla q^{(i)}+\nu\Delta
      \mathbf{w}^{(i)}, && 
      \\
      \mathbf{w}^{(i)}
      &= \mathbf{u}^n
      + k\sum_{j=1}^{s}a_{i,j}\boldsymbol{\rho}^{(j)}
      &\quad& \text{in } \Omega,
      \\
      \mathbf{w}^{(i)}\cdot\boldsymbol{\tau}&=0,
      ~\nabla \cdot \mathbf{w}^{(i)} = 0
      &\quad& \text{on } \partial\Omega,
      \\
      r^{(i)} &= r^n
      +k\sum_{j=1}^{s}a_{i,j} I_{cv}(\widetilde{\mathbf{u}}^{(j)}, \mathbf{u}^{(j)}),  & &
      \\
      \mathbf{u}^{(i)} &= \mathscr{P}\mathbf{w}^{(i)} &\quad& \text{in } \Omega, \\
      \mathbf{n}\cdot \mathbf{u}^{(i)} &= 0 &\quad& \text{on } \partial\Omega.
    \end{alignat}
  \end{subequations}
\end{defn}

As suggested by (\ref{eq:GePUPeSAVRKIntermediate}d)
 and the first equation of (\ref{eq:GePUPSAVRKFinal}),
 $\mathbf{u}^n$ is used as the initial condition of $\mathbf{w}$ 
 for time integration within the interval $[t^n,t^{n+1}]$.
 However, it is emphasized that we do \emph{not} write
 $\mathbf{w}^n=\mathbf{u}^n$
 because this would cause a notation clash to $\mathbf{w}^{n+1}$,
 which, by (\ref{eq:GePUPeSAVRKIntermediate}),
 is the evolutionary velocity before the final projection
 and thus needs not be divergence free.
In Definition \ref{def:GePUP-ES-RK},
 the approximations $\widetilde{\mathbf{u}}^{(i)}$ in $\Omega$
 may be obtained by polynomial interpolation
 based on stage values of recent time steps.
In this work,
 we fit a cubic polynomial
$\mathbf{p}$ from the known velocity
 $\hat{\mathbf{u}}_n^{(j)}$
 at time instances $t^n+\hat{c}_{j}k$ with $j=0,1,2,3$
 and then approximate $\mathbf{u}(t^n+c_ik)$ with 
\begin{equation}
  \label{eq:interpolationEAV}
  \forall i=1,\ldots,s, \qquad
  \widetilde{\mathbf{u}}_{n}^{(i)}
  = \mathbf{p} \left(c_{i}k \right)
  := \sum_{j=0}^{3} \hat{\mathbf{u}}_n^{(j)}
  \prod \limits_{\ell \neq j}
  \frac{c_{i} - \hat{c}_{\ell}}{\hat{c}_{j} - \hat{c}_{\ell}}; 
\end{equation}
{
\setlist[enumerate]{label=({EAV}-\arabic{enumi})}
\begin{enumerate}
\item for the first two time steps $n=0,1$, 
  we calculate $\hat{\mathbf{u}}^{(j)}$
  by GePUP-ERK \cite{zhang16:_GePUP},
  the explicit RK method   
  for solving the GePUP formulation,
  with time step size $\frac{1}{3}k$
  and then fit $\mathbf{p}$ in (\ref{eq:interpolationEAV})
  with $\left[ \hat{c}_{0}, \hat{c}_{1}, \hat{c}_{2}, \hat{c}_{3}
  \right]
  = \left[ 0, \frac{1}{3}, \frac{2}{3}, 1 \right]$;
\item for $n\ge 2$,
  we first calculate $\hat{\mathbf{u}}^{(3)}=\mathbf{u}^{n+1}$
  by GePUP-ERK
  with time step size $k$ 
  and then fit $\mathbf{p}$ in (\ref{eq:interpolationEAV})
  with $\left[ \hat{c}_{0}, \hat{c}_{1}, \hat{c}_{2}, \hat{c}_{3} \right]
  = \left[ -2, -1, 0, 1 \right]$.
\end{enumerate}
}
We emphasize that $\hat{\mathbf{u}}_n^{(j)}$, 
$\widetilde{\mathbf{u}}_{n}^{(i)}$,
and $\mathbf{u}(t^n+c_ik)$ are all at the same location
so that the interpolation is only in time; 
see Figure \ref{fig:interpolationEAV} for an illustration.

\begin{figure}
  \centering
  \input{tikz/interpolationEAU.tex}  
  \caption{Estimating the auxiliary velocity
    $\widetilde{\mathbf{u}}^{(i)}$ (EAV)
    in GePUP-ES-RK by (\ref{eq:interpolationEAV}).}
  \label{fig:interpolationEAV}
\end{figure}

\begin{thm}
  \label{thm:decayOfVelDiv-GePUP-SAV-RK}
  Suppose the RK method
  employed in the GePUP-ES-RK scheme (\ref{eq:GePUPSAVRKFinal})
  is algebraically stable
  in the sense of Definition \ref{def:RKAlgebraicStability}.
  Then (\ref{eq:GePUP-initialW-bdry}) implies 
  \begin{equation}
    \label{eq:decayOfVelDiv-GePUP-SAV-RK}
    \forall n\in \mathbb{N}, \quad
    \left\|\nabla\cdot \mathbf{w}^{n+1}\right\|^2
    - \left\|\nabla\cdot \mathbf{w}^n\right\|^2
    \le -2k\nu\sum_{i=1}^{s}b_i
    \left\|\nabla\left(\nabla\cdot \mathbf{w}^{(i)}\right)\right\|^2, 
  \end{equation}
  \begin{equation}
    \label{eq:decayOfNormalW-GePUP-SAV-RK}
    \forall n\!\in\! \mathbb{N}, \ 
    \forall\mathbf{x}\!\in\! \partial \Omega, ~
    \left|\mathbf{n}\!\cdot\!\mathbf{w}^{n+1}(\mathbf{x})\right|^2
    \!-\! \left|\mathbf{n}\!\cdot\!\mathbf{w}^{n}(\mathbf{x})\right|^2
    \!\leq\! -2k\lambda\sum_{i=1}^{s}b_i
    \left|\mathbf{n}\!\cdot\!\mathbf{w}^{(i)}(\mathbf{x})\right|^2\!. 
  \end{equation}
\end{thm}
\begin{proof} 
  Take divergence of (\ref{eq:GePUPeSAVRKIntermediate}c), 
  apply (\ref{eq:GePUPeSAVRKIntermediate}a) and
  the commutativity of $\Delta$ and $\nabla\cdot$, 
  and we have
  $\nabla\cdot\boldsymbol{\rho}^{(i)}=\nu\Delta\left(\nabla\cdot\mathbf{w}^{(i)}\right)$.
  Then 
  \begin{displaymath}
    \begin{array}{rl}
      &  \left\langle \nabla\cdot \mathbf{w}^{n+1},
        \nabla\cdot \mathbf{w}^{n+1}\right\rangle
      \\
      =& \left\langle
         \nabla\cdot\mathbf{u}^n+k\sum_{i=1}^{s}b_i\nabla\cdot\boldsymbol{\rho}^{(i)},
         \nabla\cdot\mathbf{u}^n +k\sum_{j=1}^{s}b_j\nabla\cdot\boldsymbol{\rho}^{(j)}\right\rangle
      \\
      =& \left\langle \nabla\cdot\mathbf{u}^n, \nabla\cdot\mathbf{u}^n\right\rangle
         +k\sum_{i=1}^{s}b_i\left\langle \nabla\cdot \mathbf{u}^n,
         \nabla\cdot\boldsymbol{\rho}^{(i)}\right\rangle
      \\
      &+k\sum_{j=1}^{s}b_j\left\langle \nabla\cdot \mathbf{u}^n,
        \nabla\cdot\boldsymbol{\rho}^{(j)}\right\rangle
        +k^2\sum_{i=1}^{s}\sum_{j=1}^{s}b_ib_j\left\langle
        \nabla\cdot\boldsymbol{\rho}^{(i)},
        \nabla\cdot\boldsymbol{\rho}^{(j)}\right\rangle
      \\
      =& \left\langle \nabla\cdot\mathbf{u}^n, \nabla\cdot
         \mathbf{u}^n\right\rangle
         +k^2\sum_{i=1}^{s}\sum_{j=1}^{s}b_ib_j\left\langle
         \nabla\cdot\boldsymbol{\rho}^{(i)},
         \nabla\cdot\boldsymbol{\rho}^{(j)}\right\rangle
      \\
      &+k\sum_{i=1}^{s}b_i\left\langle \nabla\cdot
        \mathbf{w}^{(i)}-k\sum_{j=1}^{s}a_{i,j}\nabla\cdot\boldsymbol{\rho}^{(j)},
        \nabla\cdot\boldsymbol{\rho}^{(i)}\right\rangle
      \\
      &+k\sum_{j=1}^{s}b_j\left\langle \nabla\cdot
        \mathbf{w}^{(j)}-k\sum_{i=1}^{s}a_{j,i}\nabla\cdot\boldsymbol{\rho}^{(i)},
        \nabla\cdot\boldsymbol{\rho}^{(j)}\right\rangle
      \\
      =& \left\langle\nabla\cdot \mathbf{u}^n, \nabla\cdot
         \mathbf{u}^n\right\rangle
         - 2\nu k\sum_{i=1}^{s}b_i\left\|\nabla\left(\nabla\cdot \mathbf{w}^{(i)}\right)\right\|^2 \\
       &- k^2\sum_{i=1}^{s}\sum_{j=1}^{s}m_{i,j}\left\langle \nabla\cdot\boldsymbol{\rho}^{(i)}, \nabla\cdot\boldsymbol{\rho}^{(j)}\right\rangle ,
    \end{array}
  \end{displaymath}
  where the last step follows from (\ref{eq:algebraicStabilityMatrix}),
  (\ref{eq:GreenFormula}),
  and (\ref{eq:GePUPeSAVRKIntermediate}e).
 By Definition \ref{def:RKAlgebraicStability},
  the algebraic stability implies
  the symmetric positive semi-definiteness of $M$ and thus 
  we have  $M = O\Lambda O^T$,
  where $O$ is an orthogonal matrix and
  $\Lambda$ is a diagonal matrix with
  $\lambda_{\ell}=\Lambda_{\ell, \ell}\ge 0$
  for each $\ell = 1, \ldots, s$.
 Hence \eqref{eq:decayOfVelDiv-GePUP-SAV-RK} follows from 
  \begin{displaymath}
    \begin{array}{rl}
      &\sum_{i=1}^{s}\sum_{j=1}^{s}m_{i,j}
        \left\langle\nabla\cdot\boldsymbol{\rho}^{(i)},\nabla\cdot\boldsymbol{\rho}^{(j)}\right\rangle
      \\
      =&
         \sum_{i=1}^{s}\sum_{j=1}^{s}\left(\sum_{\ell=1}^{s}o_{i,\ell}\lambda_{\ell}o_{j,\ell}\right)
         \left\langle\nabla\cdot\boldsymbol{\rho}^{(i)},
         \nabla\cdot\boldsymbol{\rho}^{(j)}\right\rangle
      \\
      =& \sum_{\ell=1}^{s}\lambda_{\ell}\left\langle
         \sum_{i=1}^{s}o_{i,\ell}\nabla\cdot\boldsymbol{\rho}^{(i)},
         \sum_{j=1}^{s}o_{j,\ell}\nabla\cdot\boldsymbol{\rho}^{(j)}\right\rangle
      \\
      =& \sum_{\ell=1}^{s}\lambda_{\ell}
         \left\|\sum_{i=1}^{s}o_{i,\ell}\nabla\cdot\boldsymbol{\rho}^{(i)}\right\|^2
         \ge 0
    \end{array}
  \end{displaymath}
  and the fact that the Leray-Helmholtz projection $\cProjLH$ do not increase the divergence of velocity,
  i.e., $\|\nabla \cdot \mathbf{u}^n\|^2 =  \|\nabla \cdot \cProjLH\mathbf{w}^n\|^2 \leq \|\nabla \cdot \mathbf{w}^n\|^2$.

  (\ref{eq:GePUPeSAVRKIntermediate}b)
  and the normal component of (\ref{eq:GePUPeSAVRKIntermediate}c)
  imply
  $\mathbf{n}\cdot\boldsymbol{\rho}^{(i)}
  =-\lambda\mathbf{n}\cdot\mathbf{w}^{(i)}$.
  Then \eqref{eq:decayOfNormalW-GePUP-SAV-RK}
  follows from arguments similar
  to those for (\ref{eq:decayOfVelDiv-GePUP-SAV-RK})
  and
  \begin{displaymath}
    \forall n \in \mathbb{N}, \
    \forall \mathbf{x} \in \partial\Omega, \
    \left|\mathbf{n}\cdot \mathbf{u}^n(\mathbf{x})\right|^2
    - \left|\mathbf{n}\cdot \mathbf{w}^n(\mathbf{x})\right|^2
    = - \left|\mathbf{n}\cdot \mathbf{w}^n(\mathbf{x})\right|^2 \leq 0.
    \hspace{8mm}\qed 
  \end{displaymath}
  
\end{proof}

In the fully discrete sense,
 we discretize each continuous operator
 in (\ref{eq:decayOfVelDiv-GePUP-SAV-RK})
 add the corresponding discretization error, 
 and obtain a discrete version
 of (\ref{eq:decayOfVelDiv-GePUP-SAV-RK}).
As discussed in Section \ref{sec:tests}, 
 this discrete version
 of (\ref{eq:decayOfVelDiv-GePUP-SAV-RK})
 is helpful in understanding the evolution
 of the discrete velocity divergence.

Interestingly,
 (\ref{eq:decayOfVelDiv-GePUP-SAV-RK})
 is useless for first-order finite difference/volume methods.
Suppose $\mathbf{w}^n=\mathbf{w}(t^n)+O(h^2)$.
Then,  due to the factor $k$ and the second derivative
 on the RHS of (\ref{eq:decayOfVelDiv-GePUP-SAV-RK}), 
 the discrete version of (\ref{eq:decayOfVelDiv-GePUP-SAV-RK})
 has an $O(1)$ error.
In contrast,
 a fourth-order method has $\mathbf{w}^n=\mathbf{w}(t^n)+O(h^4)$. 
With the added discretization error approaching zero
 as $h,k\rightarrow 0$, 
 the discrete version of (\ref{eq:decayOfVelDiv-GePUP-SAV-RK})
 indeed has control over the evolution of
 the discrete velocity divergence.
The above discussion suggests
 a crucial advantage of fourth-order methods
 over first-order methods.

\begin{coro}
  \label{coro:divFreeSemiImplicit}
  \hspace{-0.1em}Suppose the RK method
  employed in the GePUP-ES-RK scheme (\ref{eq:GePUPSAVRKFinal})
   is algebraically stable
   in the sense of Definition \ref{def:RKAlgebraicStability}.
  Then the initial condition (\ref{eq:GePUP-initialW}) implies
  \begin{align}
     \label{eq:divFreeSemiImplicit}
    \forall n\in\mathbb{N}^+,&\quad
     \left\{
       \begin{array}{l}
         \left\|\nabla\cdot \mathbf{w}^{n}\right\|=0;\\
         \forall i=1,\ldots,s,\quad
         \nabla\cdot\mathbf{w}^{(i)} = 0, 
       \end{array}
     \right.
     \\
     \label{eq:normalWonBdrySemiImplicit}
    \forall n\in\mathbb{N}^+,&\quad
     \left\{
       \begin{array}{l}
         \mathbf{n}\cdot\mathbf{w}^{n}|_{\partial\Omega}=0;\\
         \forall i=1,\ldots,s,\quad
         \mathbf{n}\cdot\mathbf{w}^{(i)}|_{\partial\Omega} = 0.
       \end{array}
     \right.
   \end{align}
 \end{coro}
\begin{proof}
 (\ref{eq:GePUP-initialW}) gives
  $\nabla\cdot \mathbf{w}(t_0) =\nabla\cdot\mathbf{u}(t_0)=0$
  and
  $\mathbf{n}\cdot\mathbf{w}(t_0)|_{\partial\Omega}=\mathbf{n}\cdot\mathbf{u}(t_0)|_{\partial\Omega}=0$. 
 Then the first clauses of (\ref{eq:divFreeSemiImplicit})
  and (\ref{eq:normalWonBdrySemiImplicit})
  follow from Theorem \ref{thm:decayOfVelDiv-GePUP-SAV-RK}
  and the first condition in Definition \ref{def:RKAlgebraicStability}.
  (\ref{eq:decayOfVelDiv-GePUP-SAV-RK})
  dictates $\|\nabla(\nabla\cdot\mathbf{w}^{(i)})\|=0$
  for each stage. 
 Then the second clause of (\ref{eq:divFreeSemiImplicit})
  follows from the boundary condition
  $\nabla\cdot\mathbf{w}|_{\partial \Omega}=0$ in (\ref{eq:GePUPe}b).
 Similarly,
 \eqref{eq:decayOfNormalW-GePUP-SAV-RK} dictates
 the second clause of (\ref{eq:normalWonBdrySemiImplicit}).
\qed 
\end{proof}

Of course Corollary \ref{coro:divFreeSemiImplicit}
 holds only in the semi-discrete sense.
Since Corollary \ref{coro:divFreeSemiImplicit}
 is used in proving Theorem \ref{thm:GePUPSAVRKEnergyDecay},
 the decay of the modified energy
 in the fully discrete sense
 depends on that of the discrete velocity divergence.

\begin{thm}
  \label{thm:GePUPSAVRKEnergyDecay}
  Suppose that the body force $\mathbf{g}$ in the GePUP-E formulation
   is conservative, that the initial condition of $\mathbf{w}$ is 
   (\ref{eq:GePUP-initialW}), 
   and that
   the employed RK method in (\ref{eq:GePUPSAVRKFinal})
   is algebraically stable
   in the sense of Definition \ref{def:RKAlgebraicStability}.
  Then the energy dissipation of the GePUP-ES-RK scheme
   (\ref{eq:GePUPSAVRKFinal})
   is governed by
   \begin{equation}
     \label{eq:GePUPSAVRKEnergyDecay}
     \mathcal{E}\left(t^{n+1}\right)-\mathcal{E}(t^n) \le
     -k\nu\sum_{i=1}^{s}b_i\left\|\nabla \mathbf{u}^{(i)}\right\|^2,
   \end{equation}
   where the modified energy is defined as 
   $\mathcal{E}(t^n) := \frac{1}{2}\left(\left\|\mathbf{u}^n\right\|^2 +
      \left|r^n\right|^2\right)$.
\end{thm}
\begin{proof}
  Denote $\boldsymbol{\sigma}^{(i)}:= \mathscr{P}\boldsymbol{\rho}^{(i)}$
  and we have from (\ref{eq:GePUPeSAVRKIntermediate}c,d) 
  and (\ref{eq:GePUPSAVRKFinal})
  \begin{displaymath}
    \begin{array}{c}
      \forall i = 1, \ldots, s, \ \ 
      \mathbf{u}^{(i)} = \mathscr{P}\mathbf{w}^{(i)}
      = \mathscr{P}\left(\mathbf{u}^n
      +k\sum_{j=1}^{s}a_{i, j}\boldsymbol{\rho}^{(j)}\right)
      = \mathbf{u}^n+k\sum_{j=1}^{s}a_{i, j}\boldsymbol{\sigma}^{(j)};
      \\
      \mathbf{u}^{n+1} = \mathscr{P}\mathbf{w}^{n+1}
      = \mathscr{P}\left(\mathbf{u}^n
      +k\sum_{i=1}^{s}b_i\boldsymbol{\rho}^{(i)}\right)
      = \mathbf{u}^n+k\sum_{i=1}^{s}b_i\boldsymbol{\sigma}^{(i)}.
    \end{array}
  \end{displaymath}
  It follows that
  \begin{displaymath}
    \begin{array}{rl}
      & \left\langle \mathbf{u}^{n+1}, \mathbf{u}^{n+1}\right\rangle
      = \left\langle
         \mathbf{u}^n+k\sum_{i=1}^{s}b_i\boldsymbol{\sigma}^{(i)},
         \mathbf{u}^n+k\sum_{j=1}^{s}b_j\boldsymbol{\sigma}^{(j)}\right\rangle
      \\
      =& \left\langle \mathbf{u}^n, \mathbf{u}^n\right\rangle
         + k\sum_{i=1}^{s}b_i\left\langle \mathbf{u}^n,
         \boldsymbol{\sigma}^{(i)}\right\rangle
         + k\sum_{j=1}^{s}b_j\left\langle
         \mathbf{u}^n,\boldsymbol{\sigma}^{(j)}\right\rangle
      \\
      &+k^2\sum_{i=1}^{s}\sum_{j=1}^{s}b_ib_j\left\langle
        \boldsymbol{\sigma}^{(i)}, \boldsymbol{\sigma}^{(j)}\right\rangle
      \\
      =& \left\langle \mathbf{u}^n, \mathbf{u}^n\right\rangle
         + k\sum_{i=1}^{s}b_i\left\langle
         \mathbf{u}^{(i)}-k\sum_{j=1}^{s}a_{i,j}\boldsymbol{\sigma}^{(j)},
         \boldsymbol{\sigma}^{(i)}\right\rangle
      \\
      &+ k\sum_{j=1}^{s}b_j\left\langle
        \mathbf{u}^{(j)}-k\sum_{i=1}^{s}a_{j,i}\boldsymbol{\sigma}^{(i)},
        \boldsymbol{\sigma}^{(j)}\right\rangle
      +k^2\sum_{i=1}^{s}\sum_{j=1}^{s}b_ib_j\left\langle
        \boldsymbol{\sigma}^{(i)},
        \boldsymbol{\sigma}^{(j)}\right\rangle
      \\
      =& \left\langle \mathbf{u}^n, \mathbf{u}^n\right\rangle
         + 2k\sum_{i=1}^{s}b_i\left\langle \mathbf{u}^{(i)},
         \boldsymbol{\sigma}^{(i)}\right\rangle
         - k^2\sum_{i=1}^{s}\sum_{j=1}^{s}m_{i,j}\left\langle \boldsymbol{\sigma}^{(i)}, \boldsymbol{\sigma}^{(j)}\right\rangle
    \end{array}
  \end{displaymath}
  where the last step is due to (\ref{eq:algebraicStabilityMatrix}). 
  The second condition in Definition \ref{def:RKAlgebraicStability}
  gives 
  \begin{displaymath}
    \begin{array}{l}
    (*):\quad
    \left\langle \mathbf{u}^{n+1},\mathbf{u}^{n+1}\right\rangle \le
    \left\langle \mathbf{u}^n,\mathbf{u}^n\right\rangle + 2k\sum_{i=1}^{s}b_i\left\langle \mathbf{u}^{(i)},\boldsymbol{\sigma}^{(i)}\right\rangle.
    \end{array}
  \end{displaymath}

  Write $\alpha_i := I_{cv}(\widetilde{\mathbf{u}}^{(i)}, \mathbf{u}^{(i)})
  = \left\langle \widetilde{\mathbf{u}}^{(i)}\cdot\nabla
    \widetilde{\mathbf{u}}^{(i)},  \mathbf{u}^{(i)}\right\rangle$
  for each $i=1, \ldots, s$. 
  Then 
  \begin{displaymath}
    \begin{array}{rl}
      &\innerProd{\boldsymbol{\sigma}^{(i)}, \mathbf{u}^{(i)}}
      = \innerProd{\mathscr{P}\boldsymbol{\rho}^{(i)}, \mathbf{u}^{(i)}}
      = \innerProd{\boldsymbol{\rho}^{(i)}, \mathbf{u}^{(i)}}
     \\
      =& \innerProd{ \mathbf{g}^{(i)}-\nabla q^{(i)}, \mathbf{u}^{(i)}}
      -r^{(i)}\alpha_i
        + \nu\innerProd{\Delta \mathbf{w}^{(i)}, \mathbf{u}^{(i)}}
        = -r^{(i)}\alpha_i -\nu\left\|\nabla \mathbf{u}^{(i)}\right\|^2,
    \end{array}
  \end{displaymath}
  where the first equality follows from
  $\boldsymbol{\sigma}^{(i)}= \mathscr{P}\boldsymbol{\rho}^{(i)}$,
  the second from (\ref{eq:cProjLH}), 
  (\ref{eq:GePUPeSAVRKIntermediate}g,h),
  and Lemma \ref{lem:orthogonalityRelation}, 
  the third from (\ref{eq:GePUPeSAVRKIntermediate}c)
  and the definition of $\alpha_i$,
  and the last from 
  $\mathbf{g}$ being conservative, 
  Lemma \ref{lem:orthogonalityRelation}, 
  (\ref{eq:innerProdLapWandU}), 
  and Corollary \ref{coro:divFreeSemiImplicit}.
  
  Since the Leray-Helmholtz projection $\mathscr{P}$
  has no control over the tangential velocity,
  in (\ref{eq:GePUPeSAVRKIntermediate}g)
  we would have $\boldsymbol{\tau}\cdot\mathbf{u}^{(i)}\ne 0$
  if $\nabla\cdot\mathbf{w}^{(i)}\ne 0$.
  Fortunately Corollary \ref{coro:divFreeSemiImplicit} dictates
  $\nabla\cdot\mathbf{w}^{(i)}= 0$
  and $\mathbf{n}\cdot\mathbf{w}^{(i)}|_{\partial\Omega}=0$,
  then in (\ref{eq:GePUPeSAVRKIntermediate}g) 
  $\mathscr{P}$ reduces to the identity. 
  Thus the second integral in the last line of
  (\ref{eq:innerProdLapWandU}) vanishes in GePUP-ES-RK.
  A related observation is that, 
  although $\nabla q$ is orthogonal to $\mathbf{w}$,
  it cannot be arbitrary as it must make $\boldsymbol{\rho}^{(j)}$
  in (\ref{eq:GePUPeSAVRKIntermediate}d) divergence-free;
  otherwise 
  in (\ref{eq:GePUPeSAVRKIntermediate}g) 
  $\mathscr{P}$ would not reduce to the identity. 
  
  Substitute 
  $\left\langle \mathbf{u}^{(i)}, \boldsymbol{\sigma}^{(i)}\right\rangle
  = -r^{(i)} \alpha_i-\nu\left\|\nabla \mathbf{u}^{(i)}\right\|^2$
  into $(*)$ and we have 
  \begin{displaymath}
    \begin{array}{l}
      \left\langle \mathbf{u}^{n+1}, \mathbf{u}^{n+1}\right\rangle
      \le \left\langle \mathbf{u}^n, \mathbf{u}^n \right\rangle-2k\sum_{i=1}^{s}b_i\alpha_ir^{(i)}-2\nu k\sum_{i=1}^{s}b_i\left\|\nabla \mathbf{u}^{(i)}\right\|^2.
    \end{array}
  \end{displaymath}
  
  Similarly, the positive semi-definiteness of $M$
  in Definition \ref{def:RKAlgebraicStability} yields 
  \begin{displaymath}
    \begin{array}{rl}
      & \left|r^{n+1}\right|^2
        = \left(r^n+k\sum_{i=1}^{s}b_i\alpha_i\right)^2
      \\
      =&
         \left|r^n\right|^2+k\sum_{i=1}^{s}b_i\alpha_ir^n+k\sum_{j=1}^{s}b_j\alpha_jr^n
      +k^2\sum_{i=1}^{s}\sum_{j=1}^{s}b_ib_j\alpha_i\alpha_j
      \\
      =& \left|r^n\right|^2
         +
         k\sum_{i=1}^{s}b_i\alpha_i\left(r^{(i)}-k\sum_{j=1}^{s}a_{i,j}\alpha_j\right)
      \\
      &+k\sum_{j=1}^{s}b_j\alpha_j\left(r^{(j)}-k\sum_{i=1}^{s}a_{j,i}\alpha_i\right)
      +k^2\sum_{i=1}^{s}\sum_{j=1}^{s}b_ib_j\alpha_i\alpha_j
      \\
      =& \left|r^n\right|^2+2k\sum_{i=1}^{s}b_i\alpha_ir^{(i)}-k^2\sum_{i=1}^{s}\sum_{j=1}^{s}m_{i,j}\alpha_i\alpha_j
      \\
      \le& \left|r^n\right|^2+2k\sum_{i=1}^{s}b_i\alpha_ir^{(i)}.
    \end{array}
  \end{displaymath}
  The proof is completed by summing up the above two inequalities.
\qed 
\end{proof}

In the last step of the above proof, 
 it is the auxiliary variable $r$ in the GePUP-ES formulation that 
 leads to the cancellation of $2k\sum_{i=1}^{s}b_i\alpha_ir^{(i)}$.

Except on staggered grids,
 the initial condition of $\mathbf{w}$ being incompressible
 in (\ref{eq:GePUP-initialW})
 cannot be exactly fulfilled in practical computations. 
However, the mechanism of divergence decay ensures that
 $\mathbf{w}$ converge to $\mathbf{u}$ sufficiently fast;
 see Theorem \ref{coro:divFreeSemiImplicit}.
This is also related to the millennium problem
 on the well-posedness of the INSE.
If the solution of the INSE blows up,
 then its divergence must blow up first.
Therefore,
 our computation only works
 when the INSE admits a bounded solution.
After all,
 one cannot expect to solve the millennium problem
 by reformulating INSE.

\subsection{Semi-implicit GePUP-ES-SDIRK schemes}
\label{sec:GePUP-SAV-SDIRK}

In light of Theorems \ref{thm:decayOfVelDiv-GePUP-SAV-RK} and
 \ref{thm:GePUPSAVRKEnergyDecay},
 one way to preserve the monotonic decrease of the kinetic energy
 and the exponential decay of the divergence residue
 is to employ an algebraically stable RK method
 in GePUP-ES-RK.
Gauss-Legendre RK methods are algebraically stable
 and have a minimal number of stages
 for a given temporal order of accuracy.
However, their employment in GePUP-ES-RK
 necessitates either the coupling of all intermediate stage values
 of $\mathbf{w}$
 or the use of complex arithmetic.
Thus we turn to singly diagonal implicit RK (SDIRK) methods
 that satisfy
 \begin{equation}
   \label{eq:SDIRK-conditions}
   a_{i,j} =
   \begin{cases}
     0 & \text{if } i<j;
     \\
     \gamma\ne 0 & \text{if } i=j, 
   \end{cases}
 \end{equation}
 aiming to design a family of GePUP-ES-SDIRK schemes that only consist of
 solving a sequence of \emph{linear} BVPs with real arithmetic,
 one intermediate stage at a time.
The core difficulty for this,
 as mentioned in Section \ref{sec:GePUP-ES},
 is the nonlinear tight coupling of $\mathbf{u}$, $\mathbf{w}$, $q$
 and $r$.
Our solution is

\begin{defn}
  \label{def:GePUP-ES-SDIRK}
  A \emph{GePUP-ES-SDIRK (GES) scheme}
   is a GePUP-ES-RK scheme where
   an algebraically stable SDIRK method is employed
   as the RK method
   and stage values for each intermediate stage
   $i=1, \ldots, s$ are decomposed as
   \begin{equation}
     \label{eq:GePUPSAVSDIRKDecomp}
     \left\{
       \begin{array}{rl}
         \mathbf{w}^{\left( i \right)} &= \mathbf{w}_0^{(i)} + \sum_{j=1}^ir^{(j)} \mathbf{w}_j^{(i)}, \\
         \mathbf{u}^{(i)} &= \mathbf{u}_0^{(i)} + \sum_{j=1}^ir^{(j)} \mathbf{u}_j^{(i)}, \\
         q^{(i)} &= q_1^{(i)} + r^{(i)}q_2^{(i)}, 
       \end{array}\right.
   \end{equation}
   where 
   $\mathbf{w}_0^{(i)}, \mathbf{w}_1^{(i)}, \ldots, \mathbf{w}_i^{(i)}$, 
   $\mathbf{u}_0^{(i)}, \mathbf{u}_1^{(i)}, \ldots, \mathbf{u}_i^{(i)}$,
   and $q_1^{(i)}, q_2^{(i)}$
   are auxiliary variables for the $i$th stage,
   which consists of steps as follows, 
   {\setlist[enumerate]{label=({GES}.\arabic{enumi})}
     \begin{enumerate}
     \item solve for $q_1^{(i)}$ from the Neumann BVP
       \begin{displaymath}
         \left\{
           \begin{array}{rll}
             \Delta q_{1}^{(i)}
             &= \nabla\cdot\mathbf{g}^{(i)}  &\text{in } \Omega,
             \\
             \mathbf{n}\cdot\nabla q_{1}^{(i)}
             &= \mathbf{n}\cdot \left( \mathbf{g}^{(i)} + \nu \Delta
               \widetilde{\mathbf{w}}^{(i)} \right) +
               \lambda\mathbf{n}\cdot\widetilde{\mathbf{w}}^{(i)}  \quad
                                             &\text{on } \partial\Omega, 
           \end{array}\right.
       \end{displaymath}
       where $\widetilde{\mathbf{w}}^{(i)}$
       is an approximation of $\mathbf{w}^{(i)}$ 
       obtained by (\ref{eq:interpolationEAV}) and (EAV-1,2).
     \item solve for $q_2^{(i)}$ from the Neumann BVP
       \begin{displaymath}
         \left\{
           \begin{array}{rll}
             \Delta q_{2}^{(i)}
             &= -\nabla\cdot \left(\widetilde{\mathbf{u}}^{(i)}
               \cdot\nabla \widetilde{\mathbf{u}}^{(i)}\right)
             & \text{in } \Omega,
             \\
             \mathbf{n}\cdot\nabla q_2^{(i)}
             &= -\mathbf{n}\cdot \left(\widetilde{\mathbf{u}}^{(i)}
               \cdot\nabla \widetilde{\mathbf{u}}^{(i)}\right)
             & \text{on } \partial\Omega.
           \end{array}\right.
       \end{displaymath}

     \item 
       solve for $\mathbf{w}_{0}^{(i)}$ and 
       $\mathbf{w}_{\ell}^{(i)}$ where $\ell=1, \ldots, i$
       from BVPs with their boundary conditions
       in (\ref{eq:GePUPeSAVRKIntermediate}e).
       \begin{displaymath}
         \left\{
           \begin{array}{rl}
             (1 - \nu\gamma k\Delta)\mathbf{w}_0^{(i)}
             
             &= \mathbf{u}^n
               + k\sum\limits_{j=1}^{i}a_{i,j}\left( \mathbf{g}^{(j)} - \nabla
               q_1^{(j)} \right)
               + \nu k\sum\limits_{j=1}^{i-1}a_{i,j}\Delta
               \mathbf{w}_0^{(j)}; 
             \\
             (1-\nu\gamma k\Delta)\mathbf{w}_{\ell}^{(i)}
             &= -ka_{i,\ell} \left(
               \widetilde{\mathbf{u}}^{(\ell)}\cdot\nabla\widetilde{\mathbf{u}}^{(\ell)}
               + \nabla q_2^{(\ell)} \right) + \nu k\sum\limits_{j=\ell}^{i-1}a_{i,j}\Delta \mathbf{w}_{\ell}^{(j)}.
           \end{array} \right.
       \end{displaymath}
     \item compute $\mathbf{u}_{\ell}^{(i)}=\cProjLH\mathbf{w}_{\ell}^{(i)}$
       for $\ell=0, 1, \ldots, i$.
     \item define $\overline{\mathbf{u}}^{(i)}:=
       \mathbf{u}_0^{(i)} +
       \sum_{\ell=1}^{i-1}r^{(\ell)}\mathbf{u}_{\ell}^{(i)}$
       and calculate $r^{(i)}$ by
       \begin{displaymath}
         \left( 1 - \gamma k
           I_{cv}\left(\widetilde{\mathbf{u}}^{(i)}, \mathbf{u}^{(i)}_i\right)
         \right)r^{(i)}
         \! = \! r^n + \gamma k
         I_{cv}\left(\widetilde{\mathbf{u}}^{(i)}\!,\!
           \overline{\mathbf{u}}^{(i)}\right)
         + k\sum_{j=1}^{i-1}a_{i,j}
         I_{cv}\left(\widetilde{\mathbf{u}}^{(j)}\!,\!
           {\mathbf{u}}^{(j)}\right).
       \end{displaymath}
     \item calculate $\mathbf{w}^{(i)}$ and $\mathbf{u}^{(i)}$ by (\ref{eq:GePUPSAVSDIRKDecomp}).
     \end{enumerate}}
\end{defn}
 
The above steps (GES.1--6)
 are direct consequences of the GePUP-ES-RK scheme
 (\ref{eq:GePUPSAVRKFinal}), 
 the property (\ref{eq:SDIRK-conditions}) of SDIRK, 
 and the decomposition (\ref{eq:GePUPSAVSDIRKDecomp}). 

Substitute the third decomposition in (\ref{eq:GePUPSAVSDIRKDecomp})
 into (\ref{eq:GePUPeSAVRKIntermediate}a,b),
 separate the terms with and without $r^{(i)}$,
 notice that $r^{(i)}$ is a scalar, and 
 we have (GES.1--2).
In step (GES.1),
 we decouple $q^{(i)}$ from $\mathbf{w}^{(i)}$
 by approximating $\mathbf{w}^{(i)}$ with
 $\widetilde{\mathbf{w}}^{(i)}$
 using the same method for calculating $\widetilde{\mathbf{u}}^{(i)}$.

For SDIRK,
 the upper bound of the summation
 in (\ref{eq:GePUPeSAVRKIntermediate}d) is $i$.
Substitute (\ref{eq:GePUPSAVSDIRKDecomp}) into
 (\ref{eq:GePUPeSAVRKIntermediate}c,d),
 separate the terms with and without $r^{(i)}$,
 switch the summation order by 
 $\sum_{j=1}^ia_{i,j}\sum_{\ell=1}^jr^{(\ell)}\Delta \mathbf{w}_{\ell}^{(j)} =
 \sum_{\ell=1}^ir^{(\ell)}\sum_{j=\ell}^ia_{i,j}\Delta
 \mathbf{w}_{\ell}^{(j)}$, 
 and we have (GES.3).
Then (GES.4) follows from (\ref{eq:GePUPeSAVRKIntermediate}g).

(GES.5) follows from
 substituting the second equation in (\ref{eq:GePUPSAVSDIRKDecomp})
 into (\ref{eq:GePUPeSAVRKIntermediate}f), 
 \begin{displaymath}
   \begin{array}{l}
   r^{(i)} = r^n + \gamma k
   I_{cv}\left(\widetilde{\mathbf{u}}^{(i)},\
     \overline{\mathbf{u}}^{(i)}+r^{(i)}\mathbf{u}_i^{(i)}\right)
   + k\sum_{j=1}^{i-1}a_{i,j}
   I_{cv}\left(\widetilde{\mathbf{u}}^{(j)},\
     {\mathbf{u}}^{(j)}\right), 
   \end{array}
 \end{displaymath}
 and moving the $r^{(i)}$-term on the RHS to the LHS.
Finally, the parenthesis on the LHS of (GES.5)
 must be no less than 1 because of
 (\ref{eq:innerProdLapWandU}),
 Lemma \ref{lem:orthogonalityRelation}, 
 and the inner product of $\mathbf{u}_i^{(i)}$ to the equation
 on $\mathbf{w}_i^{(i)}$ in (GES.3).
 
(\ref{eq:GePUPSAVRKFinal}) and (GES.1--6)
 are the complete algorithmic steps of the GES scheme.
Theorems \ref{thm:decayOfVelDiv-GePUP-SAV-RK}
 and \ref{thm:GePUPSAVRKEnergyDecay} state that
 both the velocity divergence and the modified energy
 in GES decrease monotonically
 provided that the employed SDIRK is algebraically stable; 
 an example is given in (\ref{eq:SDIRK43}). 

In this work,
we discretize the continuous spatial operators
 in Definition \ref{def:GePUP-ES-SDIRK}
 by the fourth-order collocated finite-volume operators
 in \cite[Sec. 3 \& 4]{zhang16:_GePUP}
 to obtain a fully discrete GES scheme.
However,
it is emphasized that
what we have proposed
 is not a single scheme but a space of solvers, 
 each of which can be easily constructed
 by making menu choices for `orthogonal' policies
 that span the solver space.
The orthogonal structure of this solver space 
 is very conducive to reusing \emph{in a black-box manner}
 the legacy of classical finite volume/difference methods
 and the wealth of theory and algorithms
 for numerically solving ODEs;
 see Table \ref{tab:GePUP-SAV-SDIRK-policies}.  
It is the GePUP-E formulation 
 that makes this black-box reuse possible.

It would be ideal if the conclusions of
 Theorems \ref{thm:decayOfVelDiv-GePUP-SAV-RK} and
 \ref{thm:GePUPSAVRKEnergyDecay}
 could also hold in the fully discrete case.
Unfortunately,
 for the finite volume discretization on collocated grids,
 the fully discrete counterparts of
 Theorems \ref{thm:decayOfVelDiv-GePUP-SAV-RK} and
 \ref{thm:GePUPSAVRKEnergyDecay}
 only hold \emph{asymptotically},
 i.e.,
 in the limit of $k$ and $h$ simultaneously approaching zero.
We defer to future research 
 the investigation of suitable spatial discretizations so that
 Theorems \ref{thm:decayOfVelDiv-GePUP-SAV-RK} and
 \ref{thm:GePUPSAVRKEnergyDecay}
 also hold in the fully discrete case.


 

%% file: tikz/interpolationEAU.tex
\subfigure[(EAV-1) for $n=0,1$]{
  \begin{tikzpicture}[scale = 1]
    \tikzstyle{every node}=[scale=1]
    \draw (0, 0)--(8, 0);
    \draw (0, 2)--(8, 2);
    \draw (0, 2)--(0, 2.2);
    \draw (2.67, 2)--(2.67, 2.2);
    \draw (5.34, 2)--(5.34, 2.2);
    \draw (8, 2)--(8, 2.2);
    \draw (0, 0)--(0, 0.4);
    \draw (8, 0)--(8, 0.4);
    \draw (2, 0)--(2, 0.2);
    \draw (4, 0)--(4, 0.2);
    \draw (6, 0)--(6, 0.2);
    \node [below] (1) at (0, 0) {$t^{n}$};
    \node [below] (2) at (8, 0) {$t^{n + 1}$};
    \node [above] (4) at (2, 0.2) {$\widetilde{\mathbf{u}}_{n}^{\left( 1 \right)}$};
    \node [above] (5) at (4, 0.2) {$\widetilde{\mathbf{u}}_{n}^{\left( 2 \right)}$};
    \node [above] (6) at (6, 0.2) {$\widetilde{\mathbf{u}}_n^{\left( 3 \right)}$};
    \node [below] (7) at (2, 0) {$t^{n} + c_{1}k$};
    \node [below] (8) at (4, 0) {$t^{n} + c_{2}k$};
    \node [below] (9) at (6, 0) {$t^{n} + c_{3}k$};
    \node [below] (10) at (0, 2) {$t^{n}$};
    \node [below] (11) at (2.67, 2) {$t^{n} + \frac{1}{3} k$};
    \node [below] (12) at (8, 2) {$t^{n} + \frac{3}{3} k = t^{n + 1}$};
    \node [above] (13) at (0, 0.4) {$\mathbf{u}^{n}$};
    \node [above] (14) at (8, 0.4) {$\mathbf{u}^{n + 1}$};
    \node [above] (15) at (0, 2.2) {$\hat{\mathbf{u}}_n^{\left( 0 \right)} = \mathbf{u}^{n}$};
    \node [above] (16) at (2.67, 2.2) {$\hat{\mathbf{u}}_n^{\left( 1 \right)}$};
    \node [above] (17) at (8, 2.2) {$\hat{\mathbf{u}}_n^{\left( 3 \right)}$};
    \node [below] (18) at (5.34, 2) {$t^{n} + \frac{2}{3} k$};
    \node [above] (16) at (5.34, 2.2) {$\hat{\mathbf{u}}_n^{\left( 2 \right)}$};
  \end{tikzpicture}
}

\subfigure[(EAV-2) for $n\ge 2$]{
  \begin{tikzpicture}[scale = 0.7]
    \tikzstyle{every node}=[scale=1.0]
    \draw (0, 0)--(8, 0);
    \draw (-6, 2)--(8, 2);
    \draw (0, 2)--(0, 2.2);
    \draw (-6, 2)--(-6, 2.2);
    \draw (-3, 2)--(-3, 2.2);
    \draw (8, 2)--(8, 2.2);
    \draw (0, 0)--(0, 0.4);
    \draw (8, 0)--(8, 0.4);
    \draw (2, 0)--(2, 0.2);
    \draw (4, 0)--(4, 0.2);
    \draw (6, 0)--(6, 0.2);
    \node [below] (1) at (0, 0) {$t^{n}$};
    \node [below] (2) at (8, 0) {$t^{n + 1}$};
    \node [below] (3) at (-6, 2) {$t^{n - 2}$};
    \node [above] (4) at (2, 0.2) {$\widetilde{\mathbf{u}}_{n}^{\left( 1 \right)}$};
    \node [above] (5) at (4, 0.2) {$\widetilde{\mathbf{u}}_{n}^{\left( 2 \right)}$};
    \node [above] (6) at (6, 0.2) {$\widetilde{\mathbf{u}}_n^{\left( 3 \right)}$};
    \node [below] (7) at (2, 0) {$t^{n} + c_{1}k$};
    \node [below] (8) at (4, 0) {$t^{n} + c_{2}k$};
    \node [below] (9) at (6, 0) {$t^{n} + c_{3}k$};
    \node [below] (10) at (0, 2) {$t^{n}$};
    \node [below] (11) at (-3, 2) {$t^{n - 1}$};
    \node [below] (12) at (8, 2) {$t^{n + 1} = t^n + k$};
    \node [above] (13) at (0, 0.4) {$\mathbf{u}^{n}$};
    \node [above] (14) at (8, 0.4) {$\mathbf{u}^{n + 1}$};
    \node [above] (15) at (0, 2.2) {$\hat{\mathbf{u}}_n^{\left( 2 \right)} = \mathbf{u}^{n}$};
    \node [above] (16) at (-3, 2.2) {$\hat{\mathbf{u}}_n^{\left( 1 \right)} = \mathbf{u}^{n - 1}$ };
    \node [above] (17) at (8, 2.2) {$\hat{\mathbf{u}}_n^{\left( 3 \right)}$};
    \node [above] (18) at (-6, 2.2) {$\hat{\mathbf{u}}_n^{\left( 0 \right)} = \mathbf{u}^{n - 2}$};
  \end{tikzpicture}
}


%% file: sec/tests.tex
\section{Tests}
\label{sec:tests}

In this section,
 we test a fully discrete GES scheme
 by several numerical experiments
 to confirm the analytic results in previous sections.
We employ the fourth-order finite-volume discretizations
 in \cite{zhang16:_GePUP}
 and a fourth-order, algebraically stable SDIRK method
 proposed by Du, Ju \& Lu \cite{du19:sdirk}, 
 \begin{equation}
   \label{eq:SDIRK43}   \renewcommand{\arraystretch}{1.5}
   \begin{tabular}[h]{c|ccc}
     $\gamma$ & $\gamma$ & & \\
     $\frac{1}{2}$ & $\frac{1}{2} - \gamma$ & $\gamma$ & \\
     $1-\gamma$ & $2\gamma$ & $1-4\gamma$ & $\gamma$ \\
     \hline
              & $\mu$ & $1-2\mu$ & $\mu$
   \end{tabular}, \qquad
   \left\{
     \begin{array}{l}
       \gamma = \frac{1}{\sqrt{3}}\cos \frac{\pi}{18} + \frac{1}{2},
       \\
       \mu = \frac{1}{6(2\gamma-1)^2}.
     \end{array}\right.
 \end{equation}

   \begin{table}
     \centering
     \caption{Choices of the main orthogonal policies
       that constitute the tested GES scheme
       as a fully discrete GePUP-ES-RK. 
     }
     \label{tab:GePUP-SAV-SDIRK-policies}
     \small
     \scalebox{0.96}{
     \begin{tabular}{c|c|c}
       \hline
       Orthogonal policies
       & For the tested GES scheme
       & Other possible options
       \\ \hline\hline
       time integration
       & the SDIRK method (\ref{eq:SDIRK43})
       & algebraically stable RK methods
       \\ \hline
       temporal accuracy
       & the fourth order
       & the second, third, and fifth orders
       \\ \hline
       spatial discretization
       & finite volume
       & finite difference
       \\ \hline
       spatial accuracy
       & the fourth order
       & the second order
       \\ \hline
       variable location
       & collocated (as in \cite{zhang16:_GePUP})
       & staggered 
       \\ \hline
       estimating $\widetilde{\mathbf{u}}^{(i)}$
       and $\widetilde{\mathbf{w}}^{(i)}$
       &
       & interpolate 
       \\
       in Definitions \ref{def:GePUP-ES-RK}
       and \ref{def:GePUP-ES-SDIRK}
       & (\ref{eq:interpolationEAV}) \& (EAV-1,2)
       & both in time and in space
       \\ \hline
     \end{tabular}
     }
   \end{table}
 
Along with other possible options of the major orthogonal policies,
 we show
 in Table \ref{tab:GePUP-SAV-SDIRK-policies}
 our choices that constitute the particular
 GES scheme to be tested. 
According to the analysis in Section \ref{sec:algorithms}, 
 this particular INSE solver 
 should be fourth-order accurate both in time and in space,
 preserve monotonic decrease of the total kinetic energy, 
 and have the velocity divergence well under control.
 
All numerical tests are performed
 on a rectangular domain in Cartesian coordinates, 
 where the electric boundary conditions
 in (\ref{eq:GePUPe}b)
 are enforced by homogeneous Dirichlet conditions 
 for the tangential velocity components
 and a homogeneous Neumann condition for the normal velocity.
The initial cell-averaged velocity 
 is calculated by Boole's rule, 
 a sixth-order formula of Newton-Cotes quadrature. 
Since exact solutions are unavailable,
 we define the computational error
 via Richardson extrapolation,
 i.e., by the difference of the solution
 on the current grid and
 that of the next finer grid.


For different values of the penalty parameter $\lambda=$1, 10, and 100, 
 the corresponding error norms in all numerical tests
 remain the same up to the first two significant digits; 
 these results differ from those reported in
 \cite{shirokoff11:_navier_stokes,rosales21:_high_poiss_navier_stokes}.
However, this is not surprising since
 projecting $\mathbf{w}$
 to $\mathbf{u}$ with $\mathbf{u}\cdot\mathbf{n}=0$
 in (\ref{eq:GePUPe}c,d)
 and setting 
$\mathbf{w}^0=\mathbf{u}^0$
 in (\ref{eq:GePUP-initialW})
 already imply a small magnitude of $\mathbf{n}\cdot\mathbf{w}$. 
Nonetheless, a positive-valued $\lambda$ 
 guarantees the no-penetration condition be fulfilled.
Hereafter 
 we only show results in the case of $\lambda=1$.

\subsection{Single-vortex tests}
\label{sec:single-vortex-test}

 \begin{figure}
   \centering
   \subfigure[$t$=40]{
     \includegraphics[width=0.8\textwidth]{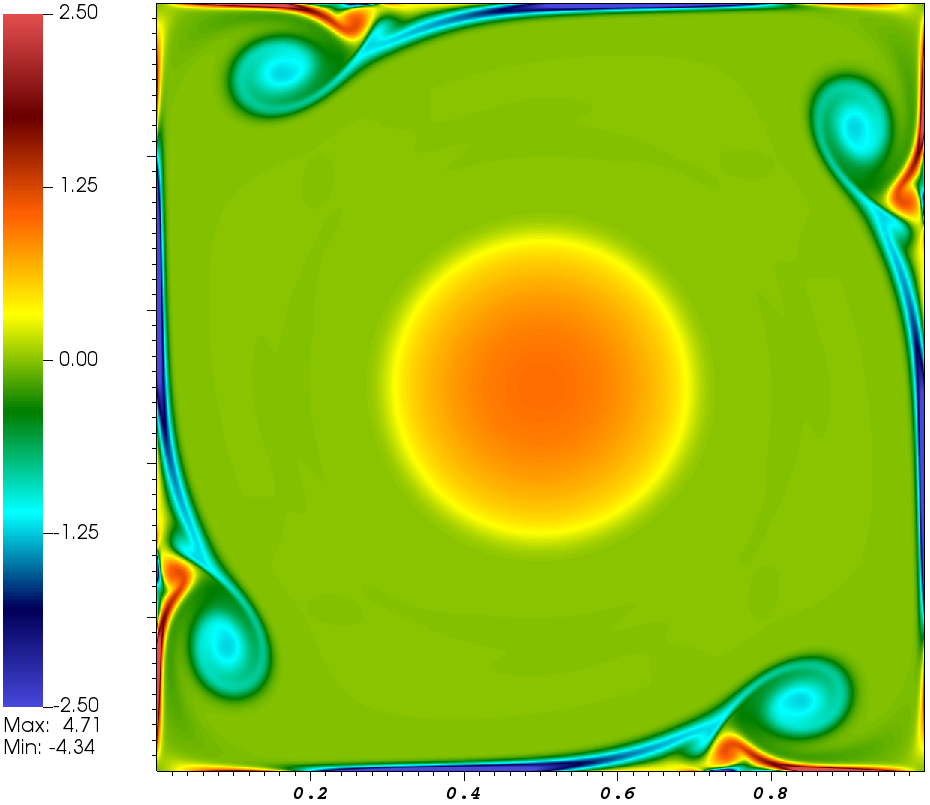}
   }
   \subfigure[$t$=60]{
     \includegraphics[width=0.8\textwidth]{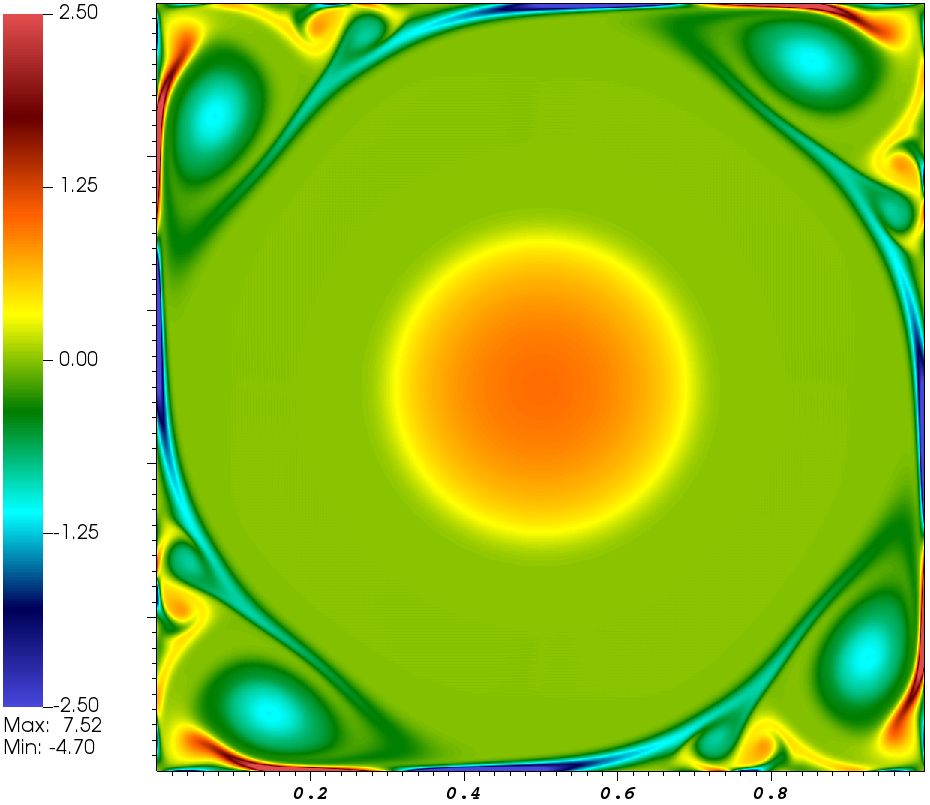}
   }
   \caption{Snapshots of vorticity for the single-vortex test
     with $\mathrm{Re}=20,000$ and $\lambda=1$
     on a uniform grid with $h=\frac{1}{1024}$ and $\mathrm{Cr}=0.5$.
     The region of each cell is filled by a single color
      that corresponds to the cell-averaged vorticity.
      No image smoothing is applied.
   }
   \label{fig:singleVortexBox}
 \end{figure}

Following
\cite{bch91:_second_order_projection_for_viscous_incompressible_flow}, 
 we define an axisymmetric velocity field
 on $\Omega=[0,1]^2$ by
\begin{equation}
  \label{eq:singleVortexBox}
  \renewcommand{\arraystretch}{1.2}
  u_{\theta}(r_v) =
  \left\{
    \begin{array}{cc}
      \Gamma(\frac{1}{2}r_v - 4r_v^3) & \text{ if } r_v<R ;
      \\
      \Gamma\frac{R}{r_v}(\frac{1}{2}R - 4R^3) & \text{ if } r_v\ge R,
    \end{array}
  \right.
\end{equation}
where $r_v$ is the distance from the domain center
 $(\frac{1}{2}, \frac{1}{2})^T$.
The choices $R=0.2$ and $\Gamma=1$
 give $\max(u_{\theta}) = 0.068$.
A small viscosity $\nu = 3.4\times 10^{-6}$
 yields a high Reynolds number $\mathrm{Re}=20,000$.
The initial velocity is obtained
 by projecting cell averages of $\mathbf{u}$ in (\ref{eq:singleVortexBox})
 ten times to make it approximately divergence-free.

 \begin{table}
   \caption{
     Errors and convergence rates of the GES scheme
     in Table \ref{tab:GePUP-SAV-SDIRK-policies}
     for solving the single-vortex test
     with \mbox{Re = 20,000},
      $t_0=0.0$, $t_e=60$, $\lambda = 1$ 
      and the Courant number $\mathrm{Cr}=0.5$.
    }
   \centering
   \input{tab/singleVortexBox_Gamma1_Re2e4}
   \label{tab:singleVortexBox}
 \end{table}

 \begin{figure}
   \centering
   \subfigure[$\mathcal{E}_h^n$ as in (\ref{eq:modifiedEnergyDiscrete})]
   {
     \includegraphics[width=0.475\textwidth]{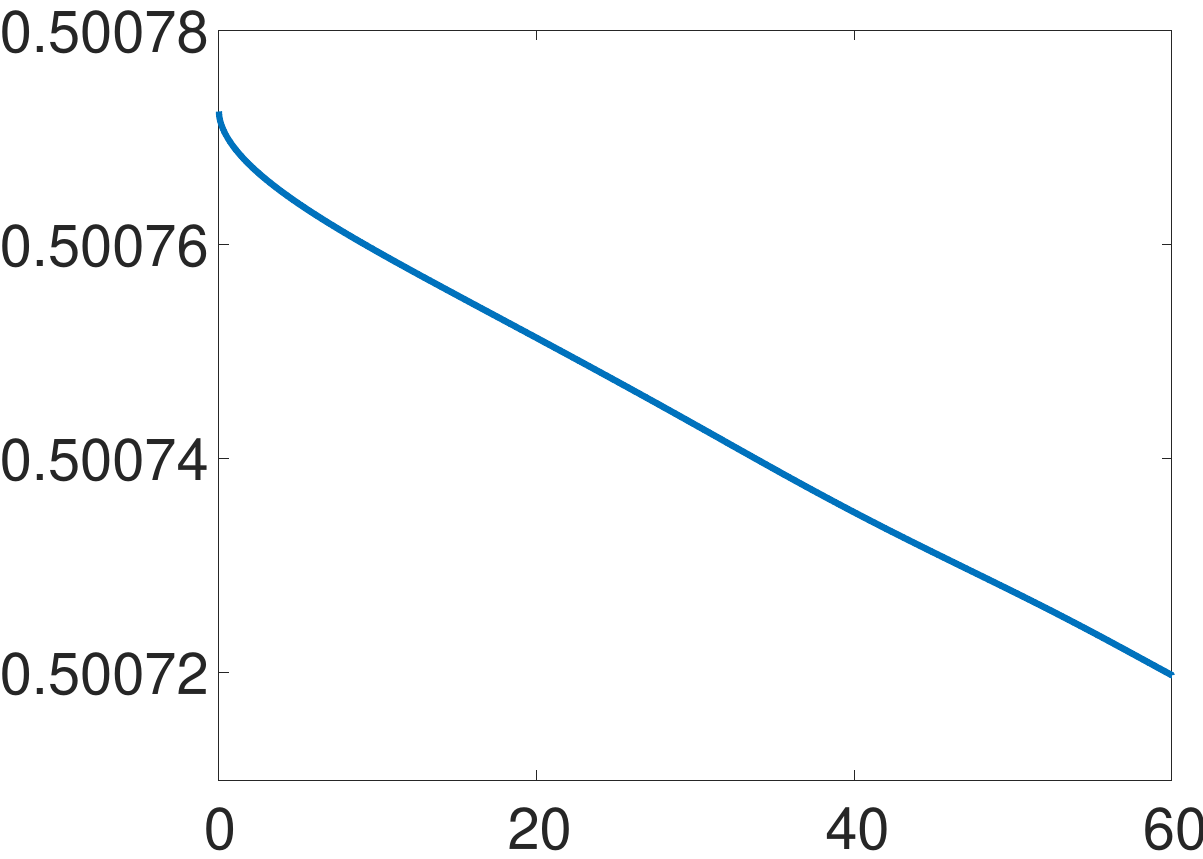}
   }
   \hfill
   \subfigure[$|r_h^n-1|$] 
   {
     \includegraphics[width=0.475\textwidth]{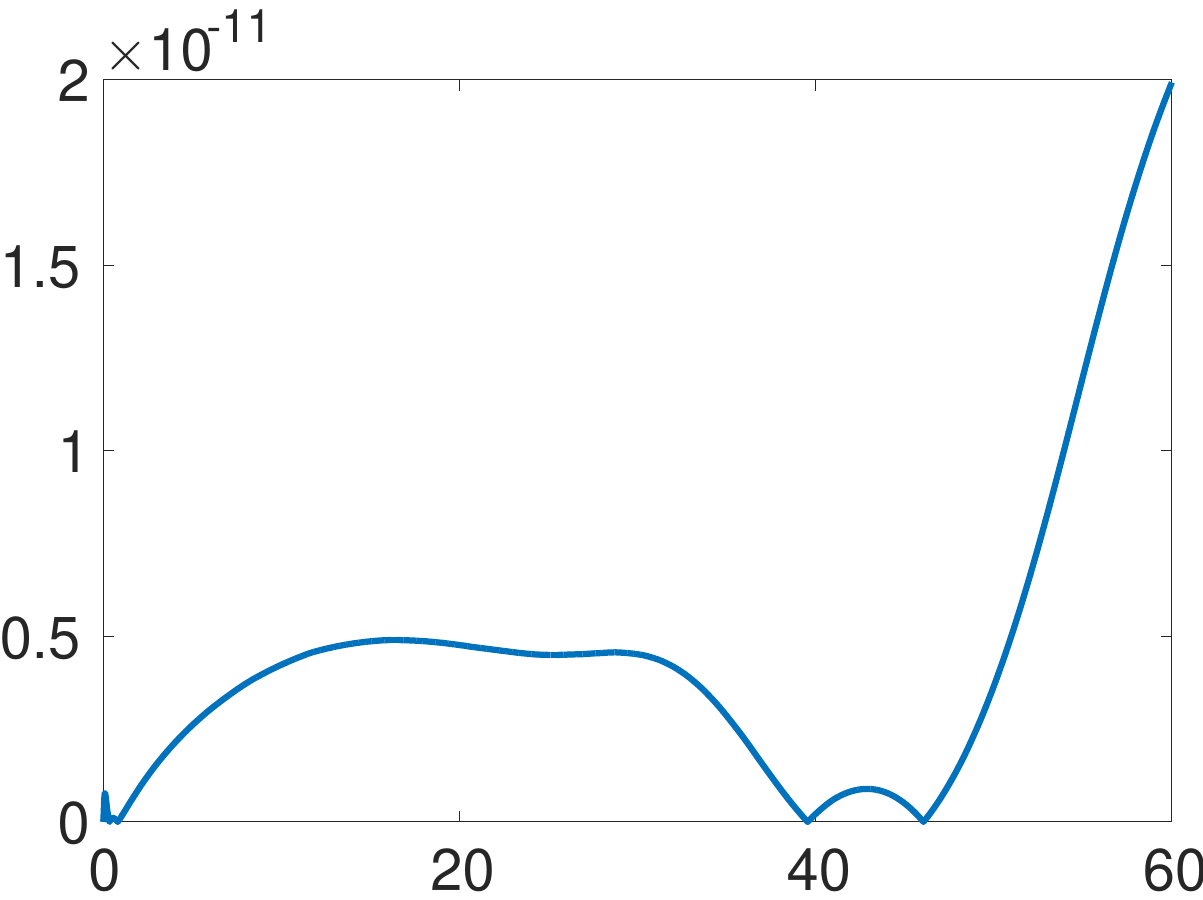}
   }

   \subfigure[$\left\|\mathbf{D} \mathbf{W}_h^{n}\right\|_{\mathcal C}$
   as in (\ref{eq:L2norm4FV})]
   {
     \includegraphics[width=0.475\textwidth]{
       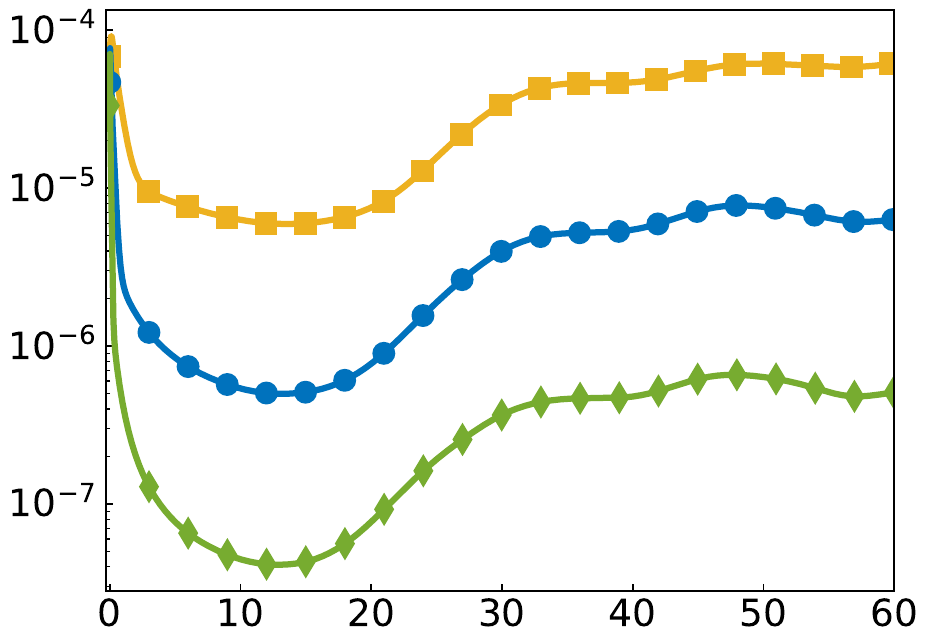}
   }
   \hfill
   \subfigure[Zoom-in of (c)]
   {
     \includegraphics[width=0.475\textwidth]{
       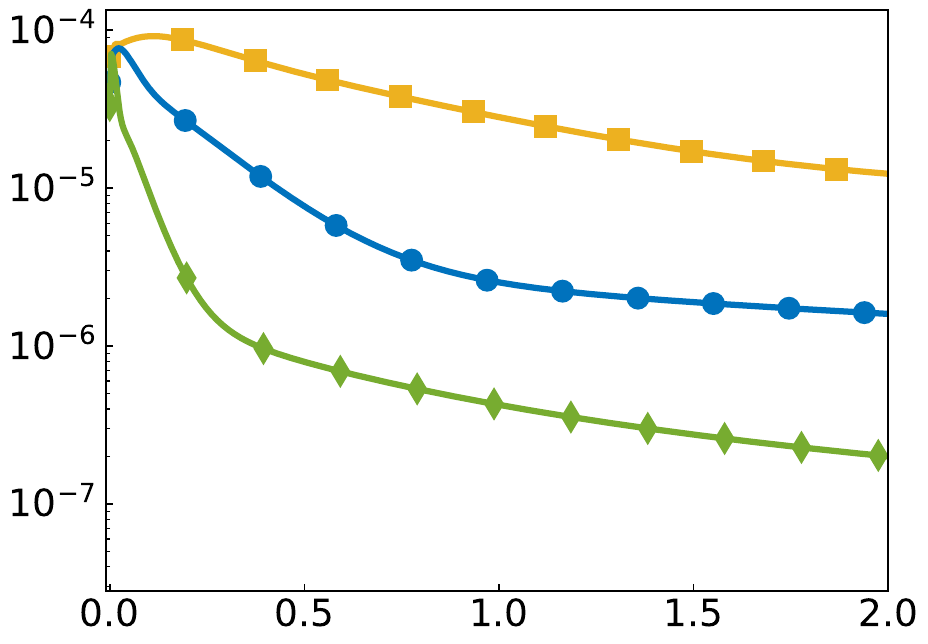}
   }
   \caption{Evolutions of the modified energy,
     the SAV, and the $L_2$-norm of $\mathbf{D} \mathbf{W}_h^{n}$
     produced by the GES scheme in Table
     \ref{tab:GePUP-SAV-SDIRK-policies}
     for solving the single-vortex test 
     with $\mathrm{Re}=2\times 10^4$, $\mathrm{Cr}=0.5$, 
     and $\lambda=1$. 
     The abscissa in all subplots is time.
     In (a) and (b), $h=\frac{1}{1024}$
     and $k\approx 7.18\times 10^{-3}$. 
     In (c) and (d),
     the curves marked by ``{\tiny $\blacksquare$},''
     ``$\bullet$,'' and ``$\blacklozenge$''
     represent the results for $h=\frac{1}{512}$, $\frac{1}{1024}$, and
     $\frac{1}{2048}$, respectively.
   }
   \label{fig:singleVortexBox-energyAndDiv}
   \vspace{-10pt}
 \end{figure}
 
The tests are performed on four successively refined grids
 with uniform grid size $h$.
The time span $[0, 60]$ is made long enough
 for the turbulent boundary layers
 to develop prominent Lagrangian coherent structures.
The uniform time step size $k$ is given by
 setting the Courant number
 $\mathrm{Cr}:=\frac{k}{h}\max(u_{\theta})$
 at 0.5. 
Snapshots of the vorticity
 at time $t = 40$ and at the final time $t = 60$
 are shown in Figure \ref{fig:singleVortexBox}, 
 where
 the essential features of vortex sheet roll-up and counter-vortices
 agree with those in 
 \cite{bch91:_second_order_projection_for_viscous_incompressible_flow}.

It is clear in Table \ref{tab:singleVortexBox}
 that convergence rates of the velocity 
 are close to 4 in all norms
 while those of the scalar $q$ and its gradient
 show order reductions, 
 which are caused by the fact
 that the Neumann boundary condition
 in (\ref{eq:GePUPSAV}g) has to be obtained
 from spatial derivatives of the velocity
 and calculating these derivatives
 incurs order reductions in finite-volume discretizations.




Let $\mathbf{U}_h^n$ and $\mathbf{W}_h^n$ denote
 finite-volume solutions that approximate
 cell averages of $\mathbf{u}$ and $\mathbf{w}$ at time $t^n$,
 respectively.
Then the \emph{$L_2$-norm for a finite-volume solution} $\mathbf{V}_h^n$
 is defined as
\begin{equation}
  \label{eq:L2norm4FV}
  \begin{array}{l}
  \|\mathbf{V}^n_h\|_{\mathcal C} :=
 \sqrt{\sum_{\mathcal{C}_{\mathbf{i}}}
 \|\mathcal{C}_{\mathbf{i}}\|\cdot|\mathbf{V}^n_{h, \mathbf{i}}|^2}
  \end{array}
\end{equation}
 where $\mathcal{C}_{\mathbf{i}}$ ranges over all control volumes.
 The modified energy is then
 \begin{equation}
   \label{eq:modifiedEnergyDiscrete}
   \begin{array}{l}
     \mathcal{E}_h^n :=
     \frac{1}{2}\left(\|\mathbf{U}^n_h\|_{\mathcal C}^2 +
     |r_h^n|^2\right), 
   \end{array}
 \end{equation}
where $r_h^n$ is the computed value of the SAV $r$ at time $t^n$.

As shown in Figure \ref{fig:singleVortexBox-energyAndDiv}(a,b),
 over the entire simulation
 $\mathcal{E}_h^n$ decreases monotonically
 and $|r_h^n-1|$ remains below $2.0\times 10^{-11}$, 
 indicating that the (unmodified) kinetic energy
 $\frac{1}{2}\|\mathbf{U}^n_h\|_{\mathcal{C}}^2$
 also decreases monotonically.

In Figure \ref{fig:singleVortexBox-energyAndDiv}(c,d), 
 the $L_2$-norm of velocity divergence
 on the coarsest grid $h=\frac{1}{512}$
 first decreases dramatically
 during the first several seconds, 
 then gradually increases,
 and oscillates within a certain range.
Furthermore, the range of oscillation decreases quickly
 as the grid is refined. 
To understand this evolution pattern, 
 we consider the fully discrete counterpart of
 \eqref{eq:decayOfVelDiv-GePUP-SAV-RK}, i.e., 
 \begin{equation}
   \label{eq:discDivEvo}
   \left\|\mathbf{D} \mathbf{W}_h^{n+1}\right\|_{\mathcal C}^2
   - \left\|\mathbf{D} \mathbf{W}_h^{n}\right\|_{\mathcal C}^2
   \le -2k\nu \sum_{i=1}^s b_i\left\|\mathbf{G}\mathbf{D}
     \mathbf{W}_h^{(i)}\right\|_{\mathcal C}^2
   + O(h^p),
 \end{equation}
 where $\mathbf{D}$ is the discrete divergence, 
 $\mathbf{G}$ the discrete gradient, 
 and $O(h^p)$ errors of spatial discretization. 
 Numerical results 
 such as those in Table \ref{tab:singleVortexBox}
 show $p>1$. 

The key difference between (\ref{eq:discDivEvo}) and
 \eqref{eq:decayOfVelDiv-GePUP-SAV-RK}
 is the extra term $O(h^p)$ in (\ref{eq:discDivEvo}),
 which explains why the discrete divergence
 does not decay monotonically.
Since the initial condition of $\mathbf{w}(t_0)$ for this test
 is not divergence free, 
 the RHS of (\ref{eq:discDivEvo}) is dominated,
 in the early simulation stage,
 by the first term that contains
 $\|\mathbf{G}\mathbf{D}\mathbf{W}^{(i)}_h\|_{\mathcal{C}}$. 
Hence (\ref{eq:discDivEvo}) dictates the decay of the discrete divergence.
However,
 as $\|\mathbf{D}\mathbf{W}^{(i)}_h\|_{\mathcal{C}}$ decreases
 continuously, 
 the RHS of (\ref{eq:discDivEvo}) eventually
 becomes dominated by $O(h^p)$. 
Then the inequality \eqref{eq:discDivEvo} loses control over
 the discrete divergence
 since $O(h^p)$ is not negative-definite. 
When the discrete divergence increases
 to the degree such that the magnitude of $O(h^p)$ is less than
 that of the other RHS term, 
 the above pattern is repeated,
 leading to the oscillation of the discrete divergence.
The bottom line is, however,
 that the discrete divergence is indirectly controlled by
 the term $O(h^p)$
 and thus the oscillation becomes less prominent
 as the grid is refined. 
 
\subsection{Viscous-box tests} 
\label{sec:viscous-box-2D}

Following \cite{bell89:_secon_order_projec_method_incom},
 we set the initial velocity on $\Omega=[0,1]^2$
 to
 \begin{equation}
   \label{eq:viscousBoxVel}
   \mathbf{u}_0(x,y) 
   = 
   \left(
     \begin{array}{c}
       \, \sin^2(\pi x) \, \sin(2 \pi y) \,
       \\
       \,  -\sin(2 \pi x) \, \sin^2(\pi y) \,
     \end{array}
   \right)
 \end{equation}
 and advance cell-averaged initial values 
 from $t_0=0$ to $t_e=0.5$ 
 on four successively refined uniform grids.
The Courant number is defined as
 $\mathrm{Cr}:=\frac{k}{h}\|\mathbf{u}_0\|_{\infty}$, 
 where $\|\mathbf{u}_0\|_{\infty}$
 is the max-norm of the initial velocity $\mathbf{u}_0$.
 
 \begin{table}
   \caption{
     Errors and convergence rates of the GES scheme
     in Table \ref{tab:GePUP-SAV-SDIRK-policies}
     for solving the 2D viscous-box test
     with Re $=10^4$,
     $t_0=0.0$, $t_e=0.5$, $\mathrm{Cr}=0.5$ and $\lambda=1$.
   }
   \centering
   \input{tab/viscousBox_Re1e4_relaxCoef1}
   \label{tab:viscousBoxRe10000}
 \end{table}

 \begin{table}
   \caption{
     Errors and convergence rates of the GES scheme
     in Table \ref{tab:GePUP-SAV-SDIRK-policies}
     for solving the 2D viscous-box test
     with $\text{Re}=10^2$,
     $t_0=0.0$, $t_e=0.5$, $\mathrm{Cr}=0.1$ and $\lambda=1$. 
   }
   \centering
   \input{tab/viscousBox_Re100_relaxCoef1}
   \label{tab:viscousBoxRe100}
 \end{table}

 \begin{figure}
    \centering
    \subfigure[$\mathcal{E}_h^n$ for Re $=10^2$ and $\mathrm{Cr}=0.1$] 
    {
      \includegraphics[width=0.43\textwidth]{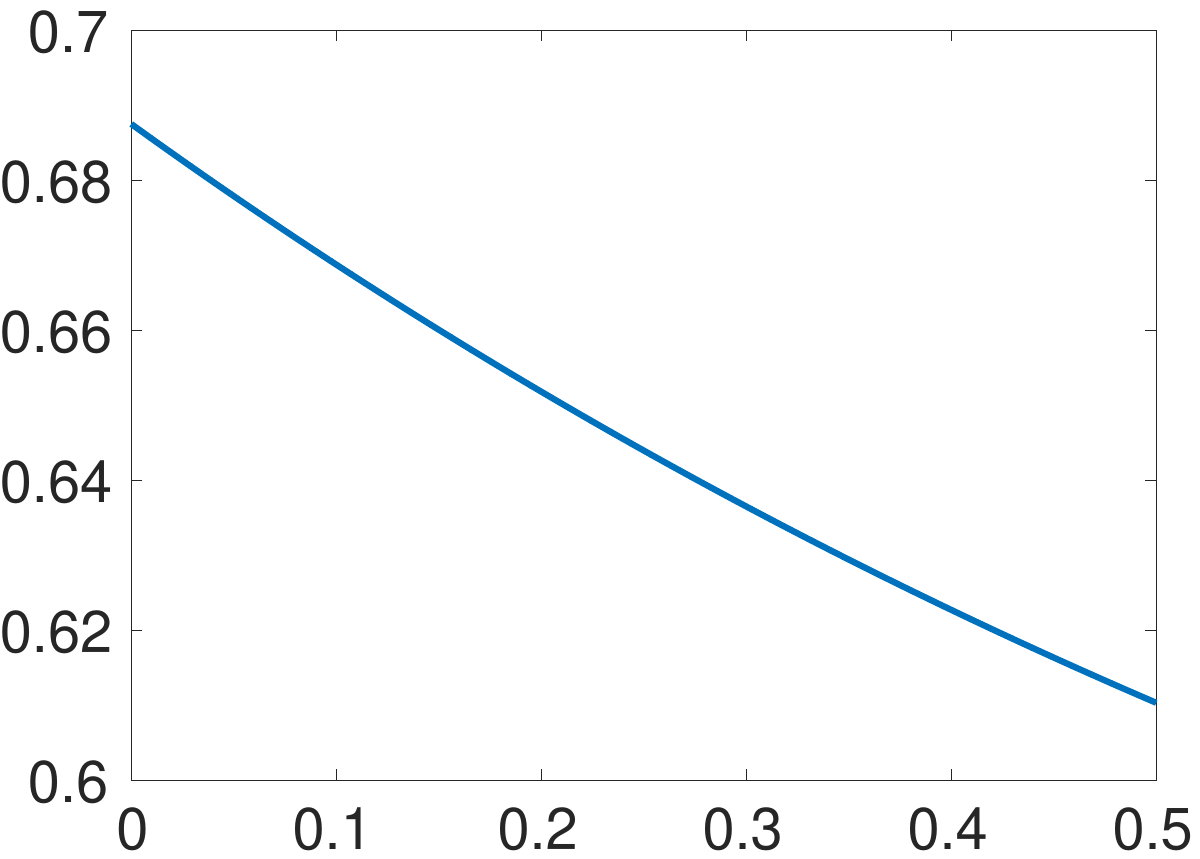}
    }
    \hfill
    \subfigure[$\mathcal{E}_h^n$ for Re $=10^4$ and $\mathrm{Cr}=0.5$]
    {
      \includegraphics[width=0.43\textwidth]{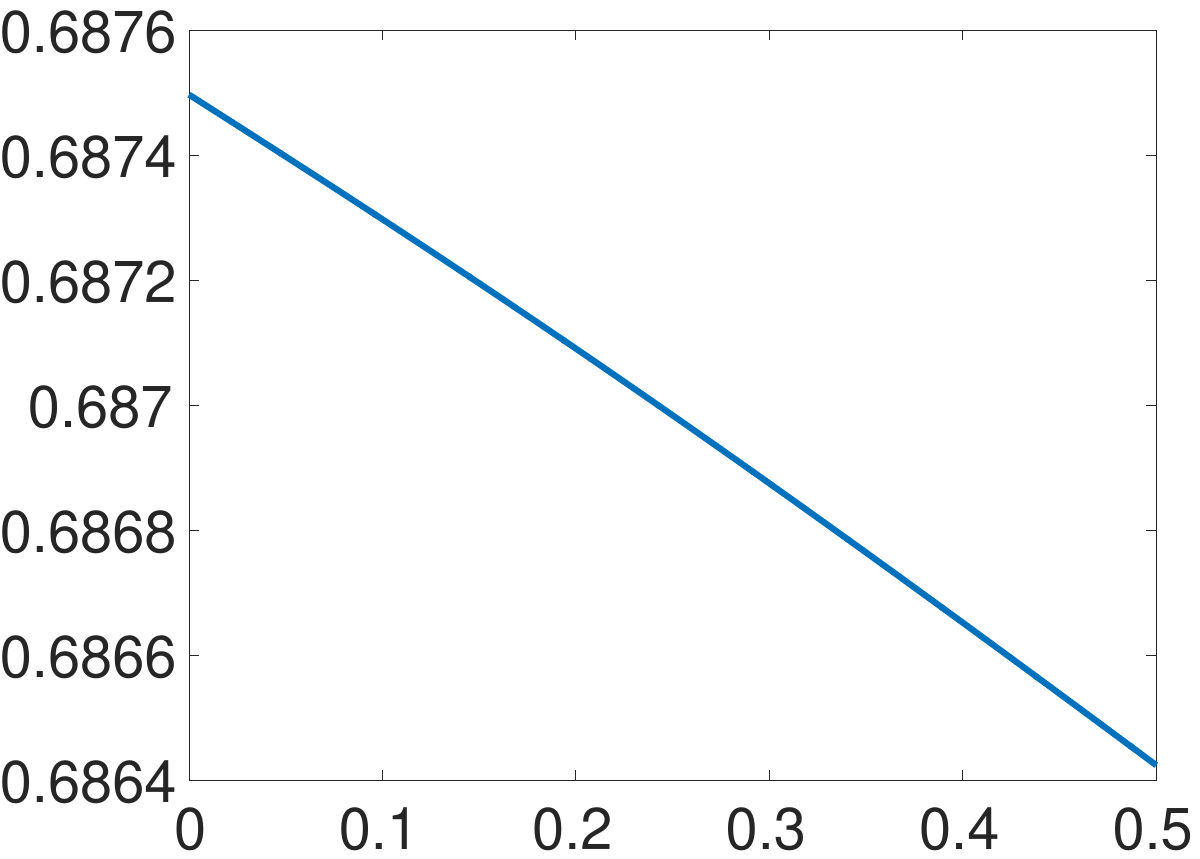}
    }

    \subfigure[$|r_h^n-1|$ for Re $=10^2$ and $\mathrm{Cr}=0.1$]
    {
      \includegraphics[width=0.44\textwidth]{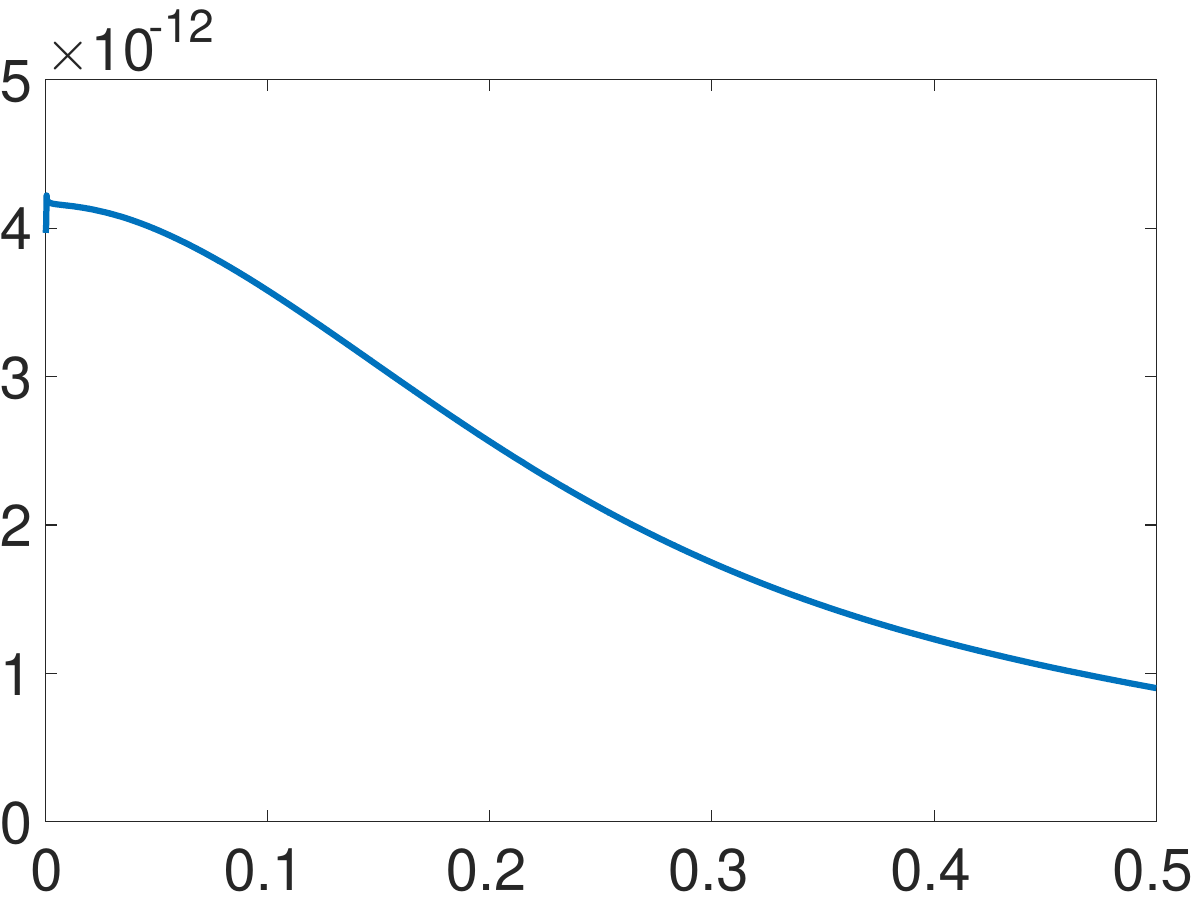}
    }
    \hfill
    \subfigure[$|r_h^n-1|$ for $\mathrm{Re}=10^4$ and $\mathrm{Cr}=0.5$]
    {
      \includegraphics[width=0.44\textwidth]{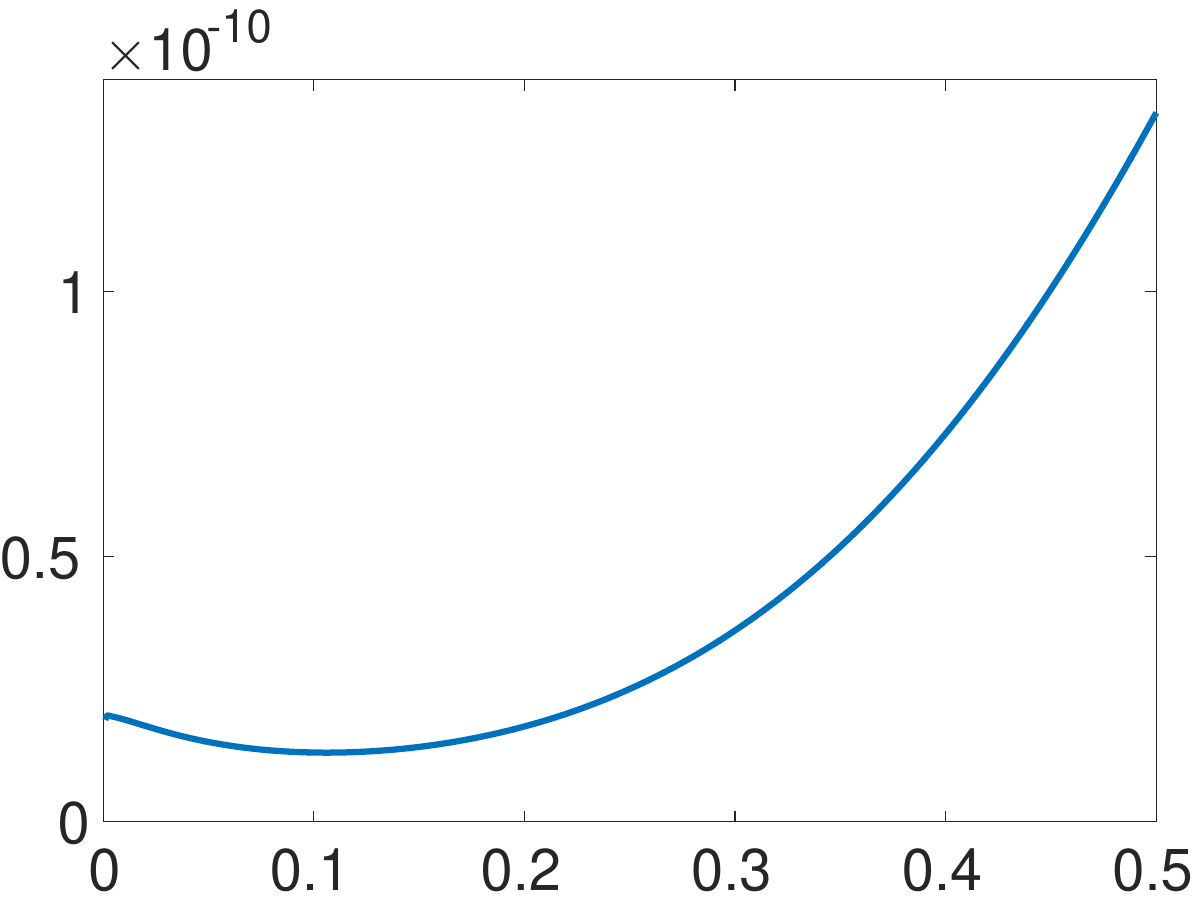}
    }

    \subfigure[$\left\|\mathbf{D} \mathbf{W}_h^{n}\right\|_{\mathcal C}$
    (solid markers) and 
    $\left\|\mathbf{D} \mathbf{U}_h^{n}\right\|_{\mathcal C}$
    (hollow markers)
    for Re $=10^2$ and $\mathrm{Cr}=0.1$]
    {
      \includegraphics[width=0.475\textwidth]{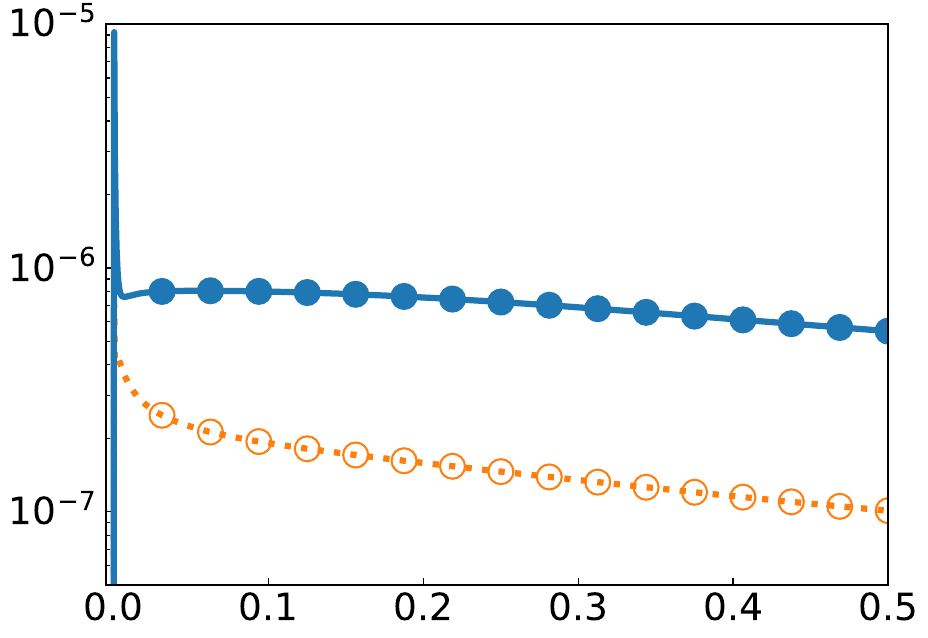}
    }
    \hfill
    \subfigure[$\left\|\mathbf{D} \mathbf{W}_h^{n}\right\|_{\mathcal C}$
    (solid markers) and 
    $\left\|\mathbf{D} \mathbf{U}_h^{n}\right\|_{\mathcal C}$
    (hollow markers)
    for Re $=10^4$ and $\mathrm{Cr}=0.5$]
    {
      \includegraphics[width=0.475\textwidth]{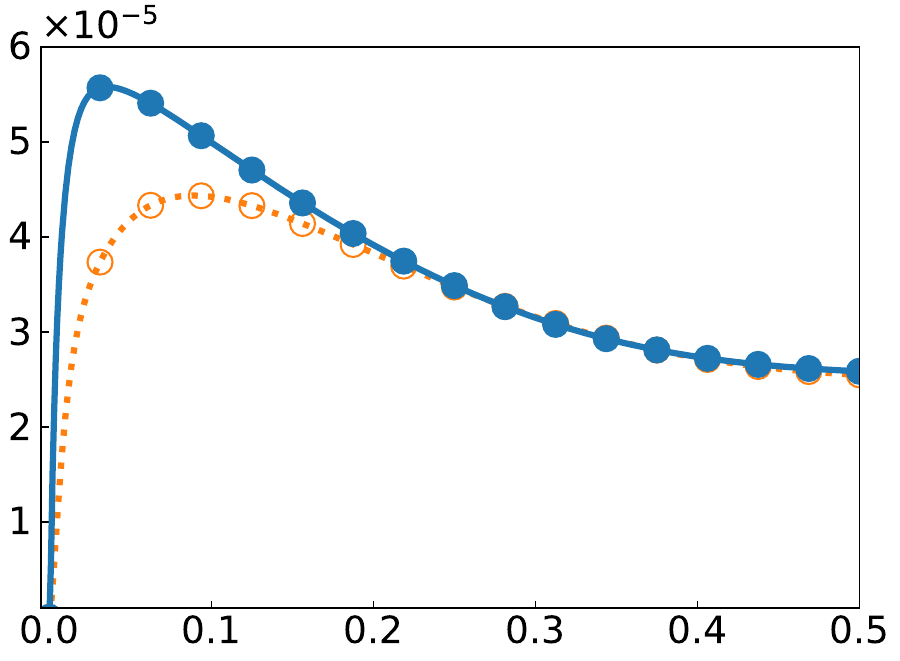}
    }
    \caption{Results of the GES scheme
      in Table \ref{tab:GePUP-SAV-SDIRK-policies}
      for solving the 2D viscous-box test
      with $h=\frac{1}{512}$ and $\lambda=1$;
      see (\ref{eq:L2norm4FV})
      and (\ref{eq:modifiedEnergyDiscrete})
      for precise definitions of
      $\left\|\cdot \right\|_{\mathcal C}$ and $\mathcal{E}_h^n$. 
      The abscissa in all subplots is time.
    }
    \label{fig:viscousBox}
  \end{figure}

Errors and convergence rates
 in the cases of $\mathrm{Re}=10^4$ and $\mathrm{Re}=10^2$
 are shown in Tables \ref{tab:viscousBoxRe10000}
 and \ref{tab:viscousBoxRe100}, respectively.
For $\mathrm{Re}=10^4$,
 convergence rates of the velocity in all norms
 are close to four.
In contrast, those for $\mathrm{Re}=10^2$
 are close to four in the $L_1$-norm and the $L_2$-norm, 
 but are around 2.5 in the $L_{\infty}$-norm.
Accordingly, 
 convergence rates of $\nabla q$ 
 in the $L_{\infty}$-norm for $\mathrm{Re}=10^4$
 are also substantially higher
 than those for $\mathrm{Re}=10^2$.

These results are not out of expectations.
We have proved in Theorem \ref{thm:wConverge2u}
 that the pressure $q$ converges to the pressure $p$
 and have shown in Section \ref{sec:stokes-pressure} that
 the pressure gradient $\nabla p$ in the INSE can be split into two parts
 $\nabla p=\nabla p_c+\nu \nabla p_s$,
 where $\nabla p_s$
 responses to the Laplace-Leray commutator.
When $\nu$ is sufficiently large,
 $\nu\nabla p_s$ dominates $\nabla p_c$
 and accounts for the bulk of $\nabla p$.
Cozzi and Pego \cite{cozzi11:_laplac_leray} showed that
 $\|\nabla p_s\|$ may not be bounded
 at a boundary point that is not $C^3$.
As a practical interpretation,
 the pressure could develop steep gradient
 at a $C^1$ discontinuity of the domain boundary
 for low-Reynolds-number flows.
Therefore, 
 we believe
 that the order reduction in the case of Re $=10^2$
 is caused by the dominance of $\nu\nabla p_s$
 and the sharp corners ($C^1$ discontinuities)
 of the square domain.
 
As shown in Figure \ref{fig:viscousBox}(a,b),
 the modified energy $\mathcal{E}_h^n$
 decreases monotonically
 over the entire simulation
 for both $\mathrm{Re}=10^4$ and $\mathrm{Re} = 10^2$.
It is also clear that $\mathcal{E}_h^n$
 decreases faster in the higher-viscosity case of Re $=10^2$, 
 which confirms Theorem \ref{thm:GePUPSAVRKEnergyDecay}.
Figure \ref{fig:viscousBox}(c,d) show that
 the deviation of SAV from 1 is at most $10^{-10}$,
 thus the (unmodified) kinetic energy
 also decreases monotonically. 

For both $\mathrm{Re}=10^4$ and $\mathrm{Re}=10^2$, 
 the evolution of the $L_2$-norm of velocity divergence
 shown in Figure \ref{fig:viscousBox}(e,f)
 has essentially the same pattern:
 the $L_2$-norm first increases to a local maximum
 and then decreases.
This pattern is different from that
 shown in Figure \ref{fig:singleVortexBox-energyAndDiv}(c,d),
 but can still be very well explained by (\ref{eq:discDivEvo}). 
Since the initial velocity (\ref{eq:viscousBoxVel}) is divergence-free,
 the magnitude of
 $\|\mathbf{G}\mathbf{D}\mathbf{W}^{(i)}_h\|_{\mathcal{C}}$ is small 
 in the early simulation stage, 
 during which the inequality \eqref{eq:discDivEvo}
 has no control over the discrete velocity divergence yet.
However, as the discrete divergence accumulates
 to the point when the RHS of (\ref{eq:discDivEvo}) gets dominated by
 its first term, 
 the inequality \eqref{eq:discDivEvo}
 takes effect and forces the discrete divergence to decrease.

\subsection{Viscous-box tests
   with initially non-solenoidal velocity} 
\label{sec:viscous-box-2D-nonDivFree}

Tests in this subsection
 are the same as those in Section \ref{sec:viscous-box-2D}
 except that the initial condition
 is set to the following non-solenoidal velocity, 
 \begin{equation}
   \label{eq:initialVelNonDivFree}
   \begin{array}{l}
   \mathbf{w}_0
   = \mathbf{u}_0 + \epsilon\nabla \phi,\quad 
   \text{ with }
   \phi=\frac{1}{2\pi^2} \sin^2(\pi x) \sin^2(\pi y), 
   \end{array}
 \end{equation}
 where $\mathbf{u}_0$ is the initially solenoidal velocity 
 in (\ref{eq:viscousBoxVel}), 
 $\epsilon \nabla \phi$ the perturbation of $\mathbf{u}_0$
 away from the divergence-free space,
 and $\epsilon>0$ the parameter to control
 the magnitude of the perturbation.

\begin{figure}
  \centering
  \subfigure[$\epsilon=1$, $t=t_0=0$]{
    \includegraphics[width=0.31\textwidth]{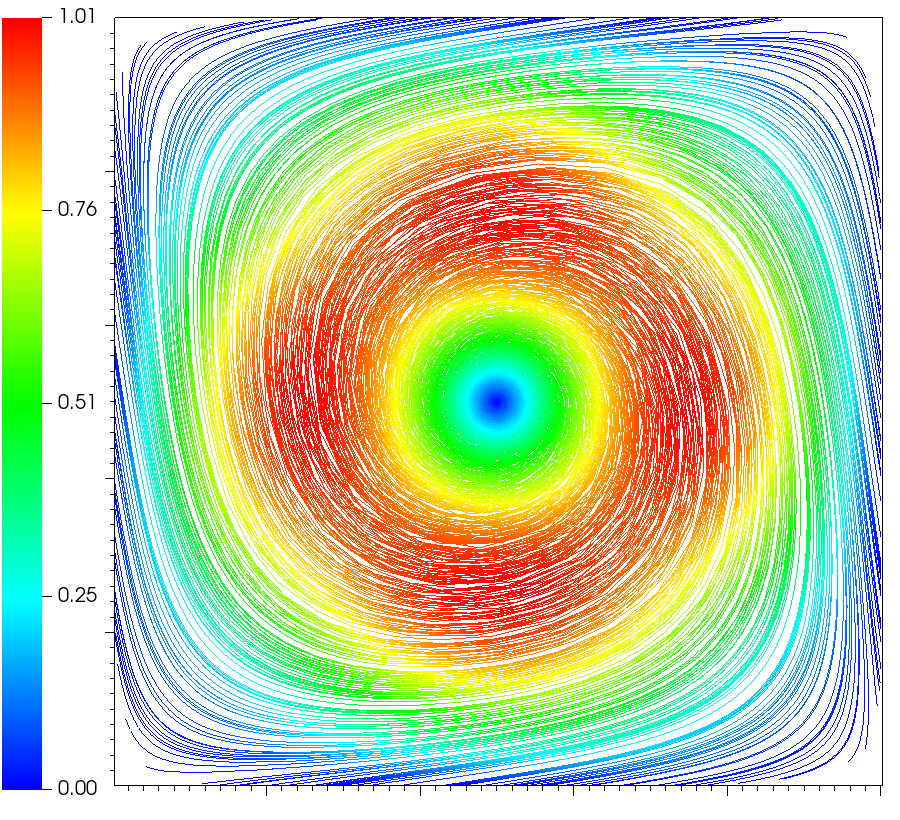}
  }
  \hfill
  \subfigure[$\epsilon=0$, $t=t_0=0$]{
    \includegraphics[width=0.31\textwidth]{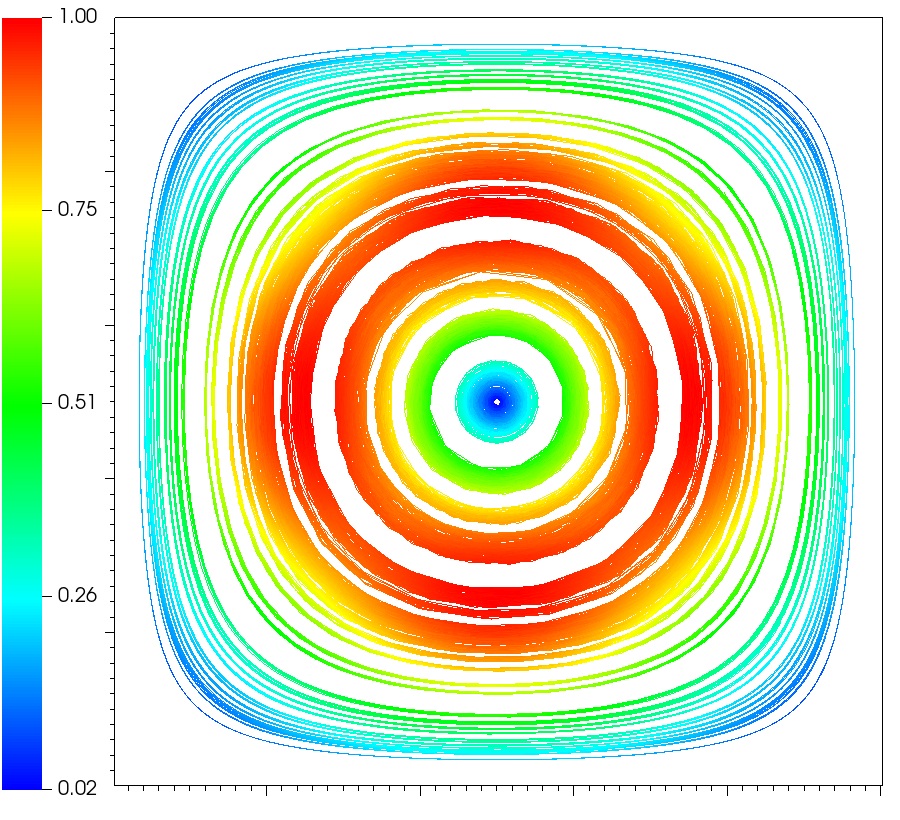}
  }
  \hfill
  \subfigure[$\epsilon\!=\!10^{-3}$, $t=t_0=0$]{
    \includegraphics[width=0.31\textwidth]{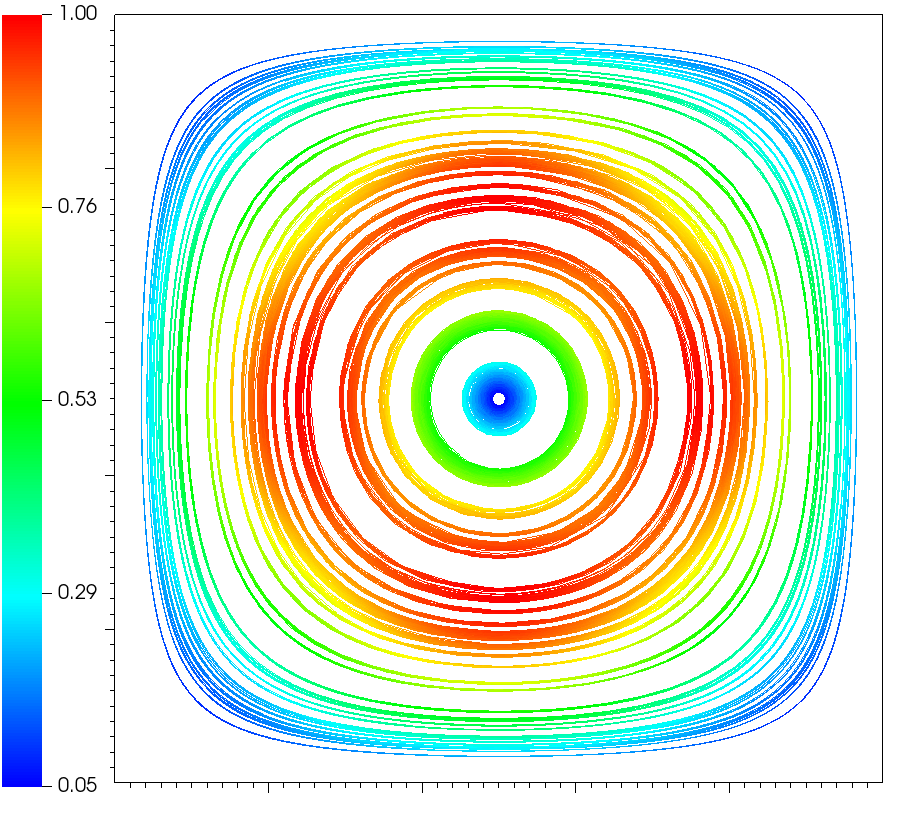}
  } 
  \\
 \subfigure[$\epsilon=1$, $t=0.5$, Re $=10^{4}$]{
    \includegraphics[width=0.31\textwidth]{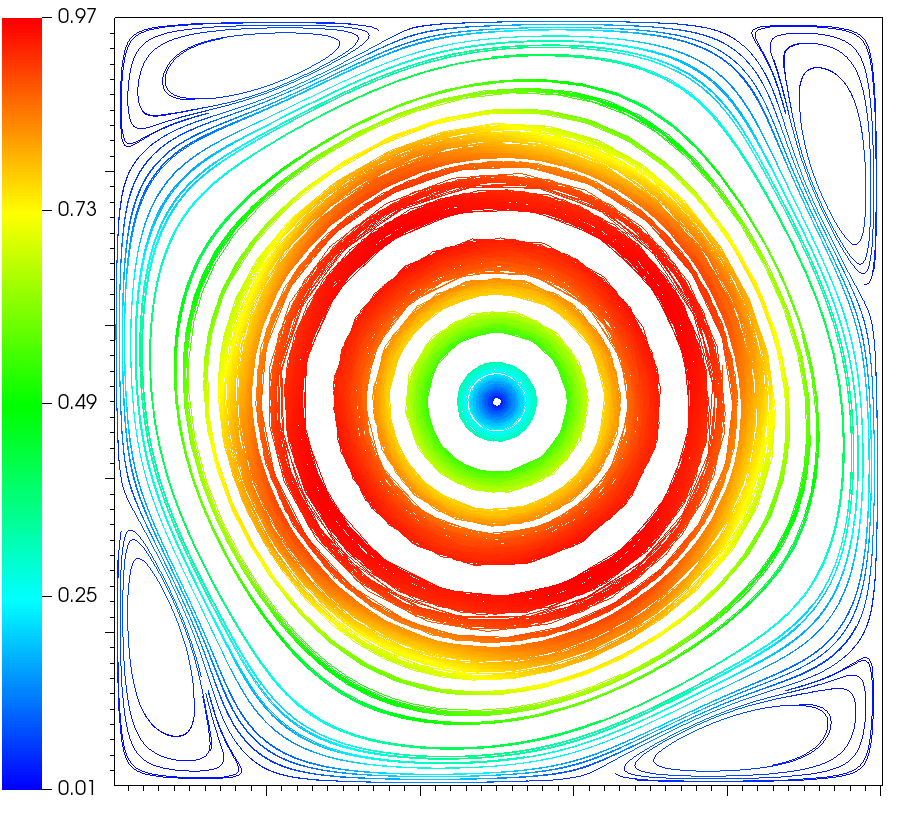}
  }
  \hfill
  \subfigure[$\epsilon=0$, $t=0.5$, Re $=10^{4}$]{
    \includegraphics[width=0.31\textwidth]{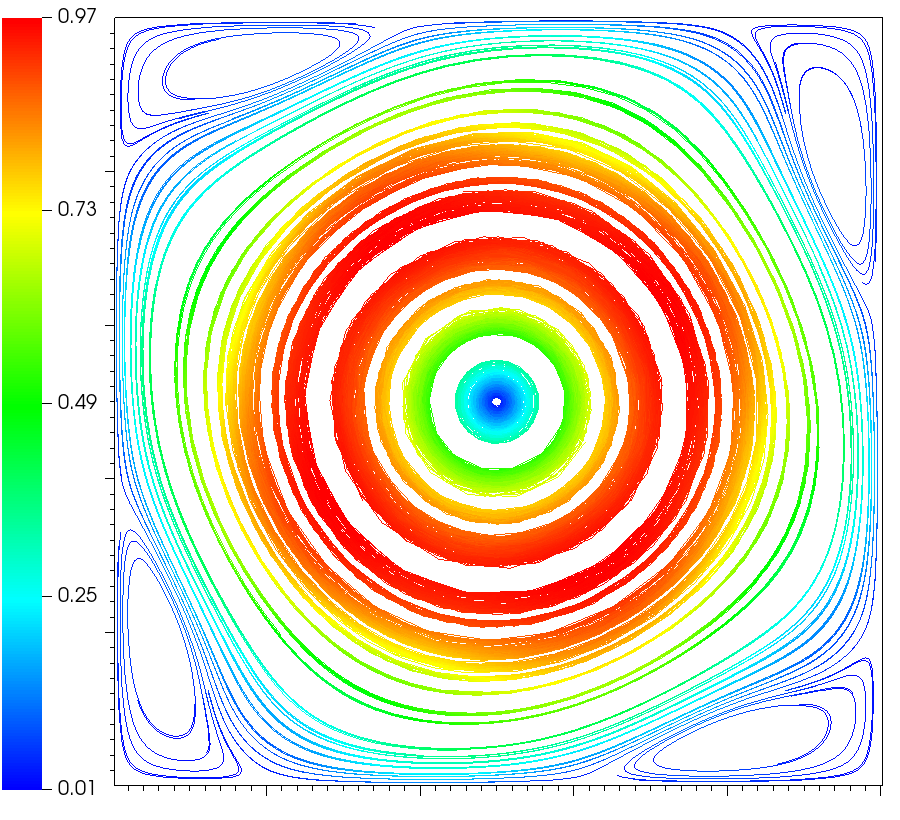}
  } 
  \hfill
  \subfigure[$\epsilon\!=\!10^{-3}$, $t\!=\!0.5$, Re $=10^{4}$]{
    \includegraphics[width=0.31\textwidth]{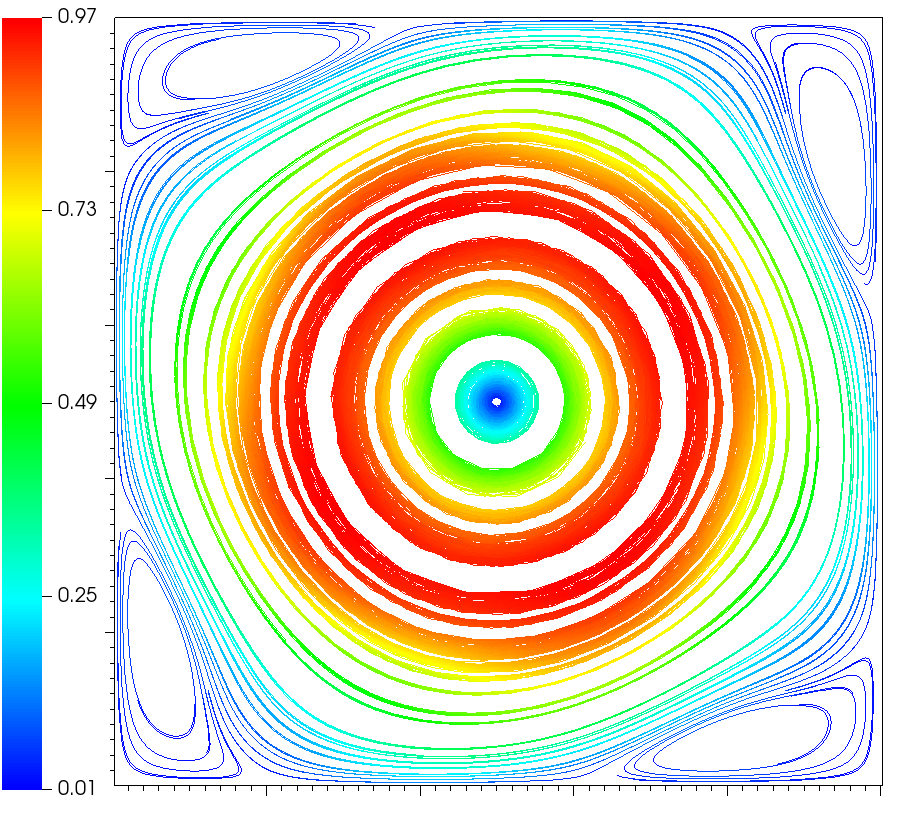}
  } 
  \\
  \subfigure[$\epsilon=1$, $t=0.5$, Re $=10^{2}$]{
    \includegraphics[width=0.31\textwidth]{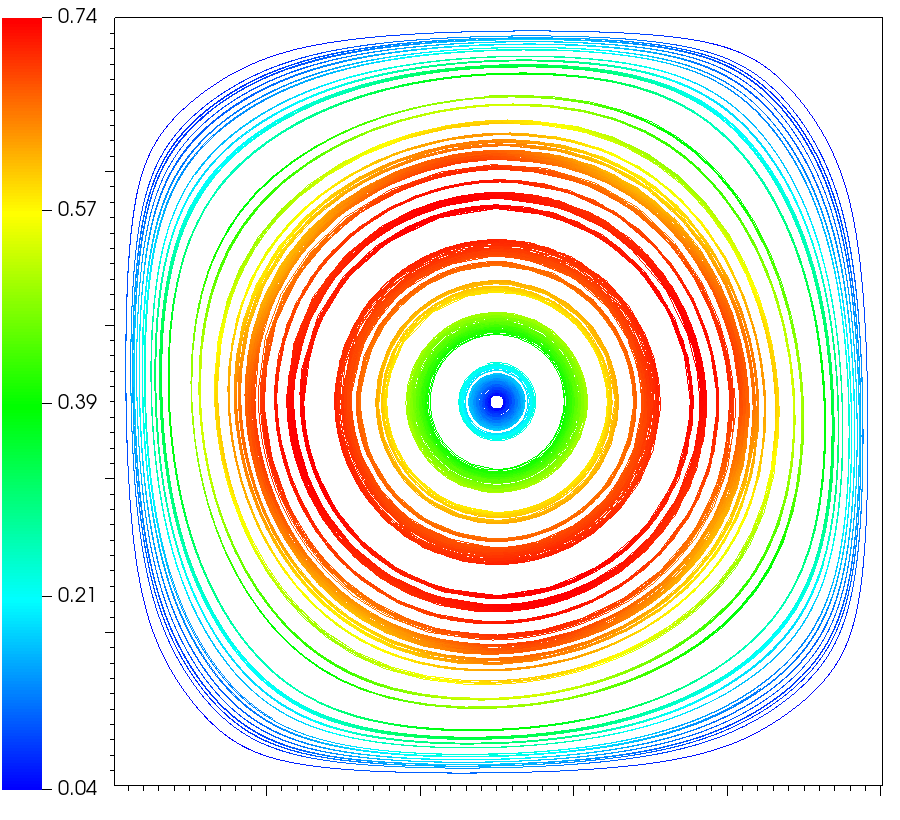}
  }
  \hfill
  \subfigure[$\epsilon=0$, $t=0.5$, Re $=10^{2}$]{
    \includegraphics[width=0.31\textwidth]{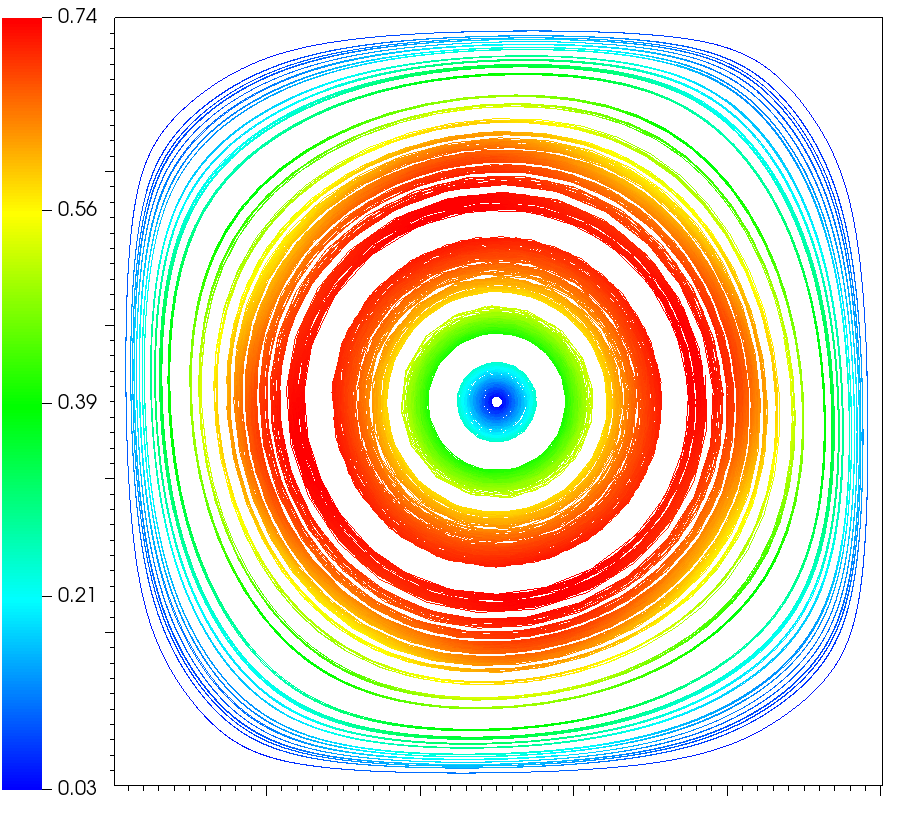}
  } 
  \hfill
  \subfigure[$\epsilon\!=\!10^{-3}$, $t=0.5$, Re $=10^{2}$]{
    \includegraphics[width=0.31\textwidth]{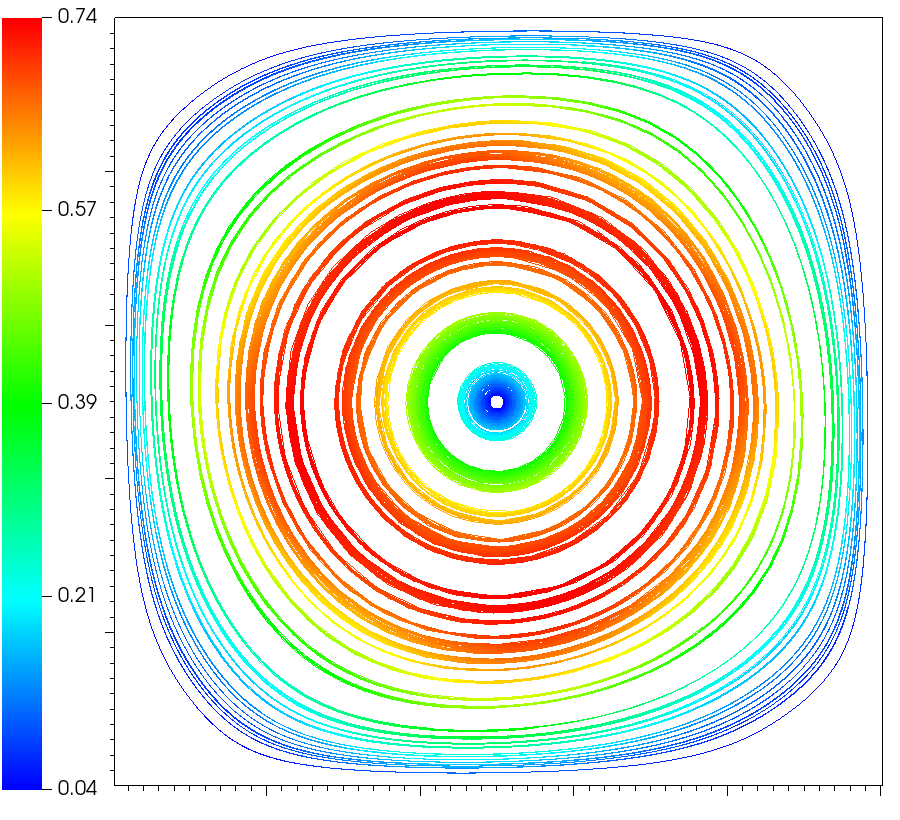}
  }
  \caption{Snapshots of streamlines of the initial velocity (the first row)
    and the final solutions 
    of the GES scheme in Table \ref{tab:GePUP-SAV-SDIRK-policies}
    for solving perturbed viscous box tests 
    in Section \ref{sec:viscous-box-2D-nonDivFree}
    with 
    Re$ = 10^{4}$ (the second row)
    and Re$=10^{2}$ (the third row)
    on a uniform grid of $h = \frac{1}{256}$.
    The unperturbed solutions $\mathbf{U}^{n}_{\epsilon=0}$
    and the perturbed solutions $\mathbf{U}^{n}_{\epsilon>0}$
    are displayed in the middle column
    and other columns, respectively.
    The span of each plot is the domain $[0,1]^2$
    and the color indicates
    the velocity magnitude. 
  }
  \label{fig:velContourfigs_Re1e4}
\end{figure}

Hereafter we denote by
 $\mathbf{U}^n_{\epsilon=0}$
 the reference solution
 obtained from the divergence-free initial velocity $\mathbf{u}_0$, 
 i.e., $\mathbf{U}^n_h$ in Section \ref{sec:viscous-box-2D}, 
 and denote by $\mathbf{U}^{n}_{\epsilon>0}$
 the perturbed solution
 computed from a perturbed initial velocity $\mathbf{w}_0$.

Streamlines of the initial velocity
 and the final solutions 
 are plotted in Figure \ref{fig:velContourfigs_Re1e4},
 where no qualitative differences can be observed
 between corresponding plots in the third and the second columns, 
 due to the small perturbation magnitude
 $\epsilon=10^{-3}$.
In contrast,
 for $\epsilon=1$,
 the perturbation is sufficient
 to generate different patterns for the initial velocity, 
 cf. Figure \ref{fig:velContourfigs_Re1e4}(a,b). 
For example,
 the streamlines in Figure \ref{fig:velContourfigs_Re1e4}(a)
 are not closed due to the large magnitude of velocity divergence.
However,
 the perturbation effects have already died out at $t=0.5$
 and the final perturbed solution
 $\mathbf{U}^{n}_{\epsilon=1}$
 is visually indistinguishable
 from the unperturbed solution
 $\mathbf{U}^{n}_{\epsilon=0}$
 for both Re $= 10^{4}$ and Re $=10^{2}$; 
 see Figure \ref{fig:velContourfigs_Re1e4}(d,e) and
 Figure \ref{fig:velContourfigs_Re1e4}(g,h).

\begin{figure}
  \centering
  \subfigure[Re $=10^2$, Cr $=0.1$]{
    \includegraphics[width=0.311\textwidth]{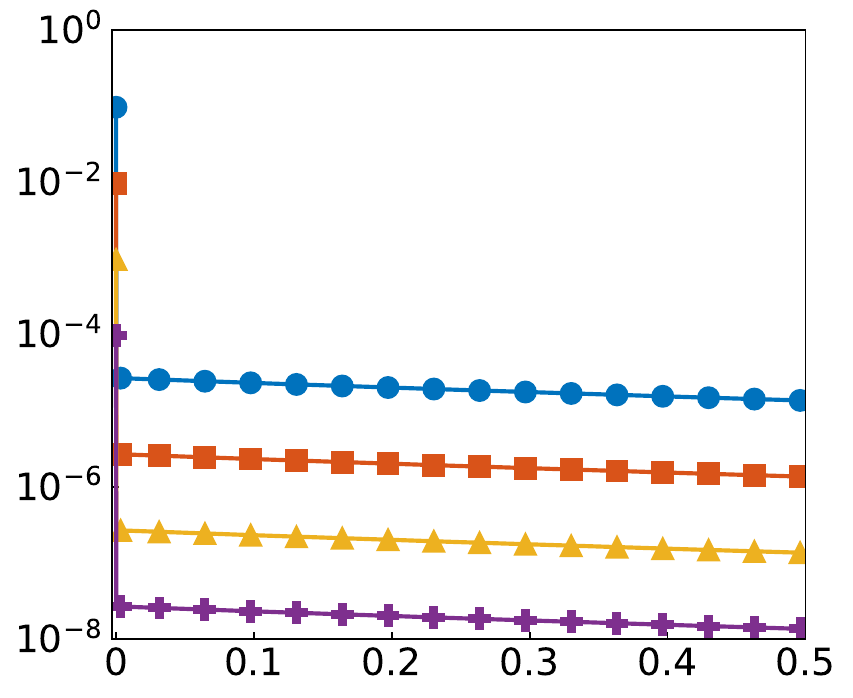}
  }
  \hfill
  \subfigure[Re $=10^2$, Cr $=0.1$]{
    \includegraphics[width=0.311\textwidth]{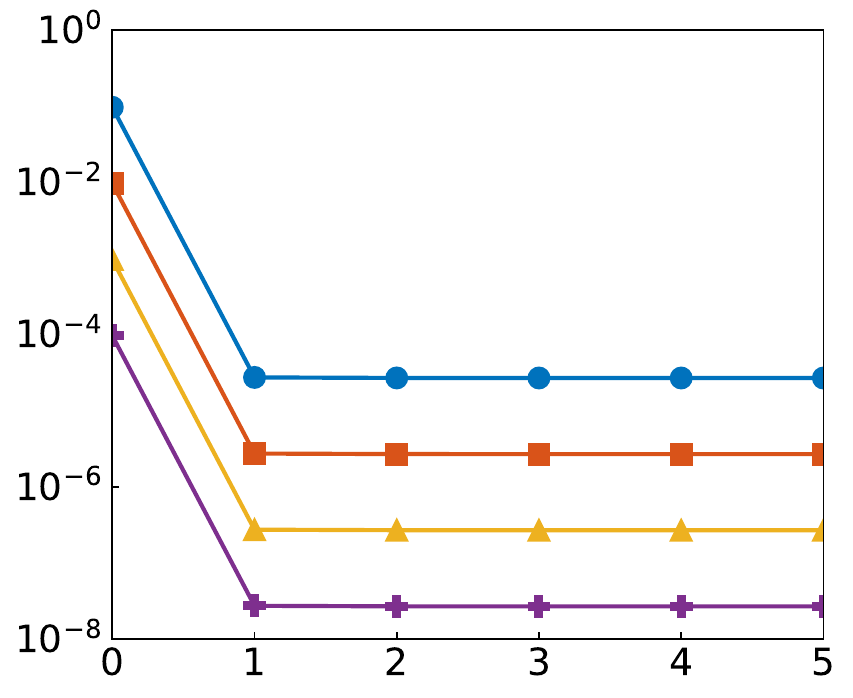}
  }
  \hfill
  \subfigure[Re $=10^2$, Cr  $=0.01$]{
    \includegraphics[width=0.311\textwidth]{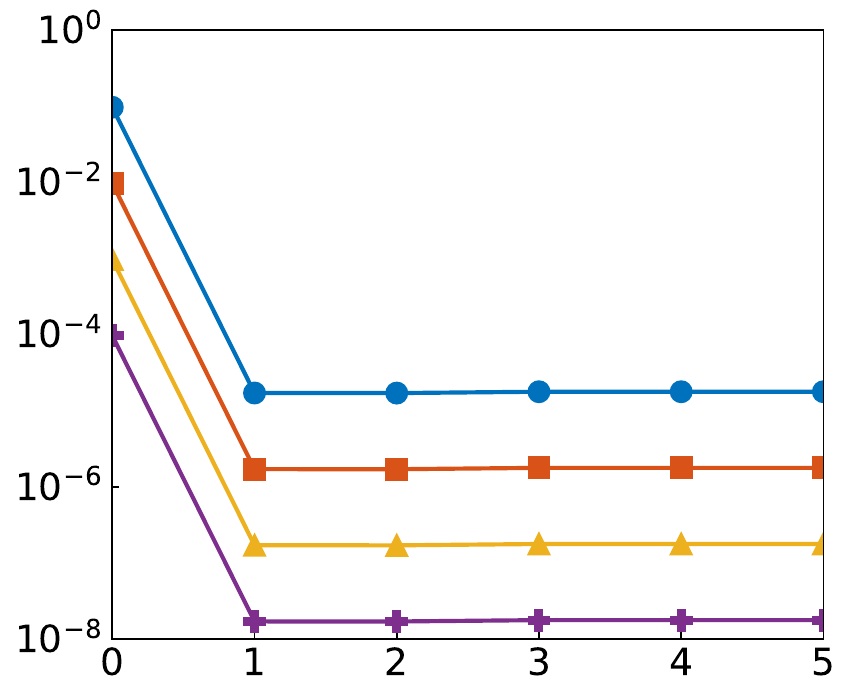}
  }
  \\
  \subfigure[Re $=10^{4}$, Cr $=0.5$]{
    \includegraphics[width=0.311\textwidth]{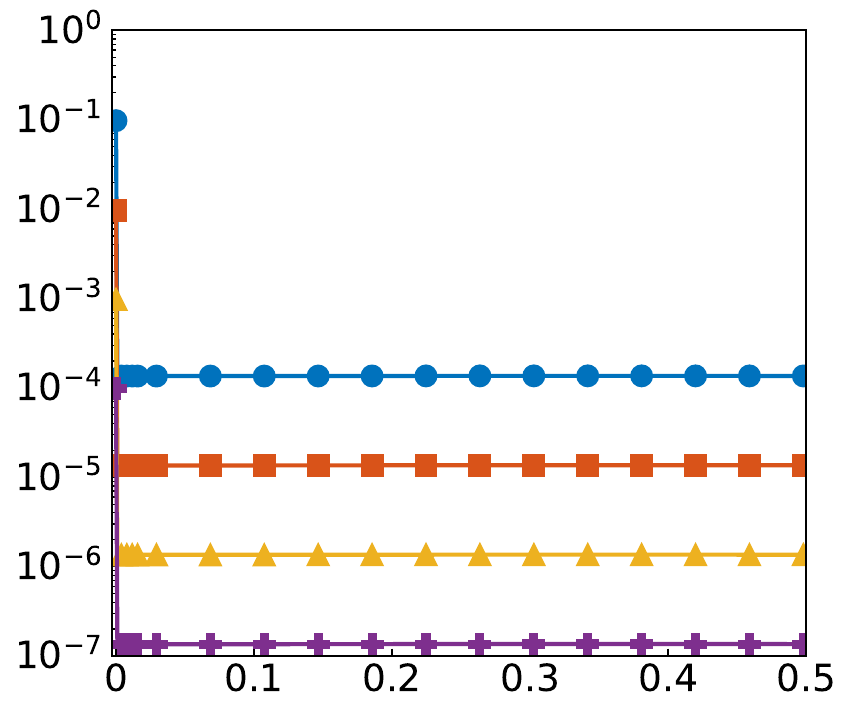}
  }
  \hfill
  \subfigure[Re $=10^{4}$, Cr $=0.5$]{
    \includegraphics[width=0.311\textwidth]{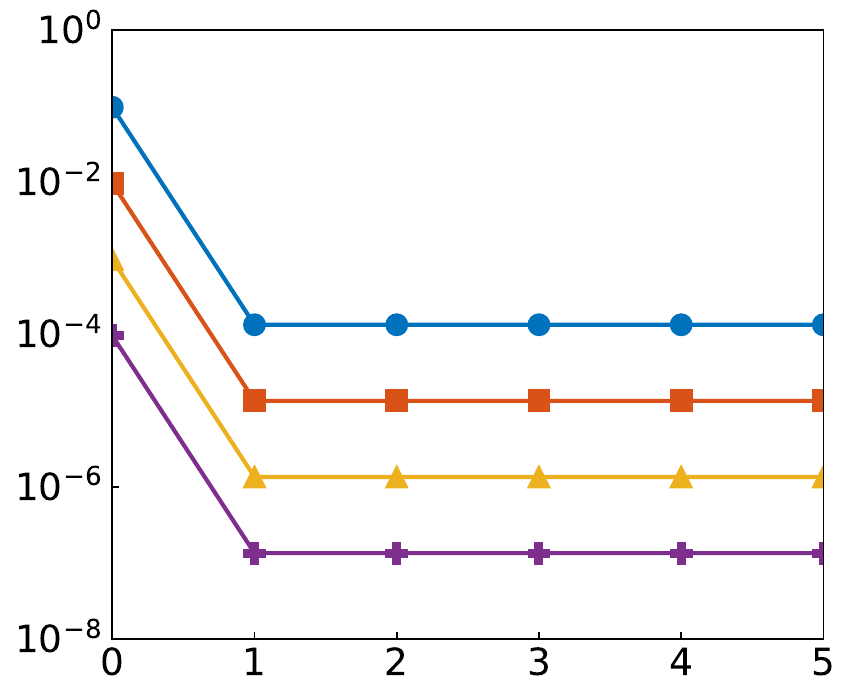}
  }
  \hfill
  \subfigure[Re $=10^{4}$, Cr $=0.05$]{
    \includegraphics[width=0.311\textwidth]{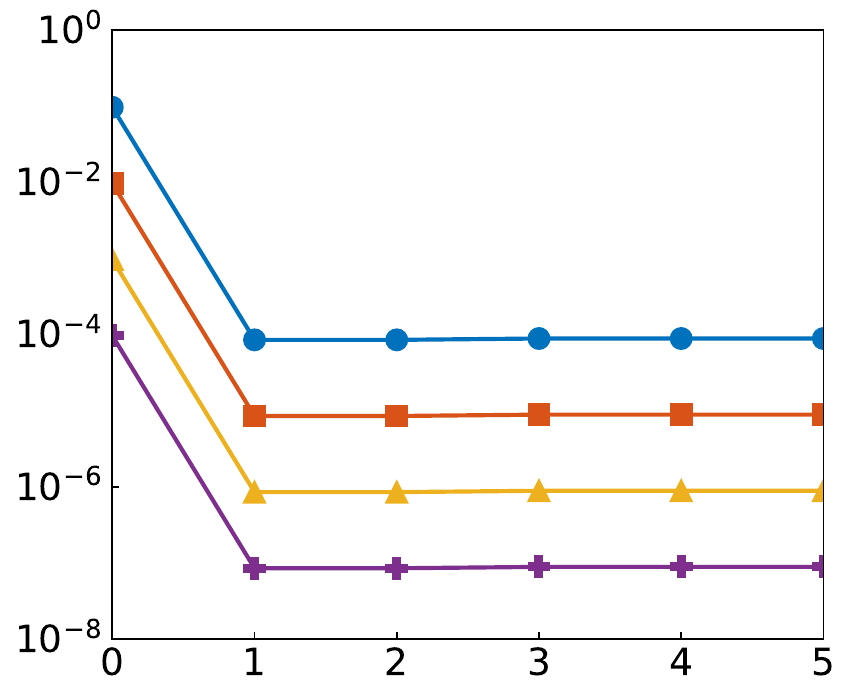}
  }
  \caption{Evolutions of $\|\mathbf{U}^n_{\epsilon>0} -
      \mathbf{U}^{n}_{\epsilon = 0}\|_{\mathcal{C}}$,
    i.e., the $L_2$-norm (\ref{eq:L2norm4FV}) of the difference
    between the computed velocities
    of the perturbed viscous-box tests
    in Section \ref{sec:viscous-box-2D-nonDivFree}
    on a uniform spatial grid $h=\frac{1}{512}$. 
    The markers
    ``$\bullet$,''  ``{\tiny $\blacksquare$},''
    ``$\blacktriangle$'' and ``$+$''
    represent values of
    $\|\mathbf{U}^n_{\epsilon>0} - \mathbf{U}^{n}_{\epsilon = 0}\|_{\mathcal{C}}$
    for $\epsilon=1$, $0.1$, $10^{-2}$, and $10^{-3}$,
    respectively.
    The abscissas in subplots (b,c,e,f) and (a,d)
    represent the index of time steps and the simultation time,
    respectively. 
   }
  \label{fig:diffVelNormsEvolution}
\end{figure}

How fast does the perturbed solution $\mathbf{U}^{n}_{\epsilon>0}$
 converge to the solution $\mathbf{U}^{n}_{\epsilon=0}$?
The answer lies in Figure \ref{fig:diffVelNormsEvolution},
 where the deviation of 
 $\mathbf{U}^{n}_{\epsilon>0}$ from $\mathbf{U}^{n}_{\epsilon=0}$
 decays by a factor of at least $10^{3}$
 within the first time step for all test cases.
Furthermore, although the norm
$\|\mathbf{U}^{n}_{\epsilon>0} -\mathbf{U}^{n}_{\epsilon=0}\|_{\mathcal{C}}$
 appears to be proportional to the perturbation magnitude, 
 the reduction rate during the first time step
 is largely independent of the time step size; 
 see Figure \ref{fig:diffVelNormsEvolution}(b,c)
 and Figure \ref{fig:diffVelNormsEvolution}(e,f). 
 
\begin{figure}
  \centering
  \subfigure[$\mathcal{E}^n_h$ for Re $=10^2$ and Cr = 0.1]{
    \includegraphics[width=0.475\textwidth]{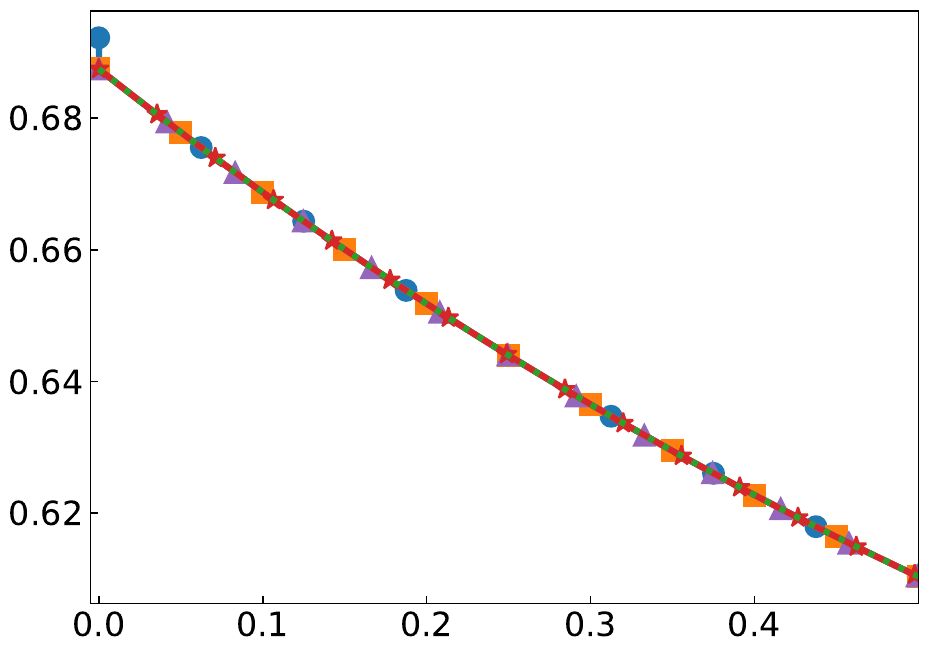}
  }
  \hfill
  \subfigure[$\mathcal{E}^n_h$ for Re $=10^{4}$ and Cr = 0.5]{
    \includegraphics[width=0.475\textwidth]{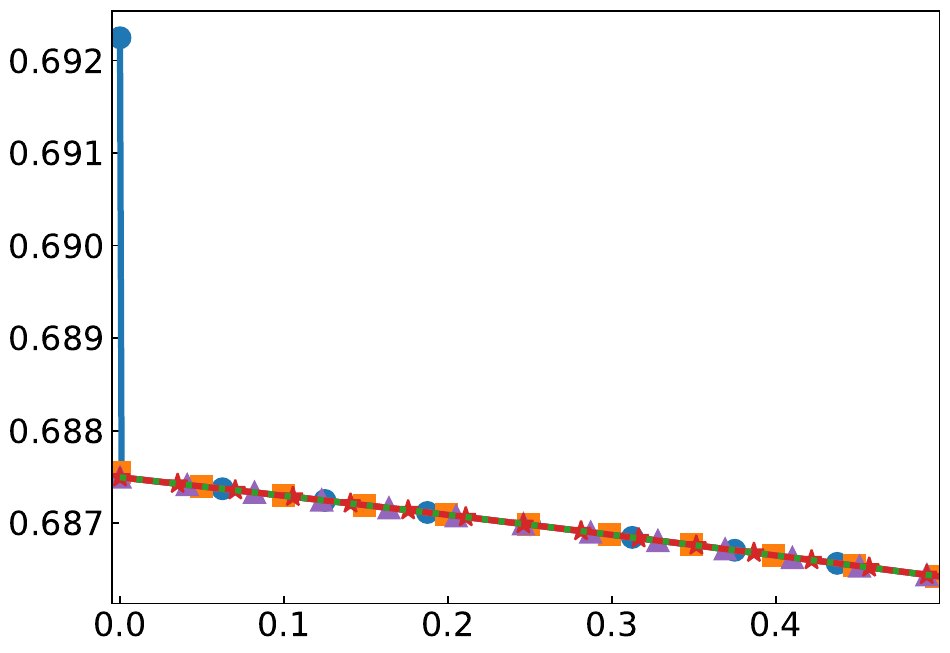}
  }
  \\
  \subfigure[$|r^n_h-1|$ for Re  $=10^2$ and Cr = 0.1]{
    \includegraphics[width=0.475\textwidth]{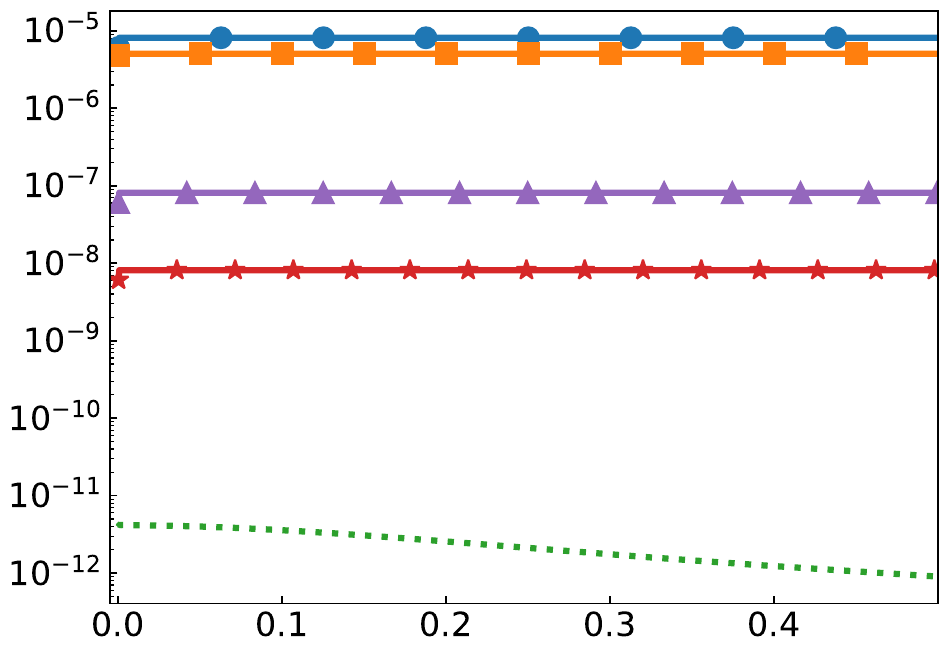}
  }
  \hfill
  \subfigure[$|r^n_h\!-\!1|$ for Re $=10^{4}$ and Cr = 0.5]{
    \includegraphics[width=0.475\textwidth]{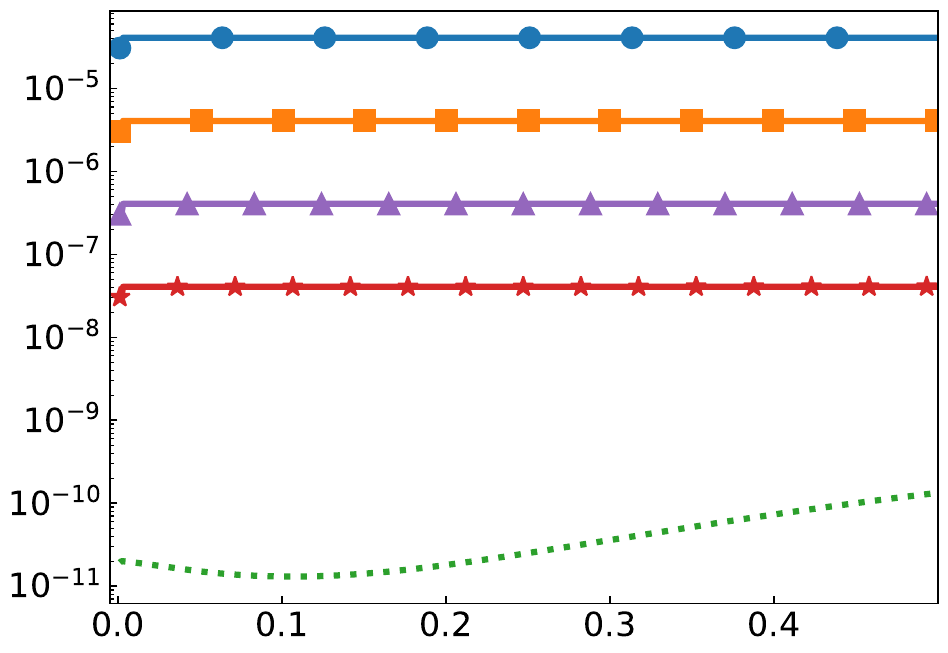}
  }
  \\
  \subfigure[$\|\mathbf{DU}^n_{\epsilon}\|_{\mathcal{C}}$ for Re $=10^2$ and Cr = 0.1]{
    \includegraphics[width=0.475\textwidth]{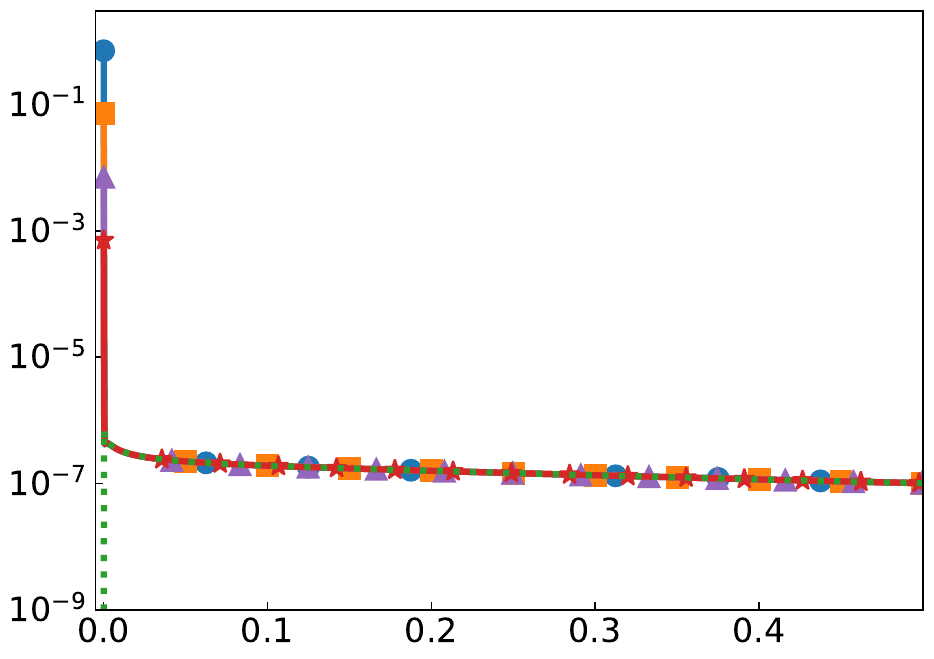}
  }
  \hfill
  \subfigure[$\|\mathbf{DU}^n_{\epsilon}\|_{\mathcal{C}}$ for Re $=10^{4}$ and Cr = 0.5]{
    \includegraphics[width=0.475\textwidth]{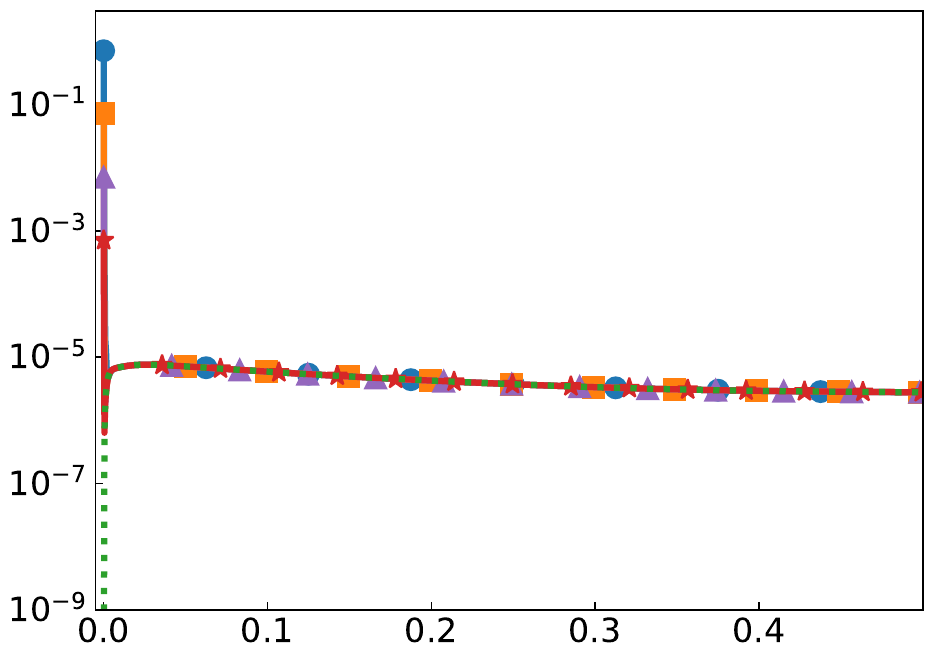}
  }
  \caption{Evolutions of $\mathcal{E}^n_h$, $|r^n_h-1|$, 
    and $ \|\mathbf{DU}^n_{\epsilon}\|_{\mathcal{C}}$
    produced by the GES scheme in Table \ref{tab:GePUP-SAV-SDIRK-policies}
    for solving the perturbed viscous box tests
    on a uniform grid of $h=\frac{1}{512}$;
    see (\ref{eq:L2norm4FV})
    and (\ref{eq:modifiedEnergyDiscrete})
    for precise definitions of
    $\left\|\cdot \right\|_{\mathcal C}$
    and $\mathcal{E}_h^n$. 
    The dashed curves represent results with $\epsilon=0$
    (i.e., those in Figure \ref{fig:viscousBox})
    while those marked by  ``$\bullet$,''
    ``{\tiny $\blacksquare$},''
    ``$\blacktriangle$,'' and ``$\bigstar$''
    correspond to $\epsilon=1$, $10^{-1}$, $10^{-2}$ and $ 10^{-3}$,
    respectively.
    The abscissa of all subplots is time.
  }
  \label{fig:GePUP-ES_energyAndDivCurves}
\end{figure}

For two curves very close to each other, 
 their tangent vectors might differ largely,
 so do other geometric quantities
 such as the curvature.
Likewise, $\|\mathbf{U}^{n}_{\epsilon>0}-\mathbf{U}^{n}_{\epsilon=0}\|_{\mathcal{C}}$
 being small does not imply
 $\|\Div\mathbf{U}^{n}_{\epsilon>0}-\Div\mathbf{U}^{n}_{\epsilon=0}\|_{\mathcal{C}}$
 being small.
This issue is addressed by Figure \ref{fig:GePUP-ES_energyAndDivCurves}, 
 which shows the perturbation effects 
 on evolutions of the modified energy, the SAV, 
 and the velocity divergence.
 
Compared to Figure \ref{fig:viscousBox}(a,b),
 Figure \ref{fig:GePUP-ES_energyAndDivCurves}(a,b)
 feature a decrease of the modified energy
 during the very first time step for $\epsilon=1$,
 indicating that most of the extra energy
 corresponding to the non-solenoidal perturbation
 has been removed.
Values of $|r_h^n-1|$ in 
 Figure \ref{fig:GePUP-ES_energyAndDivCurves}(c,d)
 are much larger than those in 
 Figure \ref{fig:viscousBox}(c,d), 
 which is still acceptable because
 the deviation of $r_h^n$ from 1 is about $10^{-5}$
 even for the case of the largest perburbation $\epsilon=1$. 
In Figure \ref{fig:GePUP-ES_energyAndDivCurves}(e,f), 
 the evolution curve of $\|\mathbf{DU}^n_{\epsilon}\|_{\mathcal{C}}$
 for each $\epsilon>0$
 overlaps with that of $\epsilon=0$
 except for the first four steps.

\begin{figure}
  \centering
  \subfigure[Re $=10^2$ and Cr = 0.1]{
    \includegraphics[width=0.475\textwidth]{
      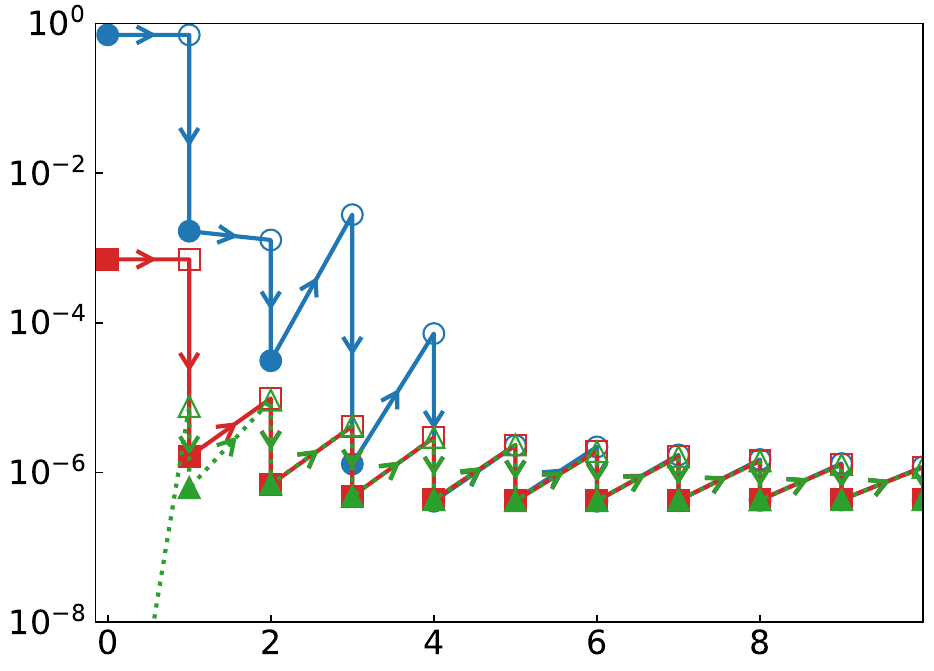}
  }
  \hfill
  \subfigure[Re $=10^{4}$ and Cr = 0.5]{
    \includegraphics[width=0.475\textwidth]{
      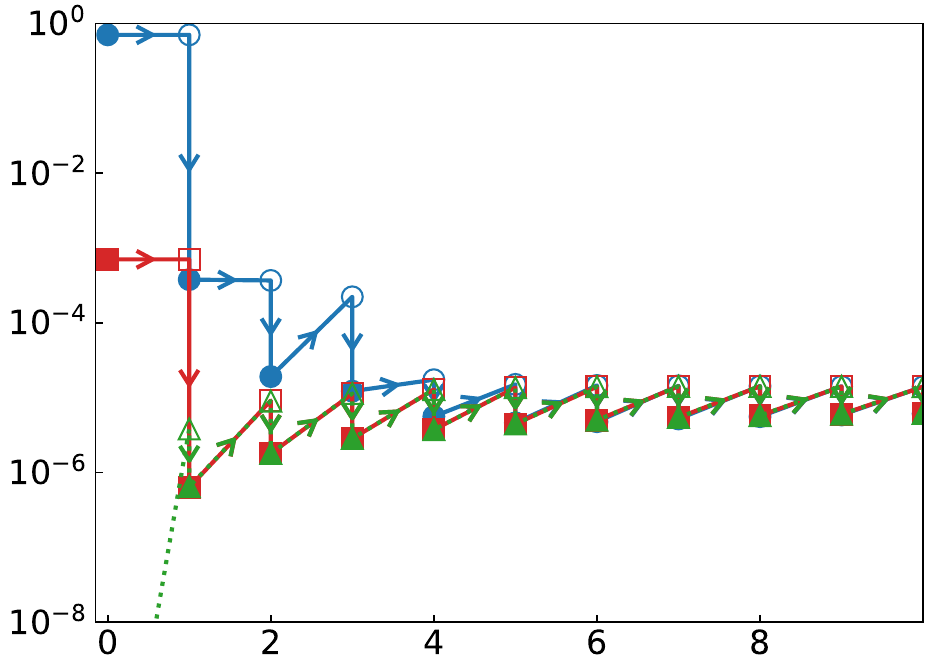}
  }

  \subfigure[Re $=10^2$ and Cr = 0.01]{
    \includegraphics[width=0.475\textwidth]{
      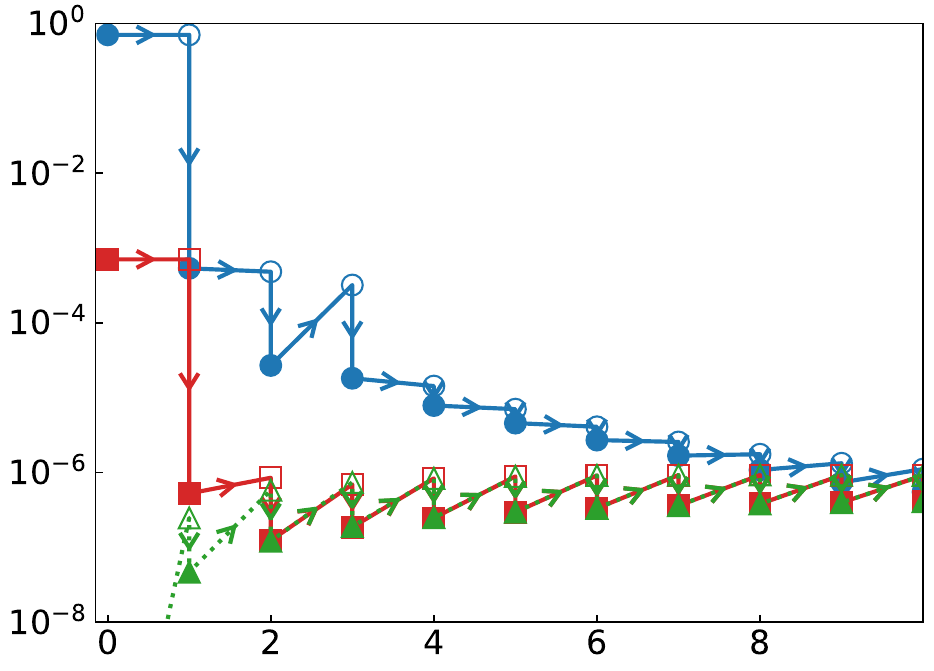}
  }
  \hfill
  \subfigure[Re $=10^{4}$ and Cr = 0.05]{
    \includegraphics[width=0.475\textwidth]{
      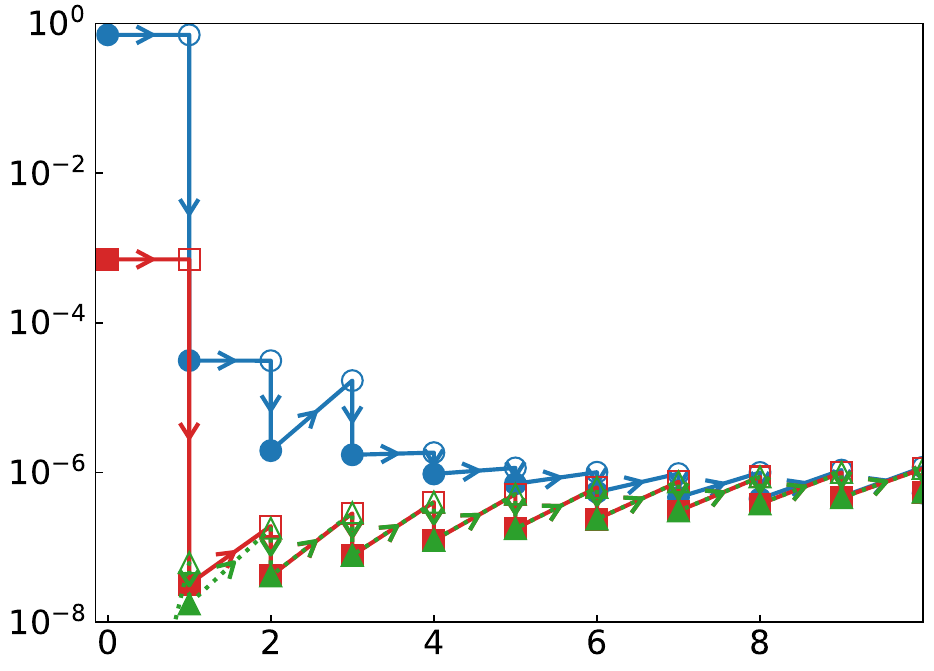}
  }
  \caption{Evolution of discrete velocity divergence
    produced by the GES scheme in Table \ref{tab:GePUP-SAV-SDIRK-policies}
    for solving the perturbed viscous box tests
    on a uniform grid of $h=\frac{1}{512}$.
    The abscissa of all subplots is the index of time steps.
    Subplots (a,b) can be considered as more detailed versions  
    of Figure \ref{fig:GePUP-ES_energyAndDivCurves} (e,f). 
    Values of $\|\mathbf{DU}^{n}_{\epsilon}\|_{\mathcal{C}}$
    for $\epsilon=0, 10^{-3}, 1$ 
    are represented respectively 
    by   solid markers ``$\blacktriangle$,''
    ``{\tiny $\blacksquare$},'' and 
    ``$\bullet$''
    while $\|\mathbf{DW}^{n}_{\epsilon}\|_{\mathcal{C}}$
    respectively by the corresponding
    the hollow markers
    ``{\tiny $\triangle$},'' ``{\tiny $\square$},'' and ``$\circ$.''
    Results of $\epsilon=10^{-1}, 10^{-2}$
    are in between those of $\epsilon=10^{-3}, 1$
    and have the same qualitative pattern. 
  }
  \label{fig:GePUP-ES_DivCurves}
\end{figure}

To examine what happens in these four steps,  
 the temporal variations of the $L_2$-norm of the velocity divergence 
 are plotted in Figure \ref{fig:GePUP-ES_DivCurves}
 for the two representative cases of $\epsilon=1,10^{-3}$.
Different from previous figures on velocity divergence,
 Figure \ref{fig:GePUP-ES_DivCurves}
 contains values of both $\|\mathbf{DW}^n_{\epsilon}\|_{\mathcal{C}}$
 and $\|\mathbf{DU}^n_{\epsilon}\|_{\mathcal{C}}$,
 which are connected by the discrete projection
 that approximates $\cProjLH$ in (\ref{eq:GePUPSAVRKFinal}). 
As illustrated in Figure \ref{fig:GePUP-ES_DivCurves}, 
 the discretization error
 might increase $\|\mathbf{DW}^{n}_{\epsilon}\|_{\mathcal{C}}$
 during a single time step, 
 but the magnitude of this increase
 approaches a constant value after about eight time steps;
 furthermore,
 the discrete projection at the end of a time step
 counteracts this increase.
In each subplot of Figure \ref{fig:GePUP-ES_DivCurves}, 
  the sequence of solid markers ``$\bullet$'' and ``{\tiny $\blacksquare$}''
 during the first four time steps
 clearly demonstrates the fast decay
 of $\|\mathbf{DU}^n_{\epsilon}\|_{\mathcal{C}}$. 
In addition,
 a comparison of the two rows of Figure \ref{fig:GePUP-ES_DivCurves}
 shows that 
 the adverse effects of the perturbation
 upon the divergence-free condition can be further reduced 
 by decreasing the size of the initial time steps. 

We have demonstrated two aspects of 
 the prominent advantage of GePUP-ES 
 in handling an initially non-solenoidal velocity.
First,
 the velocity divergence decays exponentially. 
Second,
 Figures \ref{fig:diffVelNormsEvolution}
 and \ref{fig:GePUP-ES_DivCurves} suggest that 
 the time period of adverse effects
 caused by the initially non-solenoidal velocity 
 can be very much shortened by reducing
 the size of the first several time steps.
We emphasize that theoretically there is no guarantee
 of the decay of velocity divergence
 for \emph{any} initially non-solenoidal velocity yet, 
 as the well-posedness of INSE
 is still one of the unsolved Millennium Prize Problems. 

\begin{figure}
  \centering
  \subfigure[$\mathcal{E}^n_h$]{
    \includegraphics[width=0.475\textwidth]{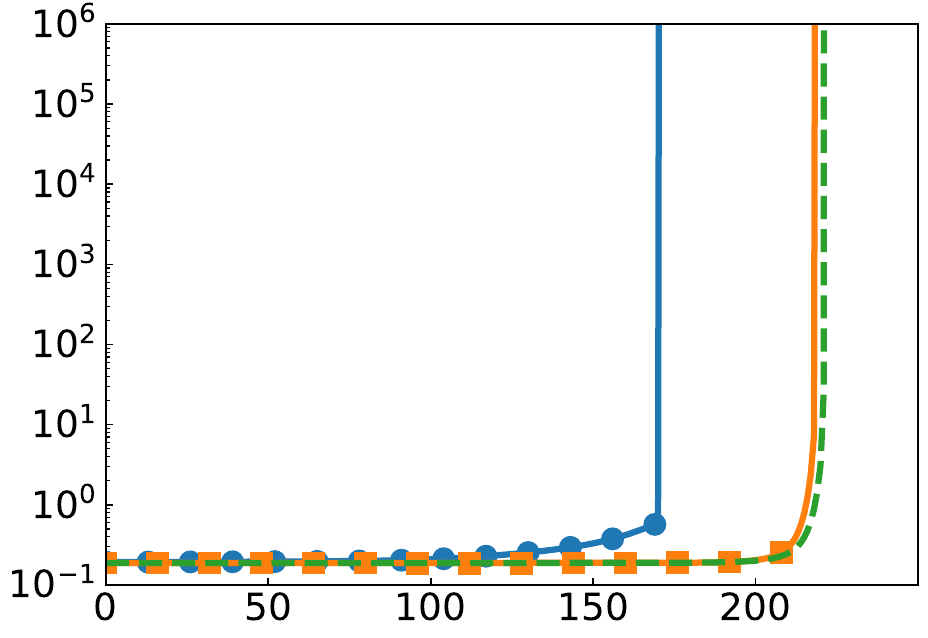}
  }
  \subfigure[$\|\mathbf{DU}^n_{\epsilon}\|_{\mathcal{C}}$]{
    \includegraphics[width=0.475\textwidth]{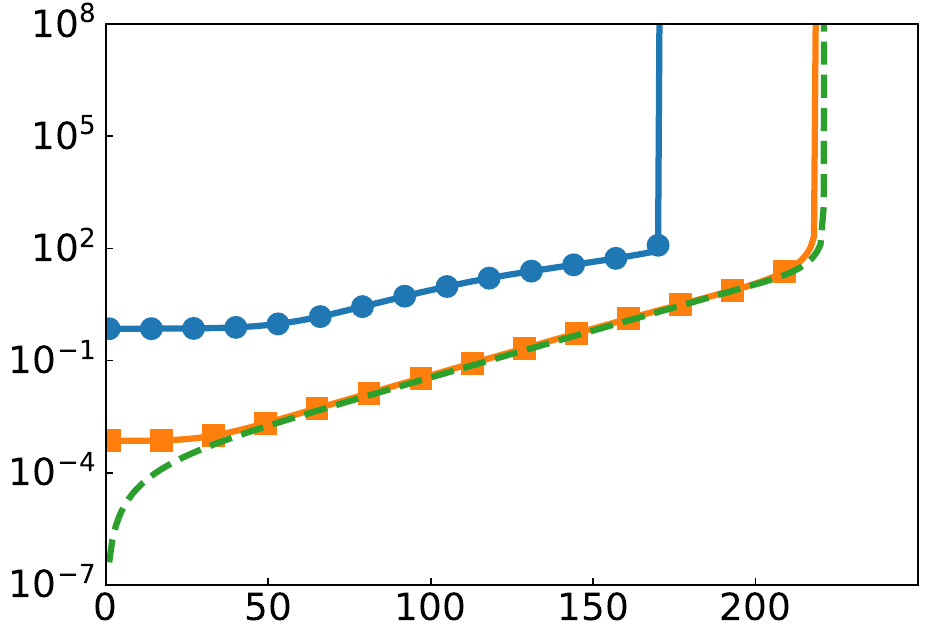}
  }
  \caption{Evolution of $\mathcal{E}^n_h$ and
    $ \|\mathbf{DU}^n_{\epsilon}\|_{\mathcal{C}}$
    produced by applying the MOL approach
    to the strong form of UPPE in (\ref{eq:UPPEstrong}) 
    for Re = $10^{4}$ , Cr = 0.5, and $h=\frac{1}{512}$.
    The curves marked by  ``$\bullet$'' and
    ``{\tiny $\blacksquare$}''
    correspond to $\epsilon=1$ and $ 10^{-3}$, respectively,
    while the dashed curve corresponds to $\epsilon=0$.
    The abscissa in all subplots is time.
  }
  \label{fig:uppe_energyAndDivCurves_Re1e4}
\end{figure}

The finite-volume based MOL scheme 
 is also applied to the strong form of UPPE in (\ref{eq:UPPEstrong})
 to obtain solutions
 of the perturbed viscous-box tests. 
As shown in Figure \ref{fig:uppe_energyAndDivCurves_Re1e4},
 both the velocity divergence and the total kinetic energy 
 blow up,
 even for the case of an initially solenoidal velocity.
It is also clear that the blow-up time becomes sooner
 as $\epsilon$ gets larger.

Results of viscous-box tests solved by GePUP,
 as shown in \cite[Table 3]{zhang16:_GePUP},
 are very close to those
 in Table \ref{tab:viscousBoxRe10000} for $\mathrm{Re}=10^{4}$.
In the case of $\mathrm{Re} = 10^{2}$, 
 the error norms of $\mathbf{u}$ produced by GePUP-ES
 on the finest grid
 are approximately 20\% smaller than those by GePUP, 
 cf. Table \ref{tab:viscousBoxRe100} 
 and \cite[Table 4]{zhang16:_GePUP}.
For the perturbed viscous box tests in this subsection, 
 results of GePUP also show fast decays 
 of discrete velocity divergence. 
However,
 GePUP-ES is a decisive advance from GePUP
 in that, in the semi-discrete case,
 the decay of velocity divergence and kinetic energy
 can be rigorously proven.

 


%% file: tab/singleVortexBox_Gamma1_Re2e4.tex
\begin{tabular}{ccccccc}
  \hline 
  \multicolumn{2}{c}{$h$} & $\frac{1}{256}-\frac{1}{512}$ & Rate & $\frac{1}{512}-\frac{1}{1024}$ & Rate & $\frac{1}{1024}-\frac{1}{2048}$
  \\ \hline \hline
  \multirow{3}*{$\mathbf{u}$} &  $L_{\infty}$ & 9.44e-04 & 3.63 & 7.61e-05 & 3.86 & 5.24e-06
  \\ & $L_1$ & 3.87e-05 & 3.76 & 2.86e-06 & 3.90 & 1.91e-07
  \\ & $L_2$ & 7.34e-05 & 3.73 & 5.52e-06 & 3.89 & 3.73e-07
  \\ \hline
  
  \multirow{3}*{$\nabla\cdot\mathbf{u}$} &$L_{\infty}$ & 1.08e-02 & 2.03 & 2.64e-03 & 2.82 & 3.75e-04
 \\ & $L_1$ & 5.48e-05 & 3.38 & 5.27e-06 & 3.98 & 3.35e-07
 \\ & $L_2$ & 3.26e-04 & 2.47 & 5.87e-05 & 3.28 & 6.04e-06
 \\ \hline

  \multirow{3}*{$q$} & $L_{\infty}$ & 8.96e-06 & 2.63 & 1.45e-06 & 2.94 & 1.89e-07
 \\ & $L_1$ & 7.91e-07 & 2.41 & 1.49e-07 & 2.50 & 2.63e-08
 \\ & $L_2$ & 1.38e-06 & 2.71 & 2.10e-07 & 2.59 & 3.50e-08
 \\ \hline

  \multirow{3}*{$\nabla q$} & $L_{\infty}$ & 4.60e-04 & 2.58 & 7.67e-05 & 3.17 & 8.51e-06
  \\ & $L_1$ & 1.78e-05 & 3.28 & 1.83e-06 & 2.89 & 2.48e-07
  \\ & $L_2$ & 3.75e-05 & 3.24 & 3.98e-06 & 3.00 & 4.96e-07
  \\ \hline


  
\end{tabular}


%% file: tab/viscousBox_Re1e4_relaxCoef1.tex
\begin{tabular}{ccccccc}
  \hline
  \multicolumn{2}{c}{$h$} 
  & $\frac{1}{128}-\frac{1}{256}$ & Rate
  & $\frac{1}{256}-\frac{1}{512}$ & Rate
  & $\frac{1}{512}-\frac{1}{1024}$
  \\ \hline \hline
  \multirow{3}*{$\mathbf{u}$}
  & $L_{\infty}$ & 1.02e-03 & 3.20 & 1.10e-04 & 4.09 & 6.49e-06
  \\ 
  & $L_1$ & 4.10e-05 & 3.53 & 3.55e-06 & 3.94 & 2.32e-07
  \\ 
  & $L_2$ & 9.85e-05 & 3.50 & 8.71e-06 & 3.96 & 5.59e-07
  \\ \hline

  \multirow{3}*{$\nabla\cdot\mathbf{u}$}
  & $L_{\infty}$ & 7.33e-03 & 1.32 & 2.94e-03 & 2.44 & 5.41e-04
  \\ 
  & $L_1$ & 3.15e-04 & 3.18 & 3.47e-05 & 3.70 & 2.68e-06
  \\ 
  & $L_2$ & 8.79e-04 & 2.22 & 1.89e-04 & 2.99 & 2.38e-05
  \\ \hline

  \multirow{3}*{$q$}
  & $L_{\infty}$ & 4.17e-04 & 1.75 & 1.24e-04 & 1.96 & 3.18e-05
  \\ 
  & $L_1$ & 6.61e-05 & 1.23 & 2.81e-05 & 1.70 & 8.69e-06
  \\ 
  & $L_2$ & 9.92e-05 & 1.39 & 3.78e-05 & 1.73 & 1.14e-05
  \\ \hline
  
  \multirow{3}*{$\nabla q$}
  & $L_{\infty}$ & 5.01e-03 & 1.70 & 1.54e-03 & 1.86 & 4.27e-04
  \\ 
  & $L_1$ & 6.82e-04 & 1.71 & 2.08e-04 & 1.86 & 5.75e-05
  \\ 
  & $L_2$ & 1.02e-03 & 1.69 & 3.15e-04 & 1.90 & 8.43e-05
  \\ \hline

  
\end{tabular}


%% file: tab/viscousBox_Re100_relaxCoef1.tex
\begin{tabular}{ccccccc}
  \hline 
  \multicolumn{2}{c}{$h$}
  & $\frac{1}{64}-\frac{1}{128}$  &  Rate
  & $\frac{1}{128}-\frac{1}{256}$ &  Rate
  & $\frac{1}{256}-\frac{1}{512}$
  \\ \hline \hline
  \multirow{3}*{$\mathbf{u}$}
  & $L_{\infty}$ & 7.84e-06 & 2.46 & 1.43e-06 & 2.78 & 2.08e-07
  \\ 
  & $L_1$ & 2.03e-06 & 3.99 & 1.28e-07 & 4.04 & 7.76e-09
  \\ 
  & $L_2$ & 2.63e-06 & 3.96 & 1.69e-07 & 3.99 & 1.07e-08
  \\ \hline

  \multirow{3}*{$\nabla\cdot\mathbf{u}$}
  & $L_{\infty}$ & 3.90e-04 & 1.91 & 1.03e-04 & 2.02 & 2.55e-05
  \\ 
  & $L_1$ & 1.55e-05 & 3.83 & 1.09e-06 & 3.75 & 8.11e-08
  \\ 
  & $L_2$ & 5.62e-05 & 3.24 & 5.96e-06 & 2.96 & 7.64e-07
  \\ \hline

  \multirow{3}*{$q$}
  & $L_{\infty}$ & 2.37e-04 & 1.83 & 6.67e-05 & 1.91 & 1.78e-05
  \\ 
  & $L_1$ & 3.33e-05 & 2.04 & 8.10e-06 & 2.05 & 1.95e-06
  \\ 
  & $L_2$ & 5.09e-05 & 2.07 & 1.21e-05 & 2.06 & 2.92e-06
  \\ \hline

  \multirow{3}*{$\nabla q$}
  & $L_{\infty}$ & 4.92e-03 & 0.93 & 2.57e-03 & 0.42 & 1.92e-03
  \\ 
  & $L_1$ & 3.24e-04 & 1.97 & 8.26e-05 & 2.01 & 2.05e-05
  \\ 
  & $L_2$ & 5.58e-04 & 1.72 & 1.69e-04 & 1.82 & 4.79e-05
  \\ \hline


\end{tabular}


%% file: sec/conclusions.tex
\section{Conclusions}
\label{sec:conclusions}

We have shown that the INSE with no-slip conditions
 can be equivalently reformulated
 as variants of the GePUP formulation \cite{zhang16:_GePUP},
 where the main evolutionary variable
 is a non-solenoidal velocity with electric boundary conditions
 whose divergence, controlled by a heat equation
 with homogeneous Dirichlet boundary conditions, 
 decays exponentially.
This GePUP-E reformulation is suitable
 for numerically solving the INSE because
\begin{itemize}
\item time integration and spatial discretization
  are completely decoupled so that
  high-order INSE solvers can be easily obtained
  from menu choices of orthogonal policies,
\item the constituting modules such as a time integrator
  are employed in a black-box manner
  so that no internal details of any module are needed
  in building the INSE solver, 
\item the influences of nonzero velocity divergence
  upon numerical stability and accuracy are clear,
\item a coupling of GePUP-E to SAV
  yields semi-discrete schemes
  with monotonically decreasing kinetic energy.
\end{itemize}

Results of numerical experiments confirm the analysis. 

We are currently augmenting the GePUP-ES solver
 to parallel computing and adaptive mesh refinement
 \cite{zhao2025FourthOrderAMR}
 for an enhanced resolution of
 the multiple time scales and length scales
 in flows at moderate or high Reynolds numbers. 
Another work in progress
 is the development of GePUP solvers
 for the INSE with irregular and moving boundaries
 via poised lattice generation \cite{zhang2024PLG}.

The next step along this line of research
 is the design of \emph{fully discrete} algorithms
 that ensure decays of velocity divergence and total kinetic energy.
We also plan to consider other types of boundary conditions
 such as the nonhomogeneous Dirichlet conditions,
 the radiation conditions,
 and mixed conditions.

